\documentclass[12pt]{article}
\usepackage{agzt_macros}
 
\usepackage{amsmath, amsfonts}
\usepackage{forest}
\usepackage{tabularx}
\usepackage{scalerel,graphicx}
 
\usepackage{mathabx, tikz-cd}
\usepackage{mathrsfs}
\usepackage{appendix}
\usetikzlibrary{decorations.pathmorphing,shapes}
\usetikzlibrary{arrows.meta}
\tikzcdset{arrow style=tikz,
	squigarrow/.style={
		decoration={
			snake, 
			amplitude=.4mm,
			segment length=2mm
		}, 
		rounded corners=.2pt,
		decorate
	}
}

\newcommand{\cotr}{\triangleleft}
\newcommand{\ZZnonneg}{\mathbb{Z}_{\geq0}}

\newcommand{\kvec}{\mathbf{k}}
\newcommand{\lvec}{\mathbf{l}}
\newcommand{\rvec}{\mathbf{r}}

\newcommand{\evec}{\mathbf{e}}
\newcommand{\btheta}{{\boldsymbol\vartheta}}
\newcommand{\bsigma}{{\boldsymbol\sigma}}
\newcommand{\btau}{{\boldsymbol\tau}}
\newcommand{\balpha}{{\boldsymbol\alpha}}
\newcommand{\bbeta}{{\boldsymbol\beta}}
\newcommand{\bmu}{{\boldsymbol\mu}}
\newcommand{\beeta}{{\boldsymbol\eta}}

\newcommand{\Alpha}{\mathcal{A}}

\newcommand{\mua}{\operatorname{m}}
\newcommand{\muu}{\mu}
\newcommand{\emptybrackets}{\ensuremath{[\cdot, \cdot]}}
\newcommand{\emptybrace}{\ensuremath{\{\cdot,\cdot\}}}

\DeclareMathOperator{\ad}{ad}
\DeclareMathOperator{\indec}{Indec}
\newcommand{\ind}[2]{\operatorname{ind}_{#1}(#2)}
\newcommand{\indmax}{\operatorname{ind}_{\max}}

 \newcommand{\abs}[1]{\lvert #1\rvert}

\newcommand{\eqlabel}[1]{\label{#1}}
\newcommand{\noproof}[1]{\!\!\!\!} 
\newcommand{\graphref}[1]{\!\!\!\!\!\!\!\!\!\!\phantom{\ref{#1}} } 
\newcommand{\nographref}[1]{\ref{#1}}    
\newcommand{\graphcaption}[1]{\unskip}                  

\title{On post-Lie structures for free Lie algebras}
\author{Annika Burmester, Ulf Kühn}
\date{ }

\begin{document}
	\maketitle
	\begin{abstract}
   We study post-Lie structures on free Lie algebras, the Grossman-Larson product on their enveloping algebras, and provide an abstract formula for its dual coproduct.  This might be of interest for the general theory of post-Hopf
algebras. Using a magmatic approach, we explore post-Lie algebras connected to multiple zeta values and their $q$-analogues. For multiple zeta values, this framework yields an algebraic interpretation of the Goncharov coproduct.  Assuming that the Bernoulli numbers satisfy the so called threshold shuffle identities, we present a post-Lie structure, whose induced Lie bracket we expect to restrict to the dual of indecomposables of multiple $q$-zeta values. Our post-Lie algebras align with Ecalle’s theory of bimoulds: we explicitly identify the ari bracket with a post-Lie structure on a free Lie algebra, and conjecture a correspondence for the uri bracket.

	\end{abstract}
	
	\tableofcontents
	\phantomsection
	
	\addcontentsline{toc}{section}{Introduction}
	\section*{Introduction}

A post-Lie algebra is a Lie algebra $(\mathfrak{g},\emptybrackets)$ equipped with an additional product $\tr$, which satisfies certain compatibility conditions with the Lie bracket $\emptybrackets$. This rather new notion was introduced by B.~Vallette \cite{Val07} in the context of operads of commutative trialgebras. 
Meanwhile, there are a multitude of applications in various areas of mathematics.

A main feature of any post-Lie algebra $(\mathfrak{g},\emptybrackets,\tr)$ is that it induces a second Lie structure $(\mathfrak{g},\emptybrace)$. 
In this article, we aim for an explicit description of the subadjacent Hopf algebra of the 
associated post-Hopf algebra $(\mathcal{U}(\mathfrak{g}),\tr)$, i.e. the
Hopf algebra structure $(\mathcal{U}(\mathfrak{g}),\glp,\co)$ on the universal enveloping algebra of $(\mathfrak{g},\emptybrace)$. Of particular interest for us is its graded dual $(\mathcal{U}(\mathfrak{g})^\vee,\shuffle,\Delta_{\glp})$.

To our knowledge, we give  
the first general description of the graded dual coproduct $\Delta_{\glp}$ to the Grossman-Larson product $\glp$. We apply these results to particular post-Lie structures arising in the theory of multiple zeta values and their $q$-analogues.  
In \cite[Prop.~3.3]{EF15}, a recursive description of the Grossman-Larson product $\glp$ on $\mathcal{U}(\mathfrak{g})$ is given. For the specific post-Lie algebras studied in this work, we present alternative descriptions that allow efficient explicit computations. 

We indicate how our work relates to multiple zeta values and their $q$-analogues. We will address these topics in forthcoming papers.

\textbf{Multiple zeta values.} The space of all $\Q$-linear combinations of multiple zeta values is an algebra $\Z$. 
Since there are no tools so far for treating transcendence issues of multiple zeta values, 
it is common to replace this subalgebra of the real numbers by some algebraic model.
We focus on the algebra of formal multiple zeta values $\Z^f$. This algebra is a quotient of the shuffle algebra $(\Q\langle x_0,x_1\rangle,\shuffle)$ by the ideal generated by the extended double shuffle relations.
A well-known conjecture states that there is an algebra isomorphism $\Z^f \cong \Z$. 

In a recent work \cite{AGZT}, we studied the post-Lie algebra structure related to the Ihara bracket $\emptybrace_I$ on $\Lie(x_0,x_1)$ and computed the Grossman-Larson product $\glp_I$ on $\Q\langle x_0,x_1\rangle$ with the construction in \cite{EF15} explicitly. We get the following two dual Hopf algebra structures
\[
\xymatrix{ 
(\Q\langle x_0,x_1\rangle,\shuffle,\Delta_G ) 
\ar@{<~>}[rrr]^{\text{graded dual}}&&&    
(\Q\langle x_0,x_1 \rangle,\glp_I,\co). }
\] 
The left hand side is a particular example of a more general construction by Goncharov \cite{Gon05}.

Racinet \cite{Ra} studied the formal multiple zeta values $\Z^f$ by introducing the subspace $\mathfrak{dm}_0\subset \Lie(x_0,x_1)$ corresponding to extended double shuffle relations. A main result of Racinet's thesis is that $(\mathfrak{dm}_0,\emptybrace_I)$ is a Lie subalgebra of $(\Lie(x_0,x_1),\emptybrace_I)$. Together with the above discussion, this implies that the coproduct $\Delta_G$  induces a Hopf algebra structure on the quotient of $\Z^f$ by $\zeta^f(2)$.

The Hopf algebra structure for motivic multiple zeta values, which define another algebraic model for multiple zeta values, is similar.  Brown used motivic multiple zeta values to obtain his celebrated theorem  
\cite{br} 
that the space of multiple zeta values is spanned by multiple zeta values $\zeta(k_1,\ldots, k_d)$, 
which have indices $k_i \in \{2,3\}$. 
Since central arguments in the proof of
rely on properties of $\Delta_G$ on $\Q\langle x_0,x_1\rangle$ only, we were able to
revisit the proof for the algebra of formal multiple zeta values \cite{AGZT} under some assumption on $\mathfrak{dm}_0$.

\textbf{$q$-analogues of multiple zeta values.} By a $q$-analogue of multiple zeta values, we mean a power series in $\Q[[q]]$ which degenerates to a multiple zeta value under the limit $q\to 1$. 
The algebra of multiple $q$-zeta values $\Z_q$ is spanned, among other models, by Schlesinger-Zudilin multiple $q$-zeta values.
These are given by 
\[\zeta_q^{\operatorname{SZ}}(s_1,\ldots,s_r)=\sum_{0<n_1<\cdots<n_r}\frac{q^{n_1s_1}}{(1-q^{n_1})^{s_1}}\cdots \frac{q^{n_rs_r}}{(1-q^{n_r})^{s_r}},\quad s_1,\ldots,s_r\in\ZZnonneg, s_r>0.\]
Numerical experiments seem to indicate that  $\Z_q$ modulo formal quasi-modular forms is equipped with a Hopf algebra structure  \cite{BK20}. 

 It follows from the work in \cite{BaBu_cMES}, \cite{Bu_balanced} that we have a surjective algebra morphism $(\Q\langle V\rangle, \shuffle)\to \Z_q$,   
where $V$ is the
infinite alphabet $V=\{ v_0,v_1,v_2,\ldots\}$. Conjecturally, the kernel of this map can be described as the invariance under a certain involution $\tau$. Thus,
as for the multiple zeta values, we replace $\Z_q$ by the algebra of formal multiple q-zeta values $\Z_q^f=(\QV, \shuffle )/(w-\tau(w)\mid w\in \QV)$. The involution $\tau$ is graded by the weight but filtered for the depth\footnote{An alternative description of $\Z_q^f$  can be found in \cite{fMES}.}. 
Inspired by Racinet's work \cite{Ra} on multiple zeta values, in \cite{BPhD} the subspace
\[
\mathfrak{bm}_0 \subset \Lie(V),
\]
corresponding to the relations in the algebra $\Z_q^f$ is introduced. We expect that $\mathfrak{bm}_0$ is a Lie algebra, this would induce a Hopf algebra structure on the quotient of $\Z_q^f$ by the formal quasi-modular forms.

The path to our promising candidate for  a Lie 
bracket $\emptybrace_u$ on $\Lie(V)$,  such that the space $\mathfrak{bm}_0$ becomes a Lie subalgebra is intricate.
It started  with an attempt in the language of Ecalle's bimoulds  \cite{K-montreal}, that was rebuilt in the spirit of \cite{Ra}  in the thesis of the first author \cite[Theorem 3.20]{BPhD}. Manchon observed that the latter approach 
was giving a Lie bracket $\emptybrace$ induced from a post-Lie algebra structure on the Lie algebra obtained from an associative, noncommutative free algebra. 
Our main motivation for this paper is the question     whether the expected Lie bracket $\emptybrace_u$ is also associated to a post-Lie structure  on $\Lie(V)$.

If the Bernoulli numbers satisfy the so-called threshold shuffle identities, then $\emptybrace_u$ is a Lie bracket on $\Lie(V)$ associated with the post-Lie algebra $\big( \Lie(V), \emptybrackets,\tr_u \big)$.   Unconditionally, we
show that the associated depth-graded bracket $\emptybrace_a$ is in fact a Lie bracket associated to the 
post-Lie algebra $\big( \Lie(V), \emptybrackets,\tr_a \big)$. 
Moreover, we derive an explicit and effective formula for the corresponding coproduct $\Delta_a$.  

In the forthcoming paper \cite{Bu_extended}, the first author shows that the space $\mathfrak{lq}\subset\Lie(V)$, which consists of elements invariant under the associated depth-graded involution of $\tau$, is a Lie subalgebra of $\big( \Lie(V), \emptybrace_a \big)$. It is $\mathfrak{lq}$  isomorphic to the Lie algebra of swap invariant, alternal bimoulds and in addition it also contains $\gr_D \mathfrak{bm}_0$. 
A conjectural description of $\gr_D \mathfrak{bm}_0$ is given in \cite{K-montreal}. Thus, if $\big(\mathfrak{bm}_0, \emptybrace_u \big)$ were a Lie algebra, then $\Delta_a$ would induce a Hopf algebra structure on the depth-graded algebra $\gr_D \Z_q^f$ modulo formal quasi-modular forms. 

From the explicit description of the coproduct $\Delta_a$ on $\Q \langle V \rangle$, we observe that it has properties very similar to that of $\Delta_G$ used in \cite{br}.  
Therefore we view the results and observations  of this paper as another step towards the conjecture that $\Z_q$ is spanned by  multiple $q$-zeta values $\zeta_q(k_1,\ldots, k_d)$ of some particular kind, which have indices $k_i \in \{1,2,3\}$ \cite[Conj. 3.2]{BK20}.

This paper is organized as follows.

Section \ref{sec:theory} covers the preliminary background on post-Lie algebras and Hopf algebra structures on their universal enveloping algebras. We focus on the coproduct $\Delta_\glp$ that is dual to the Grossman-Larson product $\glp$ induced by a post-Lie algebra. As a first step, we study  $\Delta_\glp$ in general in Theorem \ref{thm:dual_glp_fine}.

In Section \ref{sec:free_Lie} we restrict ourselves to post-Lie structures on free Lie algebras, which allows more explicit formulas for the coproduct.

{\bf Theorem \ref{thm:free_coprod2}} \textit{Let $\Alpha$ be an alphabet 
and $\big( \Lie(\Alpha), \emptybrackets, \tr \big)$ be a graded post-Lie algebra.  Then, we have for a word $A\in \Q\langle \Alpha\rangle$
\[\Delta_\glp(A)=\sum_{A=A_1\cdots A_n} (A_1\otimes\one)\shuffle_{\bullet} \cotr^{\operatorname{irr}}(A_2)\shuffle_{\bullet}\cdots \shuffle_{\bullet} \cotr^{\operatorname{irr}}(A_n).\]
Here, $\cotr^{\operatorname{irr}}$ is the dual map to the triangle map $\tr$ restricted to letters in the right factors \eqref{eq:reduced_cotr}. The product $\shuffle_{\bullet}$ on $\Q\langle \Alpha\rangle^{\otimes2}$ is the shuffle product on the left factor and concatenation on the right factor \eqref{eq:sh_bull_free}.
}

Moreover, we explain post-Lie structures on free Lie algebras from a magmatic perspective. Then in the following three sections \ref{sec:ihara}, \ref{sec:ari}, and \ref{sec:uri}, we apply the developed theory to the Ihara bracket, ari bracket, and uri bracket.

In Section \ref{sec:ihara}, we show that a natural generalization of the Ihara bracket with respect to a subset $V_0\subset V$ is induced by a post-Lie structure $(\Lie(V),\emptybrackets,\tr_I)$. Then, we give
an explicit formula for its Grossman-Larson product $\glp_I$
in Proposition \ref{prop:glp_ihara_grouplike-formula}, from which we derive a description of the dual coproduct $\Delta_I$. 

{\bf Theorem \ref{thm:coprod_ihara}} \textit{For letters $a_1,\ldots,a_n$ in the alphabet $V=\{v_0,v_1,v_2,\ldots\}$, we have 
\begin{align*}
&\Delta_I(a_1\cdots a_n)= \\
&\sum_{\substack{0\leq k \leq n \\ 0\leq i_1< j_1\leq i_2 < j_2 \leq \cdots \leq i_k< j_k \leq n} } \Big( a_1\cdots a_{i_1} \shuffle  
 \mathbb{I}\big( S(a_{i_1+1}\cdots a_{j_1-1})  
 \shuffle a_{j_1+1}\cdots a_{i_2}; a_{j_1} \big) \shuffle\cdots\\
&\hspace{10em} \shuffle  \mathbb{I}\big(S(a_{i_k+1}\cdots a_{j_k-1})\shuffle a_{j_k+1}\cdots a_n  ; a_{j_k} \big)\Big) \otimes a_{j_1}\cdots a_{j_k}.
\end{align*}
Here, $\mathbb{I}(-;v)$ is the identity for $v\notin V_0$ and maps to the constant term otherwise.
}

Finally, we relate this expression to Goncharov's formula \cite{Gon05} for the above coproduct $\Delta_G$ and obtain in this way a new approach without the use of formal iterated integrals. 

In Section \ref{sec:ari}, we study a post-Lie structure $(\Lie(V),\emptybrackets,\tr_a)$, which is motivated by Ecalle's ari bracket on bimoulds \cite{Ecalle11}. 

{\bf Theorem \ref{thm:dual_pair_ari}} \textit{
Let $V=\{v_0,v_1,v_2,\ldots\}$ be equipped with the grading $\operatorname{wt}(v_0)=1$ and $\operatorname{wt}(v_i)=i$ for $i>0$. 
The graded post-Lie algebra $\big(\Lie(V),\emptybrackets,\tr_a \big)$ from Theorem \ref{thm_ari_is_post_lie} induces a dual pair of graded Hopf algebras 
    \[
    \xymatrix{ 
     (\Q\langle V\rangle,\shuffle,\Delta_a ) 
   \ar@{<~>}[rrr]^{\text{graded dual}}&&&    
   (\Q\langle V \rangle,\glp_a,\co). 
    }
   \]
An explicit and effective formula for the Grossman-Larson product $\glp_a$ is given by Proposition \ref{prop:glp_ari_formula}, and for its dual coproduct $\Delta_a$ it is given by Theorem \ref{thm:coprod_ari}. }

We observe in Proposition \ref{prop:glpa_123} and Conjecture \ref{conj:glpa_123} that $\Delta_a$  has properties similar to that of $\Delta_G$.

In Section \ref{sec:uri}, we study a conjectural post-Lie structure which we expect to induce a Lie algebra structure on the space $\mathfrak{bm}_0$ corresponding to formal multiple $q$-zeta values.

{\bf Theorem \ref{thm:uri_is_post_lie}} \textit{The tuple $(\Lie(V),\emptybrackets,\tr_u)$, where $\tr_u$ is introduced in Definition \ref{def:uri_triangle}, is a post-Lie algebra under the assumption that the Bernoulli numbers satisfy the threshold shuffle identities from Definition \ref{def:threshold_shuffle}.
}

The associated depth-graded of $(\Lie(V),\emptybrackets,\tr_u)$ is exactly the post-Lie algebra $(\Lie(V),\emptybrackets,\tr_a)$ from the previous section. Unfortunately, it is much harder to deal with these additional terms of non-homogeneous depth. 
If the threshold shuffle identities holds, we get a pair of dual graded Hopf algebras
  \[
    \xymatrix{ 
     (\Q\langle V\rangle,\shuffle,\Delta_u ) 
   \ar@{<~>}[rrr]^{\text{graded dual}}&&&    
   (\Q\langle V \rangle,\glp_u,\co). 
    }
   \]
A closed formula for the Grossman-Larson product $\glp_u$, or for the dual coproduct $\Delta_u$ is not known yet. However, using the recursive constructions from \cite{EFM18} these maps 
can be computed numerically for small weights  \cite{Confurius:phd_thesis}.  As indicated in \cite[Remark 6.28]{BaBu_cMES} another approach to an explicit formula for the coproduct $\Delta_u$ might be via the serial expansion of combinatorial bi-multiple Eisenstein series combined with an
extension of the work of Bachmann-Tasaka \cite{BaTa}. Since we expect that $\Delta_u$ descends to a coproduct on the algebra $\Z_q$ of multiple $q$-zeta values modulo quasi-modular forms
an explicit formula is of great interest.

In Section \ref{sec:ecalle}, we relate the post-Lie structures studied in Section \ref{sec:ihara}, \ref{sec:ari}, and \ref{sec:uri} to Ecalle's theory of bimoulds \cite{Ecalle11}. First, the post-Lie algebra $(\Lie(V),\emptybrackets,\tr_a)$ from Section \ref{sec:ari} provides a new understanding of the ari bracket on alternal bimoulds in terms of post-Lie structures on free Lie algebras.  
We also indicate how the conjectural post-Lie structure $(\Lie(V),\emptybrackets,\tr_u)$ from Section \ref{sec:uri} relates to the 
uri bracket on alternil bimoulds. This section allows to compare the different approaches to a Lie structure for multiple $q$-zeta values in terms of bimoulds or in terms of post-Lie structures on free Lie algebras.

In Appendix \ref{app:comb_id}, we introduce threshold indicator functions on compositions to describe combinatorial identities among them. Those will be applied in Section \ref{sec:ari} and \ref{sec:uri}. The appendix is self contained and might be, independent of our use for post-Lie algebras, of general interest. That the Bernoulli numbers satisfy the threshold shuffle identities from Definition \ref{def:threshold_shuffle} had been checked in a lot of cases by computer.

\textbf{Acknowledgments.} We thank Dominique Manchon for his suggestion that our work is related to the theory of post-Lie algebras, Henrik Bachmann for his helpful comments and interest on various stages of this article, Niclas Confurius for providing us with examples and his computational point of view, and Hidekazu Furusho for suggestions improving Section \ref{sec:ecalle}.
This project was partially supported by the first author’s JSPS KAKENHI
Grant 24KF0150.

 \section{Dualizing the Grossman-Larson product}	
	\label{sec:theory}
	
	This section provides the general theory for post-Lie algebras with an emphasis on the implied  Hopf algebra structure on the universal enveloping algebras  given by the Grossman-Larson product. We provide  a first approach for the dual coproduct. In the following, we fix a field $k$ of characteristic $0$.
	
	\subsection{Lie algebras and universal enveloping algebras} \label{subs:liealg_univalg}
	
	Let $(\mathfrak{g},\emptybrackets)$ be a Lie algebra over $k$ and $(\mathcal{U}(\mathfrak{g}),\cdot)$ be its universal enveloping algebra. The product $\cdot$ is induced by the concatenation product on the tensor algebra of $\mathfrak{g}$, therefore we will also call $\cdot$ sometimes concatenation product.  Moreover, given $a_1,a_2 \in \mathcal{U}(\mathfrak{g})$ we write $a_1a_2$ instead of $a_1 \cdot a_2$. The universal enveloping algebra $(\mathcal{U}(\mathfrak{g}),\cdot)$ is equipped with a Hopf algebra structure, which we shortly recall. The coproduct $\Delta:\mathcal{U}(\mathfrak{g})\to\mathcal{U}(\mathfrak{g})\otimes\mathcal{U}(\mathfrak{g})$ is the $k$-linear, multiplicative map given by
	\begin{equation}\eqlabel{eq:shuffle_coproduct}
		\Delta(a) = a \otimes \one + \one \otimes a, \qquad a\in \mathfrak{g}.
	\end{equation}
	It is induced from the shuffle coproduct on the tensor algebra of $\mathfrak{g}$, so as before we will also call this product shuffle coproduct and sometimes denote it by $\Delta_\shuffle$. 
	\begin{Notation} \label{not:sweedler}
		We use Sweedler's notation for the coproduct $\Delta$ on $\mathcal{U}(\mathfrak{g})$:
		\begin{equation*}
			\Delta(A) =	A_{(1)} \otimes A_{(2)},
		\end{equation*}
		i.e. the summation symbol is omitted and the summation will always be assumed over all the occurring pairs 
		(see e.g. \cite[§~1.2]{Sweedler}). 
		We extend Sweedler's notation iteratively
		\begin{equation*}
			\Delta^{[n]}(A) = A_{(1)} \otimes \Delta^{[n-1]}(A_{(2)})  
			=A_{(1)}\otimes \dots \otimes A_{(n+1)} .
		\end{equation*}
		Precisely, $\Delta^{[n]}$ means the composition
		\[\Delta^{[n]}=(\Delta\otimes \underbrace{\operatorname{id}\otimes\cdots \otimes \operatorname{id}}_{\text{$n-1$ times}})\circ\cdots \circ (\Delta\otimes\operatorname{id})\circ \Delta.\]
		Since $\Delta$ is coassociative, this notion is well-defined and does not depend on the choice of the tensor product factor we apply the coproduct to.
	\end{Notation}
	\begin{Example}\label{exm:sweedler}
		For $A = a_1a_2\in\mathcal{U}(\mathfrak{g})$ we have
		\begin{align*}
			\Delta(a_1a_2) &= \Delta(a_1) \Delta(a_2) = (a_1\otimes\one + \one\otimes a_1) (a_2 \otimes \one + \one\otimes a_2)\\
			&= a_1a_2\otimes\one + a_1\otimes a_2 + a_2\otimes a_1 + \one\otimes a_1a_2 \\
			&= A_{(1)} \otimes A_{(2)}.
		\end{align*}
		The notation $A_{(1)} \otimes A_{(2)}$ thus corresponds to $4$ summands in this case.
		 
	\end{Example}
	The antipode $S:\mathcal{U}(\mathfrak{g})\to \mathcal{U}(\mathfrak{g})$ is the anti morphism defined by
	\begin{align} \label{eq:antipode}
		S(a)=-a,\qquad a\in \mathfrak{g}.
	\end{align}

From the definition of a Hopf algebra, we easily derive the following.
\begin{Lemma} \label{lem:antipode_rel} For $A\in \mathcal{U}(\mathfrak{g})\backslash k\one$, we have
\[S(A_{(1)})A_{(2)}=A_{(1)}S(A_{(2)})=0.\]    
\end{Lemma}

  \begin{Lemma} \label{lem:boxed_brackets_antipode} We have for $v, a_1,\ldots, a_n \in \mathfrak{g} $  and $A=a_1 \cdots a_n \in \mathcal{U}(\mathfrak{g})$
\begin{equation} \eqlabel{eq:boxed_brackets_antipode}
 [\dots[[v,a_1],a_2],\dots,a_n] = S(A_{(1)}) v A_{(2)}.
 \end{equation}
\end{Lemma}
\begin{proof} The proof follows easily by induction on $n$. 
\end{proof}	 	
	
	\subsection{Post-Lie algebras and the Grossman-Larson product}
	
	All results in this subsection are taken from \cite{EF15}.

	\begin{Definition}\label{def:post_lie_algebra} A \emph{post-Lie algebra} $\big(\mathfrak{g},\emptybrackets,\tr\big)$ is
		a Lie algebra $(\mathfrak{g},\emptybrackets)$ together with a $k$-linear map $\tr\colon \mathfrak{g}\otimes\mathfrak{g} \to \mathfrak{g}$ such that for all $x,y,z\in\mathfrak{g}$ 
		\begin{align}
			\eqlabel{eq:post-lie1}
			x\tr [y,z] &= [x\tr y, z] + [y,x\tr z]\\
			\eqlabel{eq:post-lie2}
			[x,y]\tr z &= x\tr(y\tr z) - (x\tr y)\tr z - y\tr(x\tr z) + (y\tr x)\tr z. 
		\end{align}
		The map $\tr\colon \mathfrak{g}\otimes\mathfrak{g}\to\mathfrak{g}$ is often referred to as \emph{post-Lie product}.
 	\end{Definition}

	There are many examples for post-Lie algebras in various fields of mathematics,  see for example the introduction of \cite{jaza_postlie} or of \cite{MKL13}. The main objective of this article is to present and study further examples that are motivated by 
	multiple zeta values and its $q$-analogs. First, we prepare the general algebraic framework, which we will apply in the following sections.
	
	\begin{Proposition}\label{prop:post-lie_relates_lie-algebras} 
		Given a post-Lie algebra $\big(\mathfrak{g},\emptybrackets,\tr\big) $, then 
		$(\mathfrak{g},\emptybrace)$ with
		\begin{equation*}
			\{ x,y\} \coloneqq x\tr y - y\tr x + [x,y], \qquad x,y\in\mathfrak{g},
		\end{equation*}
		is also a Lie algebra. We call this Lie bracket the \emph{post-Lie bracket}.
	\end{Proposition} 
	\graphref{def:post_lie_algebra}
	\noproof{Proposition}
 The post-Lie bracket $\emptybrace$ is also called Grossman-Larson bracket (see \cite{AKEFM22}) or composition Lie bracket (see \cite{jaza_postlie}), and the Lie algebra $(\mathfrak{g},\emptybrace)$ is also called sub-adjacent Lie algebra (see \cite{BGST24}, \cite{Li24})).
	
	\begin{Remark} \label{rem:post_lie_structures} Following \cite{Bu21}, let $\mathfrak{g}=(V, [x, y])$, $\mathfrak{n}=(V, \{x, y\})$ be two Lie algebras with the same underlying vector space $V$. A post-Lie
		algebra structure on the pair $(\mathfrak{g}, \mathfrak{n})$ is a $k$-linear product $\tr:V \otimes V\to V$ satisfying the identities\footnote{The second identity implies, that $\mathfrak{o}(z)= \{ g \in \mathfrak{g} \,| \, g \tr z= 0 \} \subset \mathfrak{g}$ is a Lie subalgebra with respect to $ \emptybrace$. }
		\begin{align*}
			x \tr y - y \tr x &= \{x, y\} - [x, y],\\
			\{x, y\} \tr z &= x \tr (y \tr z) - y \tr (x \tr z),\\
			x \tr [y, z] &= [x \tr y, z] + [y, x \tr z],
		\end{align*}
		for all $x, y, z \in V$.
		It is easy to check that, if $\mathfrak{g}= (V,\emptybrackets,\tr)$  is a post-Lie algebra with
		associated Lie algebra $\mathfrak{n}= (V,\emptybrace)$, then $\tr$ defines a post-Lie structure on the pair $(\mathfrak{g},\mathfrak{n})$.
	\end{Remark}

	\begin{Remark} \label{rem:deriv_on_CA}
		In Definition \ref{def:post_lie_algebra} the condition \eqref{eq:post-lie1} is equivalent to the requirement that $a\, \tr$ acts like  a derivation on $ \mathfrak{g}$. The assignment  $a \mapsto D_a$, where $D_a(b) =  a\tr b$ give rise to a
		homomorphism of Lie algebras $\big(\mathfrak{g}, \emptybrace\big) \to \operatorname{Der}(\mathfrak{g})$, since by \eqref{eq:post-lie2}
		\[  D_{\{a,b\} } =
		D_{[a,b]}+D_{a\tr b} - D_{b\tr a} =  D_a D_b -D_b D_a =[D_a,D_b]. 
		\]    
	\end{Remark}
	
	\begin{Definition}\label{def:ext_post-lie} \graphref{def:post_lie_algebra}
	 	Let $(\mathfrak{g},\emptybrackets,\tr)$ be a post-Lie algebra. We extend the post-Lie product to a $k$-linear map $\tr\colon \mathcal{U}(\mathfrak{g}) \otimes \mathcal{U}(\mathfrak{g}) \to \mathcal{U}(\mathfrak{g})$ via
		\begin{align}
			\eqlabel{eq:ext1_one}
			x \tr \one  &= 0 \\
			\eqlabel{eq:ext1_trivial}
			\one\tr A &= A \\
			\eqlabel{eq:ext2_asso_lie}
			xA\tr y &= x\tr(A\tr y)  - (x\tr A)\tr y \\
			\eqlabel{eq:ext3_asso_asso}
			A\tr BC &= \big(A_{(1)}\tr B\big)\big(A_{(2)}\tr C\big)
		\end{align}
		for all $A,B,C\in\mathcal{U}(\mathfrak{g})$ and $x,y\in\mathfrak{g}$.
	\end{Definition}
	
	Note that \eqref{eq:ext1_one} and \eqref{eq:ext3_asso_asso} together imply that we also have for any non-trivial product $A\in \mathcal{U}(\mathfrak{g})$ that $A\tr \one=0$. Definition \ref{def:ext_post-lie} yields a unique and well-defined product on $\mathcal{U}(\mathfrak{g})$ by \cite[Prop.~3.1]{EF15}.
	
	\begin{Remark}\label{rem:post-lie-def}
		Observe that \eqref{eq:post-lie2} is by \eqref{eq:ext2_asso_lie} in Definition \ref{def:ext_post-lie} equivalent to the identity 
		\begin{equation*}
			[x,y] \tr z = xy\tr z - yx\tr z.
		\end{equation*}
 
		 \end{Remark}

	We extend now the maps $D_a$ that were given in Remark \ref{rem:deriv_on_CA}, i.e., we set for $a\in \mathfrak{g}$ and $B\in\mathcal{U}(\mathfrak{g})$ that $D_a(B)=a\tr B$.

	\begin{Lemma}\label{lem:ext_post-lie_deriv_formula}\graphref{def:post_lie_algebra}\graphref{rem:deriv_on_CA}
		We have a $k$-linear map
		\[\mathfrak{g}  \to \operatorname{Der}\big(\mathcal{U}(\mathfrak{g})\big),\quad a \mapsto D_a.\]
	\end{Lemma}
	
	\begin{proof} We show that the images $D_a$ are indeed derivations on $\mathcal{U}(\mathfrak{g})$. We have
		\begin{equation*}
			\Delta^{[n-1]}(x) = \sum_{i=1}^n \underbrace{\one\otimes \dots \otimes \one}_{\text{$i-1$ times}} \otimes\ x \otimes \underbrace{\one \otimes \dots \otimes \one}_{\text{$n-i$ times}}.
		\end{equation*} 
		Thus, by applying \eqref{eq:ext1_trivial} and \eqref{eq:ext3_asso_asso} $n-1$ times, we deduce for $x,t_1,\dots,t_n\in\mathfrak{g}$ 
		\begin{equation} \eqlabel{eq:deriv_on U}
			x\tr(t_1\cdots t_n) = \sum_{i=1}^n t_1\cdots t_{i-1}(x\tr t_i) t_{i+1}\cdots t_n.
		\end{equation}
	\end{proof}

\begin{Lemma}\label{lem:ext_post-lie_action}
		The extended post-Lie product $\tr:\mathcal{U}(\mathfrak{g}) \otimes \mathcal{U}(\mathfrak{g})\to \mathcal{U}(\mathfrak{g})$ given in Definition 
		\ref{def:ext_post-lie} restricts to a linear map
		\[
		\tr\,:\,\,\mathcal{U}(\mathfrak{g}) \otimes \mathfrak{g}  \to \mathfrak{g}.
		\]
\end{Lemma}
\begin{proof} It suffices to show $A\tr b\in \mathfrak{g}$ for $A = a_1\cdots a_n\in\mathcal{U}(\mathfrak{g})$ and $a_1,\ldots,a_n,b\in\mathfrak{g}$. We prove the claim by induction on $n$. For $A=\one $ and $A\in \mathfrak{g}$ the statement is obvious. Thus, let $n \ge 2$. We deduce from 
		\eqref{eq:ext2_asso_lie} that
		\begin{equation*}
			A\tr b = a_1\tr(a_2\cdots a_n\tr b) - (a_1\tr a_2\cdots a_n) \tr b.
		\end{equation*}
		
		We have $a_2\cdots a_n \tr b \in\mathfrak{g}$ by the induction hypothesis. For the second summand $a_1\tr a_2\cdots a_n$ we use \eqref{eq:deriv_on U} and get 
		\begin{equation*}
			(a_1\tr a_2\cdots a_n) = \sum_{i=2}^n a_2\cdots a_{i-1}(a_1\tr a_i) a_{i+1}\cdots a_n.
		\end{equation*}
		Thus, $a_1\tr a_2\cdots a_n$ is a sum of products of $n-1$ Lie elements. Thus, we get by the induction hypothesis $(a_1\tr a_2\cdots a_n)\tr b\in \mathfrak{g}$ for each $b\in \mathfrak{g}$.
	\end{proof}

Next, we present formulas to explicitly compute the extended post-Lie product from Definition \ref{def:ext_post-lie}.
	
	\begin{Remark}\label{rem:post_lie_recursion}
		Let  
		$A,B=b_1\cdots b_m\in \mathcal{U}(\mathfrak{g})$ with $b_1,\dots,b_m\in\mathfrak{g}$. By applying \eqref{eq:ext3_asso_asso} iteratively $m-1$ times and using Notation \ref{not:sweedler}, we get
		\begin{equation*}
			A\tr B = \big(A_{(1)}\tr b_1\big) \big(A_{(2)}\tr b_2\big) \cdots \big(A_{(m)}\tr b_m\big).
		\end{equation*}
		Therefore, it is central to understand  products of  type $\mathcal{U}(\mathfrak{g}) \tr \mathfrak{g} $  in detail.  
	 		For example, if $a_1,a_2,b\in\mathfrak{g}$, then applying \eqref{eq:ext2_asso_lie} gives
			\begin{align*}
				a_1a_2\tr b &= a_1\tr (a_2\tr b) - (a_1\tr a_2) \tr b.
			\end{align*}
  	 The proof of Lemma \ref{lem:ext_post-lie_action} gives a recursive algorithm to calculate the product   
   $a_1\cdots a_n \tr b$. It is given by the sum of $n!$ monomials, each given 
			by $n$  factors, in the non-associative magma $( \mathfrak{g},\tr)$. 
   In our applications we found more efficient ways to calculate such quantities, see Propositions~\nographref{prop:post_lie_on_letter},  \nographref{prop:mult_tr_ari} and  for a partial results we refer to   \nographref{lem:uri_magma_boxed}
   and Remark \nographref{rem:uri_tr_not_efective}.  In \cite[Proposition 3.7.]{jaza_postlie} another such formula is given in in the context of derivations and regularity structures. 
 	\end{Remark}

\begin{Theorem}[{\cite[Prop.~3.3]{EF15}},{\cite[Thm.~3.4]{EF15}}] \label{thm:glp_hopf_algebra}
Let $(\mathfrak{g}, \emptybrackets, \tr)$ be a post-Lie algebra. 
\begin{enumerate}
	\item The \emph{Grossman-Larson product} on $\mathcal{U}(\mathfrak{g})$ given by
			\begin{equation}\eqlabel{eq:glp_def}
				A\glp B = A_{(1)} (A_{(2)} \tr B) \qquad \text{for }A,B\in\mathcal{U}(\mathfrak{g})
			\end{equation}
			is associative. 
	\item The triple $\big(\mathcal{U}(\mathfrak{g}),\glp,\Delta\big)$ is a Hopf algebra.
	\item Let $\bar\g = (\mathfrak{g},\emptybrace)$ be the Lie algebra from Proposition \ref{prop:post-lie_relates_lie-algebras}. Then, the map
		\begin{align*}
			(\mathcal{U}(\bar\g),\cdot,\Delta) & \to (\mathcal{U}(\g),\glp,\Delta) \\
			a_1 a_2 \cdots a_n  &\mapsto  a_1 \glp a_2 \glp \ldots \glp a_n   \qquad (a_i\in\mathfrak{g}),  
		\end{align*}
where $\glp$ is the Grossman-Larson product defined by \eqref{eq:glp_def}, is an isomorphism of Hopf algebras. 
	\end{enumerate}
		 	\end{Theorem}
			\noproof{Theorem}
	
For further details of the proof we refer to \cite{EFM18} and \cite[Thm.~2.9]{jaza_postlie}.	 
 
\begin{Remark}
The notation post-Hopf algebra for the tuple $\big(\mathcal{U}(\mathfrak{g}),\Delta\big)$ consisting of the Hopf algebra $\mathcal{U}(\mathfrak{g})$ together with the map $\tr$
has become 
standard in recent preprints  
\cite{lishta22}, \cite{catoire2024}. In those articles they refer to as subadjacent Hopf algebra for
$\big(\mathcal{U}(\mathfrak{g}),\glp,\Delta\big)$  
\end{Remark}

	\begin{Remark}\label{rem:glp_post-lie_bracket_algorithm}  
		(i)	Let $A,B=b_1\cdots b_m\in\mathcal{U}(\mathfrak{g})$. Since 
\begin{equation}\eqlabel{eq:post_lie_recursion_glp}
			A\glp B = A_{(1)} (A_{(2)} \tr b_1) (A_{(3)}\tr b_2) \cdots (A_{(m+1)}\tr b_m),
		\end{equation}
we can calculate $A\glp B$ by using the approach described in Remark \ref{rem:post_lie_recursion}. 		
	
	(ii)	Let $a,b\in\mathfrak{g}\subset \mathcal{U}(\mathfrak{g})$. Since $\Delta(a) = a\otimes \one + \one\otimes a$, we get
		\begin{align*}
			a\glp b &= ab + a\tr b,  
		\end{align*} 
		and thus		\begin{equation}\eqlabel{eq:check_lie_glp}
			\{a,b\} = a\glp b - b\glp a.
		\end{equation}
  (iii) Assume that $G$ is a grouplike element in the completion $\widehat{\mathcal{U}(\mathfrak{g})}$, i.e., we have $\Delta(G)=G\otimes G$. Then, the Grossman-Larson product simplifies for any $B\in \mathcal{U}(\mathfrak{g})$ to
  \begin{equation}
      G\glp B=G(G\tr B).
  \end{equation}
	\end{Remark}

	\subsection{Coproducts dual to Grossman-Larson products} 
	\label{subsec:duality}
	
	\begin{Definition}\label{def:duality_pairing} Let $V=\bigoplus_{n\geq0} V_n$ and $W=\bigoplus_{m\geq0} W_m$ be graded $k$-vector spaces. In the following, we will always assume that the homogeneous subspaces are finite dimensional. We call a non-degenerate pairing $(-\mid -):V\otimes W\to k$, such that
	\[ (V_n\mid W_m)=0 \quad \text{ for all } n\neq m,\]
	a \emph{graded duality pairing}. If such a pairing exists, we also say that $V$ and $W$ are \emph{graded dual}.
	\end{Definition}

	Assume that $(\mathfrak{g},\emptybrackets,\tr)$ is a graded post-Lie algebra, i.e., we have $\mathfrak{g}=\bigoplus_{n\geq0} \mathfrak{g}_n$ and $[\mathfrak{g}_n,\mathfrak{g}_m],\ \mathfrak{g}_n\tr \mathfrak{g}_m\subset \mathfrak{g}_{n+m}$ for all $m,n\geq0$. The grading on $\mathfrak{g}$ induces a grading on the universal enveloping algebras $\mathcal{U}(\mathfrak{g})$ and $\mathcal{U}(\overline{\mathfrak{g}})$ as Hopf algebras, where the homogeneous subspaces are also finite dimensional.
	
	Now, let $\mathcal{U}(\mathfrak{g})^\vee$ be a Hopf algebra, which is graded dual to $\mathcal{U}(\mathfrak{g})$ with respect to some pairing $(-\mid -)$. By abuse of notation and as a preparation for the following parts on free Lie algebras, we denote the dual of the product $\cdot$ and the coproduct $\Delta$ on $\mathcal{U}(\mathfrak{g})$ under the pairing $(-\mid -)$ by $\dec$ and $\shuffle$, respectively. We obtain the following diagram
	\begin{equation} \label{eq:diagram_big_picture_introduction}
		\begin{tikzcd}[baseline=(current  bounding  box.center)]
			\big(\mathcal{U}(\mathfrak{g})^\vee,\shuffle,\dec \big) \arrow[d,twoheadrightarrow,"\text{mod products}"'] && \big(\mathcal{U}(\mathfrak{g}),\cdot ,\Delta\big) \arrow[ll,leftrightsquigarrow,"\text{graded dual}" '] 
			\\ 
			\big(\indec(\mathcal{U}(\mathfrak{g})^\vee),\delta\big) && \big(\mathfrak{g},\emptybrackets \big).
			\arrow[ll,leftrightsquigarrow, "\text{graded dual}"'] \arrow[u, hookrightarrow]
		\end{tikzcd} \end{equation}
	Here, $\indec(\mathcal{U}(\mathfrak{g})^\vee)$ is the Lie coalgebra of the indecomposables elements of $\mathcal{U}(\mathfrak{g})^\vee$.
	
	\begin{Remark}
	If $\widehat{\mathfrak{g}}$ is a pro-nilpotent Lie algebra, the picture can be extended. Precisely, there is an exponential map from the completed Lie algebra $\widehat{\mathfrak{g}}$ into the grouplike elements $\Grp(\widehat{\mathcal{U}(\mathfrak{g})})$ of the completed Hopf algebra $(\widehat{\mathcal{U}(\mathfrak{g})},\cdot,\Delta)$. Moreover, the Hopf algebra $(\mathcal{U}(\mathfrak{g})^\vee,\shuffle,\dec)$ is the representing algebra of the affine group scheme $\Q\operatorname{-Alg}\to \operatorname{Set},\ R\mapsto \Grp(\widehat{\mathcal{U}(\mathfrak{g})\otimes R})$. Equivalently, we can view $(\mathcal{U}(\mathfrak{g})^\vee,\shuffle,\dec)$ as the Hopf algebra of functions on $\Grp(\widehat{\mathcal{U}(\mathfrak{g})})$. 
	\end{Remark}

As a dual version of Lemma \ref{lem:antipode_rel}, we have the following.

\begin{Lemma} \label{lem:dual_antipode_rel} Denote the antipode of $(\mathcal{U}(\mathfrak{g})^\vee,\shuffle,\dec)$ by $S^\vee$, and for $A^\vee\in \mathcal{U}(\mathfrak{g})^\vee$ write $\dec(A^\vee)=A_{(1)}^\vee \otimes A_{(2)}^\vee$. Then, we have
\[S^\vee(A_{(1)}^\vee)\shuffle A_{(2)}^\vee=A_{(1)}\shuffle S^\vee(A_{(2)}^\vee)=0, \qquad A^\vee\in \mathcal{U}(\mathfrak{g})^\vee_{>0}.\]
\end{Lemma}

	We analyse the change of the picture \eqref{eq:diagram_big_picture_introduction} by altering the Lie bracket with a post-Lie structure: 
	\begin{equation} \label{eq:diagram_big_picture_postlie}
		\begin{tikzcd}[baseline=(current  bounding  box.center)]
			\big(\mathcal{U}(\mathfrak{g})^\vee,\shuffle,\Delta_{\glp} \big) \arrow[d,twoheadrightarrow,"\text{mod products}"'] && \big(\mathcal{U}(\mathfrak{g}),\glp ,\Delta\big) \arrow[ll,leftrightsquigarrow,"\text{graded dual}" '] 
			\\ 
			\big(\indec(\mathcal{U}(\mathfrak{g})^\vee),\delta_{\circledast}\big) && \big(\mathfrak{g},\emptybrace \big).
			\arrow[ll,leftrightsquigarrow,"\text{graded dual}"'] \arrow[u, hookrightarrow]
		\end{tikzcd} \end{equation}
		
		We characterize the coproduct $\Delta_{\glp}$  dual to the Grossman-Larson product $\glp$ given by \eqref{eq:glp_def} with respect to a duality pairing $(-\mid -)$ as in Definition \ref{def:duality_pairing} in this  abstract setting as follows.
	
\begin{Definition}\label{def:glcoprod_unique}	\graphref{def:duality_pairing} \graphref{thm:glp_hopf_algebra}
The map $\Delta_{\glp}:\mathcal{U}(\mathfrak{g})^\vee\to \mathcal{U}(\mathfrak{g})^\vee\otimes\ \mathcal{U}(\mathfrak{g})^\vee$ is the unique map satisfying
	\[ (A\glp B\mid C^\vee)=(A\otimes B \mid \Delta_{\glp}(C^\vee)), \qquad \text{for all } A,B\in \mathcal{U}(\mathfrak{g}),\ C^\vee\in \mathcal{U}(\mathfrak{g})^\vee.\]
	Here, we set as usual $(A\otimes B\mid C_1^\vee\otimes C_2^\vee)=(A\mid C_1^\vee)(B\mid C_2^\vee)$ for all $A,B\in \mathcal{U}(\mathfrak{g})$, $C_1^\vee,C_2^\vee\in \mathcal{U}(\mathfrak{g})^\vee$.
\end{Definition}

\begin{Definition} \label{def:cotr} We define the cotriangle map $\cotr:\mathcal{U}(\mathfrak{g})^\vee\to \mathcal{U}(\mathfrak{g})^\vee\otimes\mathcal{U}(\mathfrak{g})^\vee$ to be the dual map to $\tr$ on $\mathcal{U}(\mathfrak{g})$ from Definition \ref{def:ext_post-lie},
i.e., it is 
\[(A\tr B\mid C^\vee)=(A\otimes B\mid \cotr(C^\vee)), \]
for all $A,B\in\mathcal{U}(\mathfrak{g}),\ C^\vee\in \mathcal{U}(\mathfrak{g})^\vee$.

	\end{Definition}

\begin{Remark}
If we fix a basis $\mathsf{UB}$ 
   of the universal enveloping algebra $\mathcal{U}(\mathfrak{g})$, then   
Definition \ref{def:glcoprod_unique} and \ref{def:cotr}  allow the alternative descriptions
\[\Delta_{\glp}(C^\vee)=\sum_{A,B\in \mathsf{UB}} (A\glp B\mid C^\vee)\ A^\vee\otimes B^\vee \]
and
  \[
\cotr(C^\vee) = \sum_{A ,B  \in \mathsf{UB} } 
( A \tr B \mid C^\vee )\, A^\vee \otimes B^\vee.
   \]
\end{Remark}  

Next, let $\shuffle_\bullet$ be the product on $\mathcal{U}(\mathfrak{g})^\vee\otimes\mathcal{U}(\mathfrak{g})^\vee$ given by
	\begin{align} \label{eq:shuffle_bullet}
		(A_1^\vee\otimes A_2^\vee)\shuffle_\bullet (B_1^\vee\otimes B_2^\vee)= (A_1^\vee\shuffle B_1^\vee) \otimes (A_2 B_2)^\vee
	\end{align}
	for all $A_1^\vee,A_2^\vee,B_1^\vee,B_2^\vee\in\mathcal{U}(\mathfrak{g})^\vee$.

We also use the Sweedler notation for the deconcatenation coproduct on $\mathcal{U}(\mathfrak{g})^\vee$,
   \[\dec(C^\vee)=C_{(1)}^\vee \otimes C_{(2)}^\vee, \qquad C^\vee\in \mathcal{U}(\mathfrak{g})^\vee,\]
and also the iterated version as in Notation \ref{not:sweedler}.
 
	\begin{Proposition} \label{prop:dualGLP_gen} The coproduct given in Definition \ref{def:glcoprod_unique} satisfies
		\begin{align*}
			\Delta_{\glp}(C^\vee)= (C_{(1)}^\vee\otimes \one)\shuffle_\bullet \cotr(C_{(2)}^\vee).
		\end{align*}    
	\end{Proposition}
    
	\begin{proof} Fix a basis $\mathsf{UB}$ 
   of the universal enveloping algebra $\mathcal{U}(\mathfrak{g})$.
   Then, by means of Definition \ref{def:cotr} we compute for $C^\vee\in \mathcal{U}(\mathfrak{g})^\vee$
	\begin{align*}
	\Delta_{\glp}(C^\vee)&=\sum_{A,B\in \mathsf{UB}} (A\glp B\mid C^\vee)\ A^\vee\otimes B^\vee =\sum_{A,B\in \mathsf{UB}}  \big(A_{(1)}(A_{(2)}\tr B)\mid C^\vee\big)\ A^\vee\otimes B^\vee \\
    &=\sum_{A,B\in \mathsf{UB}}  \big(A_{(1)}\otimes (A_{(2)}\tr B)\mid \dec (C^\vee) \big)\ A^\vee\otimes B^\vee \\
	&=\sum_{A,B\in \mathsf{UB}} (A_{(1)}\mid C_{(1)}^\vee)(A_{(2)}\tr B\mid C_{(2)}^\vee)\ A^\vee\otimes B^\vee \\
	&=\sum_{A,B\in \mathsf{UB}} (A_{(1)}\mid C_{(1)}^\vee)(A_{(2)}\otimes B\mid \cotr (C_{(2)}^\vee))\ A^\vee \otimes B^\vee. \end{align*}
	Since $\shuffle$ is by definition the dual coproduct to $\Delta(A)=A_{(1)}\otimes A_{(2)}$, we deduce
	\begin{align*} 		\Delta_{\glp}(C^\vee)&=\sum_{A_1,A_2,B\in \mathsf{UB}} (A_1\mid C_{(1)}^\vee)(A_2\otimes B\mid \cotr (C_{(2)}^\vee))\ (A_1^\vee\shuffle A_2^\vee)\otimes B^\vee \\
	&=\sum_{A_1,A_2,B\in \mathsf{UB}} (A_1\mid C_{(1)}^\vee)(A_2\otimes B\mid \cotr (C_{(2)}^\vee))\ (A_1^\vee\otimes 1) \shuffle_\bullet (A_2^\vee\otimes B^\vee)\\
	&=(C_{(1)}^\vee\otimes 1) \shuffle_\bullet \cotr(C_{(2)}^\vee). \qedhere
	\end{align*}
	\end{proof}

If we assume that there exist unique product decompositions in $\mathcal{U}(\mathfrak{g})$, we can give a refinement of Proposition \ref{prop:dualGLP_gen}. 
    
\begin{Notation}\label{not:CU_has_UDB} We say $\mathcal{U}(\mathfrak{g})$ has a basis with unique decomposition, if the following holds.    There is a basis $\mathsf{UB}$ of $\mathcal{U}(\mathfrak{g})$ and a generating set $\mathsf{GS}$ of the Lie algebra $\mathfrak{g}$, such that any element $B\in\mathsf{UB}$ can be uniquely written as a product of elements in $\mathsf{GS}$:
    \[ B=b_1\cdots b_n,\qquad b_i\in \mathsf{GS}.\]
As before, we denote for $B\in \mathsf{UB}$ the dual basis element in $\mathcal{U}(\mathfrak{g})^\vee$ by $B^\vee\in \mathsf{UB}^\vee$. Then, the above is equivalent to requiring that for each $B^\vee\in \mathsf{UB}^\vee$ there is a unique $n\geq0$ such that $\dec^{[n-1]}(B^\vee)$ contains a summand $b_1^\vee\otimes\cdots\otimes b_n^\vee$ with $b_1,\ldots,b_n\in \mathsf{GS}$.
\end{Notation}	
	 Observe that we must have $\mathsf{GS}\subset \mathsf{UB}$.

    \begin{Definition} \label{def:cotr_reduced} Assume $\mathcal{U}(\mathfrak{g})$ has a basis with unique decomposition as in Notation \ref{not:CU_has_UDB}. 
    We define the reduced triangle map for $C^\vee\in \mathcal{U}(\mathfrak{g})^\vee$ by
    \[\cotr^{\operatorname{irr}}(C^\vee)=\sum_{\substack{A \in \mathsf{UB} \\ b\in \mathsf{GS}}} \big(A\otimes b\mid \cotr(C^\vee)\big) A^\vee\otimes b^\vee.\]
    \end{Definition}

	Note that we get directly $\cotr^{\operatorname{irr}}(\one)=0$.

    \begin{Theorem} \label{thm:dual_glp_fine} Assume $\mathcal{U}(\mathfrak{g})$ has a basis with unique decomposition as in Notation \ref{not:CU_has_UDB}, then for any $C^\vee\in\mathcal{U}(\mathfrak{g})^\vee$, we have
    \[\Delta_{\glp}(C^\vee)=\sum_{n\geq1} (C_{(1)}^\vee\otimes\one)\shuffle_\bullet \cotr^{\operatorname{irr}}(C_{(2)}^\vee)\shuffle_\bullet \cdots \shuffle_\bullet \cotr^{\operatorname{irr}}(C_{(n)}^\vee). \]
    \end{Theorem}

    \begin{proof} We obtain  
    \begin{align*}
    \Delta_{\glp}(C^\vee)&=\sum_{A,B\in\mathsf{UB}} (A\glp B\mid C^\vee)\ A^\vee\otimes B^\vee \\
    &= \sum_{n\geq1} \sum_{\substack{A\in\mathsf{UB} \\ b_2,\ldots,b_n\in\mathsf{GS}}} \big(A_{(1)}(A_{(2)}\tr b_2)\cdots (A_{(n)}\tr b_n)\mid C^\vee\big)\ A^\vee\otimes (b_2 \cdots b_n)^\vee,
    \end{align*}
    since any $B\in\mathsf{UB}$ is required to have a unqiue expression $B=b_2\cdots b_n$ with $b_i\in\mathsf{GS}$ and due to formula \eqref{eq:post_lie_recursion_glp}. Thus, we get
    \begin{align*}
    &\Delta_{\glp}(C^\vee) = \sum_{n\geq1} \sum_{\substack{A\in\mathsf{UB} \\ b_2,\ldots,b_n\in\mathsf{GS}}} (A_{(1)}\mid C_{(1)}^\vee)(A_{(2)}\tr b_2\mid C_{(2)}^\vee)\cdots (A_{(n)}\tr b_n\mid C_{(n)}^\vee)\\
    &\hspace{12cm} \cdot A^\vee\otimes (b_2\cdots b_n)^\vee \\
    &=\sum_{n\geq1} \sum_{\substack{A\in\mathsf{UB} \\ b_2,\ldots,b_n\in\mathsf{GS}}} (A_{(1)}\mid C_{(1)}^\vee)(A_{(2)}\otimes b_2\mid \cotr(C_{(2)}^\vee))\cdots (A_{(n)}\otimes b_n\mid \cotr(C_{(n)}^\vee))\\
    &\hspace{12cm} \cdot A^\vee\otimes (b_2\cdots b_n)^\vee \\
    &=\sum_{n\geq1} \sum_{\substack{A_1,\ldots,A_n\in\mathsf{UB} \\ b_2,\ldots,b_n\in\mathsf{GS}}} (A_1\mid C_{(1)}^\vee)(A_2\otimes b_2\mid \cotr (C_{(2)}^\vee))\cdots (A_n\otimes b_n\mid \cotr(C_{(n)}^\vee))\\
    &\hspace{8,5cm} \cdot (A_1^\vee\shuffle \cdots \shuffle A_n^\vee)\otimes (b_2\cdots b_n)^\vee \\
    &=\sum_{n\geq1} \sum_{\substack{A_2,\ldots,A_n\in\mathsf{UB} \\ b_2,\ldots,b_n\in\mathsf{GS}}} (A_2\otimes b_2\mid \cotr (C_{(2)}^\vee))\cdots (A_n\otimes b_n\mid \cotr(C_{(n)}^\vee))\\
    &\hspace{5cm} \cdot (C_{(1)}^\vee\otimes \one)\shuffle_\bullet (A_2^\vee\otimes b_2^\vee)\shuffle_\bullet \cdots \shuffle_\bullet (A_n^\vee\otimes b_n^\vee),
    \end{align*}
    as $\cotr$ is dual to $\tr$ by Definition \ref{def:cotr} and the shuffle product $\shuffle$ is dual to the coproduct $\Delta_{\glp}(A)=A_{(1)}\otimes A_{(2)}$. Applying the Definition \ref{def:cotr_reduced} of the reduced triangle, we get the claimed formula.
    \end{proof}

\section{Applications to free Lie algebras}	\label{sec:free_Lie}

\subsection{Post-Lie structures on free Lie algebras}
	
From now on, we will study post-Lie structures on free Lie algebra in an abstract setting. We fix a countable set $\Alpha$, which we call \emph{alphabet} and whose elements we call \emph{letters}. The set of words $w$ with letters in $\Alpha$, including the empty word $\one$, will be denoted by $\Alpha^*$. We denote by $(\Q\langle\Alpha\rangle,\cdot)$ the non-commutative free algebra equipped with concatenation. With the shuffle coproduct $\co:\Q\langle\Alpha\rangle\to \Q\langle\Alpha\rangle \otimes\Q\langle\Alpha\rangle$ on the letters given by $\co(a)=a\otimes \one +\one \otimes a$, we obtain a Hopf algebra structure on $(\Q\langle\Alpha\rangle,\cdot)$.
	
We write $(\Lie(\Alpha),\emptybrackets)$ for the free Lie algebra generated by $\Alpha$, which we view as a subset of $\Q\langle\Alpha\rangle$ via the assignment $[u,w]\mapsto uw-wu$. In fact, we can identify  $(\Q\langle\Alpha\rangle,\cdot,\co)$ with the universal enveloping algebra $(\mathcal{U}(\Lie(\Alpha)),\cdot,\Delta)$. This allows, to apply the previously presented results for post-Lie algebras in this special context.

Let $(\Lie (\Alpha),\emptybrackets,\tr)$ be a post-Lie algebra, where $\emptybrackets$ is the usual Lie bracket on $\Lie (\Alpha)$. Since $(\Q\langle \Alpha\rangle,\cdot)$ is the universal enveloping algebra of $(\Lie (\Alpha),\emptybrackets)$, the Grossman-Larson product $\glp$ from Theorem \ref{thm:glp_hopf_algebra} defines another product on $\Q\langle \Alpha\rangle$. Applying Remark \ref{rem:post_lie_recursion} to the free setup yields the following formula for the Grossman-Larson product.

\begin{Proposition} \label{prop:glp_free_lie} For any $B\in \Q\langle \Alpha\rangle$ and letters $a_1,\ldots,a_n \in \Alpha$, we have
	\[B\glp a_1\cdots a_n= B_{(1)}(B_{(2)}\tr a_1)\cdots (B_{(n+1)}\tr a_n).\]
		\end{Proposition}
  \graphref{rem:post_lie_recursion}   
  \graphref{thm:glp_hopf_algebra}
  \noproof{Proposition}

 Moreover, from Theorem \ref{thm:glp_hopf_algebra}, we immediately deduce the following.
	
\begin{Corollary}\label{cor:free_coprod1}
	$(\Q\langle \Alpha\rangle,\glp,\co)$ is a  Hopf algebra.
	\end{Corollary}
 \graphref{thm:glp_hopf_algebra}
\noproof{Corollary}
 	If $\tr$ is nontrivial, the two Hopf algebra structures $(\Q\langle \Alpha\rangle,\glp,\co)$ and $(\Q\langle \Alpha\rangle,\operatorname{conc},\co)$ differ.

	Define a degree map $\deg:\Alpha\to \mathbb{Z}_{\geq1}$ on the alphabet $\Alpha$, such that only finitely many letters have the same degree. This induces a grading on the Lie algebra $(\Lie(\Alpha),\emptybrackets)$ and on the Hopf algebra $(\Q\langle \Alpha \rangle,\cdot,\co)$, where the homogeneous subspaces are finite dimensional. In the following, we view $\Q\langle \Alpha \rangle$ to be graded dual to itself via the pairing
	\begin{equation}\eqlabel{eq:dual_pairing} 
		\begin{aligned}
			\Q\langle \Alpha\rangle \otimes \Q\langle \Alpha\rangle &\longrightarrow \Q\langle \Alpha\rangle, \\
			A\otimes B &\longmapsto (A\mid B),
		\end{aligned}
	\end{equation}
	where we denote by $(A\mid w)$ the coefficient of a word $w\in \A^*$ in $A$ and extend this by $\Q$-linearity.

\begin{Remark} \label{rem:shuffle}
	The shuffle product $\shuffle$ on dual to $\co$ also has the following recursive expression. Set $\one \shuffle A = A = A\shuffle \one$ for all $A\in\Q\langle \Alpha\rangle$, and for $x,y\in\Alpha$, $A,B\in\Q\langle \Alpha\rangle$, let
		\begin{equation}\eqlabel{eq:def_shuffle_product}
			xA \shuffle yB = x(A\shuffle yB) + y(xA\shuffle b)
		\end{equation}
		and extend this by $\Q$-linearity.
\end{Remark}

	Assume that the post-Lie algebra $(\Lie (\Alpha),\emptybrackets,\tr)$ respects the previously defined grading, i.e., also $\tr$ is a homogeneous map. In the free setup, there is a canonical choice for the basis $\mathsf{UB}$ and the generating set $\mathsf{GS}$ in Notation \ref{not:CU_has_UDB}. Namely, we choose $\mathsf{UB}=\Alpha^*$ to be the set of all words and $\mathsf{GS}=\Alpha$ to be the set of all letters. Then, the reduced triangle map from Definition \ref{def:cotr_reduced} is given by
	\begin{align} \eqlabel{eq:reduced_cotr}
	\cotr^{\operatorname{irr}}(A)=\sum_{w\in \Alpha^*,\ a\in \Alpha} (w\otimes a\mid \cotr(A))\ w\otimes a,\qquad A\in \Q\langle\Alpha\rangle,
	\end{align}
    and the product $\shuffle_\bullet$ on $\QV\otimes \QV$ given in \eqref{eq:shuffle_bullet} simplifies to
    \begin{align} \label{eq:sh_bull_free}
    (A_1\otimes A_2)\shuffle_\bullet (B_1\otimes B_2)=(A_1\shuffle B_1)\otimes(A_2B_2).
    \end{align}
	From Theorem \ref{thm:dual_glp_fine}, we obtain the following.
	\begin{Theorem} \label{thm:free_coprod2} For a word $A\in \Alpha^*$, we have
	\[\Delta_{\glp}(A)=\sum_{n\geq1} \sum_{A=A_1\cdots A_n} (A_1\otimes \one)\shuffle_\bullet \cotr^{\operatorname{irr}}(A_2)\shuffle_\bullet \cdots \shuffle_\bullet \cotr^{\operatorname{irr}}(A_n).\]
	\end{Theorem}
\graphref{thm:dual_glp_fine} 
\noproof{Corollary}    

We summarize the previous explained results and considerations on free post-Lie algebras in the following theorem.

\begin{Theorem}\label{thm:dual_pair_free_lie}
Given a graded, free post-Lie algebra $\big( \Lie(\Alpha), \emptybrackets, \tr \big)$, 
    then we get a dual pair of graded Hopf algebras 
    \[
    \xymatrix{ 
     (\Q\langle \Alpha\rangle,\shuffle,\Delta_{\glp}) 
   \ar@{<~>}[rrr]^{\text{graded dual}}&&&    
   (\Q\langle \Alpha\rangle,\glp,\co). 
    }
   \]
Here the Grossman-Larson product $\glp$ is given by Proposition \ref{prop:glp_free_lie} and its dual coproduct $\Delta_{\glp}$ is given by Theorem \ref{thm:free_coprod2}.     
\end{Theorem}
\noproof{Theorem}
  	
	\subsection{A magmatic point of view}
	\label{subsec:magma_post_lie}

	In the following we consider a particular way to attach post-Lie algebras 
	  to free Lie algebras. Our approach is similar to that of Foissy in \cite{Foissy_post_lie}. 
	
	Recall a set with a binary operation is called a magma. We consider here the free magma $M(\Alpha)$ generated by some alphabet $\Alpha$ with the operation 
	\[\star: M(\Alpha) \times M(\Alpha) \to M(\Alpha).\] 
	The $\Q$-vector space $M(\Alpha)_\Q$ spanned by $M(\Alpha)$ is nothing else then the free noncommutative, nonassociative algebra generated by $\Alpha$. 
	A derivation is a linear map $d: M(\Alpha)_\Q \to M(\Alpha)_\Q$ such that $d(a\star b) = da\star b+a\star db$.

\begin{Definition}\label{def:free_tr_magma}	Assume now we are given another map 
	\begin{align*} 
		\tr:  M(\Alpha)_\Q \otimes \Alpha & \to M(\Alpha)_\Q\\
		(t,a) &\mapsto t \tr a, \nonumber
	\end{align*}
	then for each $t \in M(\Alpha)_\Q$ we extend $t \tr$ to a derivation on $(M(\Alpha)_\Q,\star)$. This allows us to view  $\tr$  as a product 
	\begin{align} \eqlabel{eq:tr_ext_M}
	\tr:M(\Alpha)_\Q\otimes M(\Alpha)_\Q\to M(\Alpha)_\Q.
	\end{align}
	Let $\big( \QQ\langle M(\Alpha)\rangle,\, \cdot,\, \Delta \big)$ be the free associative, noncommutative $\Q$-algebra generated by $M(\Alpha)$. We extend the pairing $\tr$ on $M(\Alpha)_\Q$ further to a  linear map
	\begin{align} \eqlabel{eq:tr_ext_AM}
	\tr\colon \QQ\langle M(\Alpha)\rangle \otimes \QQ\langle M(\Alpha)\rangle \to \QQ\langle M(\Alpha)\rangle
	\end{align}
	via the requirements
	\begin{align}
		A \tr \one  &= (A \mid \one ) \eqlabel{eq:nonass_tr_1} \\
		\one\tr A &= A \eqlabel{eq:nonass_tr_2}\\
		(x\cdot A) \tr y &= x\tr(A\tr y)  - (x\tr A)\tr y \eqlabel{eq:nonass_tr_3}\\
		A\tr ( B\cdot C )&= \big(A_{(1)}\tr B\big)\big(A_{(2)}\tr C\big)\eqlabel{eq:nonass_tr_4}
	\end{align}
	for all $A,B,C\in \QQ\langle M(\Alpha)\rangle$ and $x,y\in M(\Alpha)_\Q$.
\end{Definition}
	
	\begin{Proposition} The extension $\tr\colon \QQ\langle M(\Alpha)\rangle \otimes \QQ\langle M(\Alpha)\rangle \to \QQ\langle M(\Alpha)\rangle$ is well-defined.
	\end{Proposition}
	\begin{proof} The above requirements are exactly as in  Definition \ref{def:ext_post-lie}  and the claim follows with the same considerations as  \cite[Prop.~3.1]{EF15}.\end{proof}
	
	The coproduct satisfies
	$\Delta( t ) =t \otimes 1 + 1 \otimes t$ for all $t \in M(\Alpha)$, thus the following proposition holds.
	
	\begin{Lemma}\label{lem:mv_acts_derivation}
		For $t\in M(V)_\Q$, $t_1,\ldots,t_n\in \QQ\langle M(\Alpha)\rangle$, we have
		\[
		t \tr ( t_1 \cdot t_2 \cdot ... \cdot t_n ) = \sum_i  t_1 \cdot ...\cdot (t \tr t_i ) \cdot ... \cdot t_n.
		\]
	In particular, $t\tr$ is a derivation on $\QQ\langle M(\Alpha)\rangle$.
	\end{Lemma}	
	\begin{proof}
		This follows by induction from the requirement \eqref{eq:nonass_tr_4} in Definition \ref{def:free_tr_magma}.    
	\end{proof}
	
	\begin{Lemma} The map $\tr$ restricts to an $\Q$-linear map
	\[\tr:\QQ\langle M(\Alpha)\rangle\otimes M(\Alpha)_\Q\to M(\Alpha)_\Q.\]
	\end{Lemma} 
	\begin{proof}
		This follows directly the recursive requirement \eqref{eq:nonass_tr_3} in Definition \ref{def:free_tr_magma}.    
	\end{proof}
	
	The identity on $\Alpha$ extends naturally to a surjective homomorphism 
	\begin{align*}
		\Lie:\quad M(\Alpha)_\Q &\to \Lie(\Alpha) \\
		a \star b & \mapsto [a,b],
	\end{align*}
	whose kernel is the ideal $I_{\Lie}$ generated by $x\star x$ and by
	\[ 
	J(x,y,z) =
	(x\star y )\star z+(y\star z)\star x+ (z\star x)\star y.
	\] 
	In particular, we have for all $a,b \in M(\Alpha)$ that
	\[
	a \star b + b \star a \in I_{\Lie}.
	\]
	We use the same notation for an element $x \in M(\Alpha)_\Q $ and its equivalence class in $\Lie(\Alpha)\simeq M(\Alpha)_\Q/I_{\Lie}$.

	\begin{Lemma}\label{lem:lie_descend} If $I_{\Lie} \tr \Alpha \subseteq I_{\Lie}$, then the map $\tr$ on $M(\Alpha)_\Q$ descends to 
	\[\tr:\Lie(\Alpha)\otimes\Lie(\Alpha)\to \Lie(\Alpha).\]
	\end{Lemma} 
	\begin{proof} It suffices to show that $I_{\Lie} \tr M(\Alpha)_\Q \subseteq I_{\Lie}$   and $  M(\Alpha)_\Q \tr I_{\Lie} \subseteq I_{\Lie}$. The first formula follows from $I_{\Lie} \tr \Alpha \subseteq I_{\Lie}$ and Lemma \ref{lem:mv_acts_derivation}. The second holds as $t\tr$ is a derivation for all $t\in M(\Alpha)_\Q$.
	Indeed, it is easy to see that any derivation $d$ on $M(\Alpha)_\Q$ satisfies $d(I_{\Lie}) \subseteq I_{\Lie}$, since
		\begin{align*}
			d(x\star x) &= x \star d(x) + d(x) \star x \in I_{\Lie}, \\
			d(J(x,y,z)) &= J(d(x),y,z) + J(x,d(y),z) +J(x,y,d(z)) \in I_{\Lie}. \qedhere
		\end{align*}
	\end{proof}

In other words, if $\tr$ descends to a pairing on $\Lie(\Alpha)$, then the operations $\Lie$ and $\tr$ commute, i.e. for all $x,y \in M(\Alpha)_\Q$ we then have 
 \begin{equation}
\Lie(x\tr y)=\Lie(x) \tr \Lie(y).
 \end{equation}

	\begin{Proposition}\label{prop:tr_descends}
	Assume that $\tr$ descends to a map $\tr: \Lie(\Alpha)\otimes\Lie(\Alpha)\to\Lie(\Alpha)$ and that for all $ x,y,z \in M(\Alpha)_\Q$
		\begin{align*}
			\Lie \Big( \big( x \star y -x\cdot y + y \cdot x \big) \tr z     \Big ) = 0.
		\end{align*}
	Then, $\big( \Lie(\Alpha), \emptybrackets, \tr\big)$ is a post-Lie algebra.
	\end{Proposition}
	\begin{proof}  
	Using the linearity of the maps $\tr$, $\Lie$, and the formula \eqref{eq:nonass_tr_3} in Definition \ref{def:free_tr_magma}, we derive
		\begin{align*}
			0 &=\Lie \Big( \big( x \star y -x\cdot y + y \cdot x \big) \tr z     \Big )\\
			&= [x,y] \tr z - x\tr ( y \tr z) +(x\tr y) \tr z +  y \tr ( x\tr z ) - (y \tr x) \tr z.
		\end{align*}
	This is precisely the condition \eqref{eq:post-lie2} of Definition \ref{def:post_lie_algebra}. Moreover, $t\tr$ is a derivation by construction, therefore the condition \eqref{eq:post-lie1} in Definition \ref{def:post_lie_algebra} is also satisfied.
	\end{proof}

\section{The Ihara bracket} \label{sec:ihara}

The Ihara bracket $\emptybrace_I$ defines another Lie bracket on the free Lie algebra $\Lie(x_0,x_1)$. It origins from Ihara's work on Galois representations on $\mathbb{P}^1 \setminus \{ 0, 1, \infty \}$. The double shuffle Lie
algebra $\mathfrak{dm}_0$ is a Lie  subalgebra, which is central for the theory of multiple zeta values. 

We relate the Ihara bracket to a post-Lie algebra. In fact, we consider a generalisation on $\Lie(V)$, where $V=\{v_0,v_1,v_2,\ldots\}$.
We  use the magmatic approach given in Subsection \ref{subsec:magma_post_lie}
to construct this post-Lie structure.
As an application, we derive a new approach to the Goncharov coproduct on $\Q\langle x_0, x_1 \rangle$.

	\subsection{The Ihara bracket}
	
	\begin{Definition} \label{def:magma_tr_ihara}
	Let $V_0\subset V$ be a subset and consider the triangle map
		\begin{align*}
			\tr_I:  M(V)_\Q \otimes V & \to M(V)_\Q\\
			(t,v) &\mapsto \begin{cases}    0  & v \in V_0  \\     v \star t \quad & v \in V \setminus V_0\, .   \end{cases}
		\end{align*}
		We denote by $\tr_I$ also its extensions to products on $M(V)_\Q$ and $\Q\langle M(V)\rangle $, as explained in \eqref{eq:tr_ext_M} and \eqref{eq:tr_ext_AM}.
	\end{Definition}
	
	\begin{Lemma}\label{lem:ihara_descends}
		The map $\tr_I$ on $M(V)_\Q$ descends to $\tr_I:\Lie(V)\otimes\Lie(V)\to \Lie(V)$.
	\end{Lemma}
	\begin{proof}
	Due to Lemma \ref{lem:lie_descend}, this follows directly from Definition \ref{def:magma_tr_ihara}.
	\end{proof}
	
	\begin{Lemma}\label{lem:ihara_magma_boxed}
	We have for $v \in V \setminus V_0 $  
	\[(t_1 \cdots t_n ) \tr_I v = (\ldots (( v \star t_1) \star t_2 ) \ldots)\star t_n\,.\]
	\end{Lemma}
	
	\begin{proof}
		We prove the claim by induction. Assume it holds for $n-1$ factors, then by Lemma \ref{lem:mv_acts_derivation} and by Definition \ref{def:magma_tr_ihara}
		\begin{align*}
			&(t_1 \cdots t_n ) \tr_I v  = 
			t_1 \tr_I \Big(  (t_2 \cdots t_n ) \tr_I v  \Big) -  \big( t_1 \tr_I (t_2 \cdots t_n )\big)  \tr_I v \\
			&=t_1\tr_I \big(\ldots((v\star t_2)\star t_3)\ldots \star t_n\big)-\Big(\sum_{i=2}^n t_2\cdots (t_1\tr_I t_i)\cdots t_n \Big) \tr_I v\\
			&=  (\ldots((  t_1 \tr_I v ) \star t_2 ) \ldots)\star t_n  + \sum_{i=2}^n  (\ldots(\ldots(( v \star t_2)\star t_3)\ldots\star ( t_1 \tr_I t_i ) )\ldots)\star t_n \\
			&\qquad -  \sum_{i=2}^n   (\ldots(\ldots(( v \star t_2)\star t_3)\ldots\star ( t_1 \tr_I t_i ) )\ldots)\star t_n \\
			&= (\ldots(( v \star t_1) \star t_2 ) \ldots)\star t_n . \qedhere
		\end{align*}
	\end{proof}

	\begin{Theorem}
    \label{thm_ihara_is_post_lie} 
    The triple $( \Lie(V),\emptybrackets, \tr_I)$ is a post-Lie algebra.
	\end{Theorem}
	\begin{proof} By Lemma \ref{lem:ihara_descends} the map $\tr_I$ descends to $\Lie(V)$. Therefore, it remains to verify the second condition in Proposition \ref{prop:tr_descends}. 
		We have by definition of $\tr_I$ 
		\[( t_1 \star t_2  ) \tr_I v = v \star (t_1 \star t_2 ) \]
		and by Lemma \ref{lem:ihara_magma_boxed}
		\[( t_1 \cdot t_2  - t_2 \cdot t_1 ) \tr_I v = ( v \star t_1) \star t_2  - ( v \star t_2) \star t_1.\]
		Applying the Lie-map gives by means of the Jacobi identity
		\begin{align*}
			& \Lie \Big( \big( t_1 \star t_2 -t_1\cdot t_2 + t_2 \cdot t_1 \big) \tr_I v     \Big ) \\
			&\quad =\Lie\Big( v \star (t_1 \star t_2 )   -\big(  ( v \star t_1) \star t_2  - ( v \star t_2) \star t_1    \big)   \Big)\\
			&\quad = 
			[ v ,[t_1 , t_2 ]]   - [[ v , t_1],  t_2]  + [[ v , t_2],  t_1] \\  
			&\quad = 0. \qedhere
		\end{align*}
	\end{proof}
	
	\begin{Definition} We call the from the post-Lie algebra $\big(\Lie(V), \emptybrackets, \tr_I \big)$ induced post-Lie bracket
	\[\{x,y\}_I=x\tr_I y - y \tr_I x +[x,y]\]
	the \emph{Ihara bracket}. 
	\end{Definition}

	\begin{Example}\label{exm:post-lie1}
	We can also consider the restriction of the Ihara bracket to any subsets $X\subset V$ and $X_0=X\mathfrak{g}p V_0$. In particular, the choice $X=\{v_0,v_1\}$ and $X_0=\{v_0\}$ leads to
	\[a \tr_I v_i=\begin{cases} 0, &\quad i=0 \\ [v_1,a], &\quad  i=1 \end{cases} \qquad\qquad (a\in \Lie(X)).\]
	So, the post-Lie bracket $\emptybrace_I$ restricted to $X$ is the usual Ihara bracket, as for example considered in \cite{Ra}, \cite[§~3.10.5]{BGF}. 
	In those works, the Ihara bracket is made with the special derivations \[d_A:\QQ\langle X\rangle \to \QQ\langle X\rangle, \]
 which were for $A\in \QQ\langle X\rangle$ given by $d_A(v_0)=0$ and $d_A(v_1)=[v_1,A]$. Since the derivation $d_A$ is $\Q$-linear by definition, the maps $A \tr_I$ and $d_A$ agree on $\Lie(X)$. But we want to emphasize on the fact that  $A\tr_I$ and $d_A$ differ on $\QQ\langle X\rangle$ outside of $\Lie(X)$. For example,
\begin{align*}
	d_{v_0v_1}(v_1)& = [v_1,v_0v_1],  \\
	v_0v_1 \tr_I v_1 & = v_0 \tr_I ( v_1 \tr_I v_1 ) - (v_0 \tr_I v_1) \tr_I v_1 
	= - [ v_1, [v_1,v_0]] .
\end{align*}
	\end{Example}

	\subsection{The Grossman-Larson product for the Ihara bracket}
	
	The universal enveloping algebra of $\Lie(V)$ is isomorphic to the free non-commutative algebra $\QV$ and by Theorem \ref{thm_ihara_is_post_lie} the extension in Definition \ref{def:ext_post-lie} yields a linear map
	\begin{align}\eqlabel{eq:def_trI}
		\tr_I: \QV \otimes \QV  \to \QV  .
	\end{align}
	
	\begin{Definition} We write the Grossman-Larson product on $\QV$ corresponding to $\tr_I$ as
	\[A\glp_I B= A_{(1)}(A_{(2)}\tr_I B).\]
	\end{Definition} 

	By Theorem \ref{thm:glp_hopf_algebra}, the tuple $(\QV,\glp_I,\co)$ is a Hopf algebra.

	In order to understand the Grossman-Larson product $\glp_I$  associated with the Ihara bracket, it suffices because of  Remark \ref{rem:post_lie_recursion}  to understand $A \tr_I v$ for $A\in \QV$ and $v \in V$.

 \begin{Definition}\label{def:ihara_chi}
For $A \in V^*$ and $v \in V$ we define 
 \[
\chi(A,v) = \begin{cases}
				\delta_{A,\one}, & \quad \text{if} \quad v \in V_0,    \\
				1,  & \quad \text{else},   
			\end{cases}
\]    
 \end{Definition}

\begin{Proposition}\label{prop:post_lie_on_letter}
\graphref{rem:post_lie_recursion}	 
 Let $A\in V^*$ and  $ v \in V$, then 
 we have 
\[
A \tr_I v = \chi(A,v)\, S\big(A_{(1)}\big) v A_{(2)} .\graphref{def:ihara_chi}
\]		
\end{Proposition}

\begin{proof} Because of Definition \ref{def:magma_tr_ihara} there is nothing to show for $v \in V_0$.  
By Lemma \ref{lem:ihara_magma_boxed}, we have for $v\in V \setminus V_0$ and for $A = a_1\cdots a_n\in\QV$ with $a_1,\ldots, a_n \in \Lie(V)$,  that
		\begin{equation}\eqlabel{eq:trI_explizit}
			A\tr_I v = [\dots[[v,a_1],a_2],\dots,a_n].
		\end{equation}
The claim follows by Lemma \ref{lem:boxed_brackets_antipode}.
\end{proof}
	
\begin{Remark}
Because of \eqref{eq:trI_explizit} the above Proposition \ref{prop:post_lie_on_letter} holds also if we replace $\chi(A,v)$ by  $\chi_0(A,v)$, where
\[
\chi_0(A,v) = \begin{cases}
				 				0, & \quad \text{if} \quad A = v a_2 \cdots a_n  ,    \\
			\chi(A,v)	,  & \quad \text{else},   
			\end{cases}
\]  
\end{Remark}	 

\begin{Lemma}\label{lem:tr_on_I-complement}
		For $A,B\in\QV$ and $w\in V_0^*$, we have
		\begin{equation*}
			A\tr_I w B = w (A\tr_I B).
		\end{equation*}
	\end{Lemma}
	
	\begin{proof}
		By \eqref{eq:ext3_asso_asso}, we have
		\begin{equation*}
			A\tr_I w B = \big(A_{(1)}\tr_I w\big) \big(A_{(2)}\tr_I B\big).
		\end{equation*}
	Definition \ref{def:magma_tr_ihara} together with \eqref{eq:ext3_asso_asso} implies that $\big(A_{(1)}\tr_I w\big)$ vanishes unless $A_{(1)} = \one$. Hence we have
		\begin{equation*}
			A\tr_I w B = (\one\tr_I w) (A \tr_I B) = w (A\tr_I B). \qedhere
		\end{equation*}
	\end{proof}

	\begin{Proposition}\label{prop:glp_ihara_grouplike-formula}
		Let $A\in\QV$ and $w = w_1v_{i_1}\cdots w_{d}v_{i_d}w_{d+1}\in V^*$ with $w_j \in V_0^*$ and $v_{i_1},\dots,v_{i_d} \in V \setminus V_0$. Then we have the formula
		\begin{equation*}
			A \glp_I w = A_{(1)} w_1 S\big(A_{(2)}\big) v_{i_1} A_{(3)} w_2\cdots w_{d} S\big(A_{(2d)}\big) v_{i_d} A_{(2d+1)} w_{d+1}.
		\end{equation*}
	\end{Proposition}

	Here, we use the iterated Sweedler notation from Remark \ref{not:sweedler}.
	
	\begin{proof}
		We deduce from the Definition of the Grossman-Larson product \eqref{eq:glp_def}, Proposition \ref{prop:glp_free_lie}, and Lemma \ref{lem:tr_on_I-complement} that
		\begin{align*}
			A\glp_I w &= A_{(1)} (A_{(2)}\tr_I w) = A_{(1)} \big(A_{(2)}\tr_I w_1v_{i_1}\big) \cdots \big(A_{(d+1)}\tr_I w_{d} v_{i_d}\big) \big(A_{(d+2)}\tr_I w_{d+1}\big) \\
			&= A_{(1)} w_1 \big(A_{(2)}\tr_I v_{i_1}\big) \cdots w_{d} \big(A_{(d+1)}\tr_I v_{i_d}\big) w_{d+1}.  
		\end{align*}
	Apply Proposition \ref{prop:post_lie_on_letter} to get the claimed formula.
	\end{proof}

	\subsection{The coproduct for the Ihara bracket} \label{sec:coprod_ihara}
	
	Set $\operatorname{wt}(v_0)=1$ and $\operatorname{wt}(v_i)=i$ for $i\geq1$. This defines a degree map $\operatorname{wt}:V\to \mathbb{Z}_{\geq1}$, which extends to a grading on the Hopf algebra $(\QV,\glp_I,\co)$. We refer to this grading as the \emph{weight grading}. In this subsection, we determine the Hopf algebra \[(\QV,\shuffle,\Delta_I)\] graded dual to $(\QV,\glp_I,\co)$ with respect to the pairing $(-\mid-)$ as given in \eqref{eq:dual_pairing}. To make use of the formula in Theorem \ref{thm:dual_glp_fine}, we will first compute the reduced triangle map from \eqref{eq:reduced_cotr} explicitly.

\begin{Definition}\label{def:ihara_I}
Let $\mathbb{I}:\QV \times V \to \QV$ be the $\Q$-linear map defined  
for   words $w \in V^*$ and a letter $v\in V$ by
\[
\mathbb{I}(w; v) = 
\chi(w,v)  \,  w ,
\]
where the map $\chi$ is given in Definition \ref{def:ihara_chi}. 
\end{Definition}

	\begin{Proposition} \label{prop:redtr_Ihara} For $a_1,\ldots,a_n\in V$, we have \graphref{def:ihara_I}
		\begin{align*}
		\cotr_I^{\operatorname{irr}}(a_1\cdots a_n)=\sum_{j=1}^n \mathbb{I}\big(S(a_1\cdots a_{j-1})\shuffle a_{j+1}\cdots a_n; \, a_j\big) \otimes a_j.
	\end{align*}
	\end{Proposition}

\begin{proof}		
 For $ a \in V$ we have indeed
\begin{align*}
		\cotr_I^{\operatorname{irr}}(a )=  1 \otimes a = \mathbb{I}(1;a) \otimes a ,
\end{align*}
and for $A=a_1\cdots a_n$, we get
 with Proposition \ref{prop:post_lie_on_letter}
\begin{align*}
	\cotr_I^{\operatorname{irr}}(A)
 &=\sum_{w\in V^*,\, a\in V} \big(w\tr_I a \,\big | A\big)\, w\otimes a \\
	&=\sum_{w\in V^*,\, a\in V }   \chi(w,a)\,\big(S(w_{(1)})a w_{(2)} \,\big |\, A\big)\, w\otimes a \\
	&=\sum_{ j=1 }^n \sum_{w\in V^*}  \chi(w,a_j)\,  \big(S(w_{(1)}) a_jw_{(2)} \,\big |\, A\big)\, w \otimes a_j \\
	&=\sum_{ j=1  }^n \sum_{w\in V^*}  \chi(w,a_j)\,\big(\co(w) \,\big |\, S(a_1\cdots a_{j-1})\otimes a_{j+1}\cdots a_n\big) \, w \otimes a_j \\
	&=\sum_{ j=1  }^n \sum_{w\in V^*}  \chi(w,a_j)\,\big(w \,\big |\, S(a_1\cdots a_{j-1})\shuffle a_{j+1}\cdots a_n\big)\, w \otimes a_j \\
	&=\sum_{ j=1  }^n   \mathbb{I}\big(S(a_1\cdots a_{j-1})\shuffle a_{j+1}\cdots a_n; \, a_j\big) \otimes a_j. \qedhere
	\end{align*}
 \end{proof}

With the the explicit formula for the reduced triangle map in Proposition \ref{prop:redtr_Ihara}, we are able to describe the dual coproduct $\Delta_I$ for the Grossman-Larson product $\glp_I$ corresponding to the Ihara bracket.

\begin{Theorem}\label{thm:coprod_ihara} For $a_1,\ldots,a_n\in V$, we have
\begin{align*}
&\Delta_I(a_1\cdots a_n)= \\
&\sum_{\substack{0\leq k \leq n \\ 0\leq i_1< j_1\leq i_2 < j_2 \leq \cdots \leq i_k< j_k \leq n} } \Big( a_1\cdots a_{i_1} \shuffle  
 \mathbb{I}\big( S(a_{i_1+1}\cdots a_{j_1-1})  
 \shuffle a_{j_1+1}\cdots a_{i_2}; a_{j_1} \big) \shuffle\cdots\\
&\hspace{10em} \shuffle  \mathbb{I}\big(S(a_{i_k+1}\cdots a_{j_k-1})\shuffle a_{j_k+1}\cdots a_n  ; a_{j_k} \big)\Big) \otimes a_{j_1}\cdots a_{j_k}.
\end{align*}
\end{Theorem} 

\begin{proof} 
 Let $A=a_1\cdots a_n$ and $a_i\in V$. Then, we have by Theorem \ref{thm:free_coprod2} and Proposition \ref{prop:redtr_Ihara}
\begin{align*}
&\Delta_I(A)=\sum_{k\geq1}\sum_{A=A_1\cdots A_k} (A_1\otimes\one)\shuffle_\bullet \cotr_I^{\operatorname{irr}}(A_2)\shuffle_\bullet \cdots \shuffle_\bullet \cotr_I^{\operatorname{irr}}(A_k) \\
&=\sum_{\substack{0\leq k\leq n \\ 0\leq i_1<\cdots<i_k<n}} (a_1\cdots a_{i_1}\otimes\one) \shuffle_\bullet \cotr_I^{\operatorname{irr}}(a_{i_1+1}\cdots a_{i_2})\shuffle_\bullet \cdots \shuffle_\bullet \cotr_I^{\operatorname{irr}}(a_{i_k+1}\cdots a_n)  \\
&=\sum_{\substack{0\leq k\leq n \\ 0\le i_1<\cdots<i_k< n}} \sum_{j_1=i_1+1}^{i_2}\cdots \sum_{j_k=i_k+1}^{n} \Big( a_1\cdots a_{i_1} \shuffle  
 \mathbb{I}\big( S(a_{i_1+1}\cdots a_{j_1-1})  
 \shuffle a_{j_1+1}\cdots a_{i_2}; a_{j_1} \big)\\
 &\hspace{8em}\shuffle\cdots \shuffle  \mathbb{I}\big(S(a_{i_k+1}\cdots a_{j_k-1})\shuffle a_{j_k+1}\cdots a_n  ; a_{j_k} \big)\Big) \otimes a_{j_1}\cdots a_{j_k}. \qedhere
\end{align*}
\end{proof}

We introduce the following notation.  
	\begin{Notation}\label{not:I_ihara}
		For $v_i, v_j\in V$ and any $f\in \mathcal{U}(\mathfrak{g})$ we set
		\begin{equation*}
			I(v_i; f; v_j) = 
			\begin{cases}
				f, & (v_i,v_j)\in (V\setminus V_0) \times V_0, \\
				S(f) & (v_i,v_j)\in V_0 \times (V\setminus V_0), \\
				(f\mid \one), & (v_i,v_j)\in (V\setminus V_0)^2 \mathcal{U}p V_0^2.
			\end{cases}
		\end{equation*}
In the usual notation for the Goncharov coproduct, as given in  \cite{Gon05}, the evaluation of a symbol  $I( v_i; f; v_j)$ depends on the actual start point $v_i$ and end point $v_j$, whereas here it depends only on being an element in $V_0$ or $V \setminus V_0$.
        
	\end{Notation} 

 \begin{Theorem} \label{thm:coprod_similar_Gon} For $a_1,\ldots,a_n\in V$, we have
\begin{align*}
\Delta_I(a_1\cdots a_n) = \sum_{\substack{0\leq k \leq n \\ 0< j_1<\cdots <j_k\leq n}} \Bigshuffle_{p=0}^k I(a_{j_p}; a_{j_p+1}\cdots a_{j_{p+1}-1};a_{j_{p+1}}) \otimes a_{j_1}\cdots a_{j_k},
\end{align*}
where we formally set $a_{j_0}=a_0\in V\backslash V_0$ and $a_{j_{k+1}}=a_{n+1}\in V_0$.
 \end{Theorem}

Note that this formula is very similar to the well-known formula for the Goncharov coproduct\footnote{In fact, we have not dealt here with the question wether Goncharov's Hopf algebras $\mathscr{S}_\bullet(S)$ and $\tilde{\mathscr{S}}_\bullet(S)$ are dual to the universal enveloping algebra of some post-Lie algebra or not.}.    In Remark \ref{rem:gon_coprod_form} we consider the special case related to multiple zeta values, where  our Notation \ref{not:I_ihara} equals the usual one.

\begin{proof} We first observe that we have for all $A,B\in \QV$ and $v\in V$
\[ \mathbb{I}( A \shuffle B ; v) = \mathbb{I}( A  ; v) \shuffle    \mathbb{I}( B ; v). \]
From Theorem \ref{thm:coprod_ihara}, we derive for $a_1,\ldots,a_n\in V$
\begin{align*}
&\Delta_I(a_1\cdots a_n)= \\
&\sum_{\substack{0\leq k \leq n \\ 0\leq i_1< j_1\leq i_2 < j_2 \leq \cdots \leq i_k< j_k \leq n} } \Big( a_1\cdots a_{i_1} \shuffle  
 \mathbb{I}\big( S(a_{i_1+1}\cdots a_{j_1-1}); a_{j_1} \big)  
 \shuffle \mathbb{I}\big(a_{j_1+1}\cdots a_{i_2}; a_{j_1} \big) \shuffle\\
&\hspace{7em}\cdots \shuffle  \mathbb{I}\big(S(a_{i_k+1}\cdots a_{j_k-1}); a_{j_k} \big)\shuffle \mathbb{I}\big(a_{j_k+1}\cdots a_n  ; a_{j_k} \big)\Big) \otimes a_{j_1}\cdots a_{j_k}.
\end{align*}

For $j_p\in\{0,\ldots,n\}$, we set
\[f(j_p)=\sum_{i_{p+1}=j_p}^{j_{p+1}-1} \mathbb{I}(a_{j_p+1}\cdots a_{i_{p+1}};j_p)\shuffle \mathbb{I}(S(a_{i_{p+1}+1}\cdots a_{j_{p+1}-1};a_{j_{p+1}}) \]
so $f(j_p)$ is the factor which contains the letters $a_{j_p+1}, \ldots, a_{j_{p+1} -1}$. Then, we have with the convention $a_{j_0}=a_0\in V\setminus V_0$ and $a_{j_{k+1}}=a_{n+1}\in V_0$
\[\Delta_I(a_1\cdots a_n)=\sum_{\substack{0\leq k\leq n \\ 0<j_1<\cdots<j_k\leq n}} \Bigshuffle_{p=0}^k f(j_p) \otimes a_{j_1}\cdots a_{j_p}.\]
Thus, it is left to show that \[f(j_p)=I(a_{j_p};a_{j_p+1},\ldots,a_{j_{p+1}-1};a_{j_{p+1}}).\]
If $a_{j_p}, a_{j_{p+1}} \in V_0$, then $f(j_p)$ is only nonzero when $j_p=j_{p+1}-1$. In this case, we have $f(j_p)=1$. In particular, $f(j_p)=I(a_{j_p};a_{j_p+1},\ldots
,a_{j_{p+1}-1};a_{j_{p+1}})$.

If $a_{j_p} \in V\setminus V_0$ and $a_{j_{p+1}} \in  V_0$, then the right of $f(j_p)$ is only nonzero if $i_{p+1}=j_{p+1}-1$. Thus, we get
$f(j_p)=a_{j_p+1}\cdots a_{j_{p+1}-1}=I(a_{j_p};a_{j_p+1},\ldots,a_{j_{p+1}-1};a_{j_{p+1}})$. 

If $a_{j_p}\in V_0$ and $a_{j_{p+1}}\in V \setminus V_0$, then the left factor of $f(j_p)$ is only nonzero if $i_{p+1}=j_p$. Thus, we get $f(j_p)=S(a_{j_p+1}\cdots a_{j_{p+1}-1})=I(a_{j_p};a_{j_p+1},\ldots,a_{j_{p+1}-1};a_{j_{p+1}})$.

If $a_{j_p},a_{j_{p+1}} \in V\setminus V_0$, then we get by applying the antipode relation in the shuffle Hopf algebra, i.e., a special case of Lemma  \ref{lem:dual_antipode_rel},
\begin{align*} 
f(j_p)&=\sum_{i_{p+1}=j_p}^{j_{p+1}-1} a_{j_p+1}\cdots a_{i_{p+1}}\shuffle S(a_{i_{p+1}+1}\cdots a_{j_{p+1}-1}) =\begin{cases} 1 & \text{ if }   j_p+1 = j_{p+1} \\
 0 & \text{ else}
 \end{cases} \\
&=I(a_{j_p};a_{j_p+1},\ldots,a_{j_{p+1}-1};a_{j_{p+1}}). \qedhere
 \end{align*}
\end{proof}

	\begin{Example} Let $V_0=\{v_0\}$. Then, we compute
		\begin{equation*}
			\begin{aligned}
				\Delta_I(v_1v_2v_3v_0) &= \one \otimes v_1v_2v_3v_0 + v_0 \otimes v_1v_2v_3 + v_3 \otimes v_1v_2v_0 + v_3v_0 \otimes v_1v_2 \\
				&\quad + v_2v_3 \otimes v_1v_0 + v_1v_2v_3 \otimes v_0 + v_2v_3v_0 \otimes v_1 + v_1v_2v_3v_0 \otimes \one.
			\end{aligned}
		\end{equation*}
	\end{Example}

\begin{Corollary} (i) For $a_1,\ldots,a_n\in V_0$, we have
\[\Delta_I(a_1\cdots a_n)= \sum_{k=0}^n a_1\cdots a_k \otimes a_{k+1}\cdots a_n =\dec(a_1\cdots a_n). \]
(ii) For $a_1,\ldots,a_n\in V\backslash V_0$, we have
\[\Delta_I(a_1\cdots a_n)= \sum_{k=0}^n a_{k+1}\cdots a_n\otimes a_1\cdots a_k =t\circ \dec(a_1\cdots a_n),\]
where $t:\QV^{\otimes2}\to \QV^{\otimes 2}$ permutes the tensor product factors.
\end{Corollary}
\noproof{Corollary}
 
By Theorem \ref{thm_ihara_is_post_lie} 
the Ihara bracket belongs to a post-Lie algebra,  this yields the following 
an explicit version of \ref{thm:dual_pair_free_lie}.

\begin{Theorem}\label{thm:dual_pair_ihara}
Let $V=\{v_0,v_1,\ldots\}$ be equipped with the grading $\operatorname{wt}(v_0)=1$ and $\operatorname{wt}(v_i)=i$ for $i>0$. 
Consider the graded post-Lie algebra $\big(\Lie(V),\emptybrackets,\tr_I \big)$  associated with the Ihara bracket relative to a subset $V_0 \subset V$. Then
  we get a dual pair of graded Hopf algebras 
    \[
    \xymatrix{ 
     (\Q\langle V\rangle,\shuffle,\Delta_I ) 
   \ar@{<~>}[rrr]^{\text{graded dual}}&&&    
   (\Q\langle V \rangle,\glp_I,\co). 
    }
   \]
An explicit and effective formula for the Grossman-Larson product $\glp_I$ is given by Proposition \ref{prop:glp_ihara_grouplike-formula}  and for its dual coproduct $\Delta_I$ it is given by Theorem \ref{thm:coprod_ihara}.     
\graphref{thm_ihara_is_post_lie} 
\graphref{thm:dual_pair_free_lie}  
\end{Theorem}
\noproof{Theorem}
	
\begin{Remark} \label{rem:gon_coprod_form}	Observe that Theorem \ref{thm:dual_pair_ihara} also holds for any subsets $X \subset V$ and $X_0 = X \mathfrak{g}p V_0$. In particular by the previous discussion, for the choice $X=\{v_0,v_1\}$ and $X_0=\{v_0\}$ the coproduct $\Delta_I$ equals the Goncharov coproduct \cite{Gon05}.
\end{Remark}

\section{The ari bracket} \label{sec:ari}
	
The ari bracket $\emptybrace_a$ origins from Ecalle's theory of moulds and bimoulds \cite{Ecalle02}.  In this section, we study a post-Lie structure on  $\Lie(V)$, where $V=\{v_0,v_1,v_2,\ldots\}$, that 
we identify in Section \ref{sec:ecalle} with a post-Lie algebra structure on alternal bimoulds associated to the ari bracket. 

The description of the ari bracket $\emptybrace_a$ as a post-Lie bracket allows to study the universal enveloping algebra and its dual with the tools developed in Section \ref{sec:free_Lie}. In particular, we derive an explicit, effective formula for the coproduct $\Delta_a$ graded dual to the Grossman-Larson product $\glp_a$.

This Lie bracket plays an important role in the study of the associated depth-graded algebra of the algebra $\Z_q$ of multiple $q$-zeta values.
 
\subsection{The ari multiplicities}\label{subsec:ari_key_lemma}
	
\begin{Notation}\label{nota:indices}
Let $d\geq0$. For $\kvec=(k_1,...,k_d) \in \ZZnonneg^d$, we set 
\[\abs\kvec=k_1+...+k_d,\quad\ell(\kvec)=d,\] 
and $\abs{\emptyset}=\ell(\emptyset)=0$.
For $1\leq i_1<\cdots<i_n\leq d$, we write \[\widehat\kvec^{i_1,\ldots,i_n} = (k_1,...,k_{i_1-1},k_{i_1+1},\ldots, k_{i_n-1},k_{i_n+1},\ldots,k_d)\] for the index where the components $k_{i_1},\ldots,k_{i_n}$ of $\kvec$ are omitted. 
We let $\evec_i = (0,...,0,1,0,...,0)$, thus we have for $s\in \ZZnonneg$
\[\kvec + s \evec_i =   (k_1,...,k_{i-1},k_i+s, k_{i+1},...k_d).\] 
The concatenation of two indices $\kvec=(k_1,\ldots,k_d),\ \lvec=(l_1,\ldots,l_r)$ is denoted by
\[(\kvec,\lvec)=(k_1,\ldots,k_d,l_1,\ldots,l_r).\]
In particular, if for example $\lvec$ is the empty index, then $(\kvec,\lvec)=\kvec$. 
We let $\one =(1,...,1) \in \ZZnonneg^d$. For $\kvec=(k_1,\ldots,k_d)$, $\lvec=(l_1,\ldots,l_d)$ write $\lvec \le \kvec$ if $l_i \le k_i$ for all $i=1,...,d$, and we set
\[\binom{\kvec-\one}{\lvec-\one}=\binom{k_1-1}{l_1-1}\cdots \binom{k_d-1}{l_d-1}.\]
\end{Notation} 
	
	\begin{Definition}\label{def:ari_mult}
		For vectors $\kvec = (k_1,\ldots, k_d),\,
		\lvec =(l_1,\ldots, l_d)\in \ZZnonneg^d$, we define their \emph{ari multiplicity} by
		\begin{equation*}
			\mua_{\kvec,\lvec}  = (-1)^{\abs\kvec + \abs\lvec   } \binom{\kvec-\one}{\lvec-\one}.
		\end{equation*}
    Here, for $k,l\in \ZZnonneg$ the binomial coefficients are understood as follows
    \begin{align*}
        \binom{k-1}{l-1}=\begin{cases} 1, &\quad k=l=0, \\ 0, &\quad k>0, l=0 \text{ or } k=0, l>0.
        \end{cases}
    \end{align*}
    Furthermore, we set for $\kvec=\lvec=\emptyset$ that $\mua_{\kvec,\lvec}=1$.
	\end{Definition}
	
	\begin{Remark}\label{rem:ari_mult_mult}
		\begin{itemize}
			\item[(i)]
			It is 
			\begin{equation}\eqlabel{eq:ari_mult_nonzero}
				\mua_{\kvec,\lvec}  \neq 0
			\end{equation} if and only if  $\lvec \le \kvec$.
			\item[(ii)]
			Let $\kvec,\lvec \in\ZZnonneg^d$. Let $\kvec_1,\lvec_1\in\ZZnonneg^{d_1}$ and $\kvec_2,\lvec_2\in\ZZnonneg^{d_2}$ be any tuples, such that $d_1+d_2=d$ and the $d$-tuples $(\kvec_1,\kvec_2)$, $(\lvec_1,\lvec_2)$ equal a permutation of $\kvec$ and $\lvec$ respectively. Then we have the following multiplicative property
			\begin{equation}\eqlabel{eq:ari_mult_mult}
				\mua_{\kvec,\lvec} = \mua_{\kvec_1,\lvec_1} \cdot \mua_{\kvec_2,\lvec_2}.
			\end{equation}
		\end{itemize}
	\end{Remark}
	
	\begin{Lemma} \label{lem:ari_mult_klr} For all $\kvec,\lvec,\rvec \in \ZZnonneg^d$ we have
		\begin{equation}\eqlabel{eq:ari_mult_klr}
			\mua_{\kvec, \rvec} \cdot \mua_{ \rvec,\lvec} =  (-1)^{\abs\kvec + \abs\rvec} \cdot \mua_{\kvec,\lvec} \cdot \mua_{\kvec-\lvec+\one, \rvec- \lvec+\one}  .
		\end{equation}
	\end{Lemma}	
	\begin{proof} We first observe, that only for $\lvec \le \rvec \le \kvec$ non-zero multiplicities occur.
		The claim follows from the identities $\binom{k-1}{r-1}\binom{r-1}{l-1} = \binom{k-1}{l-1} \binom{k-l}{r-l}$ for $k,l,r\in\ZZnonneg$.
	\end{proof}
	
	\begin{Lemma} \label{lem:ari_mult_identities}	
		For all $a,b \in \ZZnonneg$ and $\kvec \in \ZZnonneg^d$, we have
		\begin{align}\eqlabel{eq:lem_ari_mult_identities}
			\sum_{{\bf r}\in \ZZnonneg^d}   (-1)^{\abs{{\bf k}}+\abs{{\bf r}}} 
			\mua_{{\bf k}+\one,{\bf r} +\one}  \mua_{ a +|{\bf k}-{\bf r}|  , b }    
			&= \mua_{ a , b  - |  {\bf k}|    }
		\end{align}
	\end{Lemma}

	\begin{proof}	
		We prove \eqref{eq:lem_ari_mult_identities} by induction on $d$. The case $d=1$  is the 
        binomial identity\footnote{A proof is given in Corollary \ref{cor:binomial_identities}} 
        \[\sum_{r\in \ZZnonneg} (-1)^r\binom{k}{r}\binom{a+k-r-1}{b-1}=\binom{a-1}{b-k-1}.\]
     	Now assume $d \ge 2$ and \eqref{eq:lem_ari_mult_identities} holds for $d-1 \in\N$. Then we  deduce the claim by applying to the left hand side of \eqref{eq:lem_ari_mult_identities}
		first the case $d=1$ and then the case $d-1$
		\begin{align*}
			\sum_{{\bf r}\in \N_0^d}   (-1)^{\abs{{\bf k}}+\abs{{\bf r}}} &
			\mua_{{\bf k}+\one,{\bf r} +\one}  \mua_{ a +|{\bf k}-{\bf r}|  , b }  \\
			&=\sum_{\widehat{\rvec}^d\in \N_0^{d-1}}  (-1)^{\abs{\widehat{\kvec}^d} +\abs{\widehat{\rvec}^d} } \,\mua_{\widehat{\kvec}^d+\one, \widehat{\rvec}^d+\one}\mua_{a+\abs{\widehat{\kvec}^d - \widehat{\rvec}^d}  , b-k_d} \\ 
			&=  \mua_{a, b-\abs{\kvec}}. \qedhere
		\end{align*}
	\end{proof}
	
	\begin{Remark}
		We could also have made an axiomatic definition for the ari multiplicity by requiring \eqref{eq:ari_mult_nonzero},
		\eqref{eq:ari_mult_mult}, \eqref{eq:ari_mult_klr} and \eqref{eq:lem_ari_mult_identities},
		since in the following we only use these properties. 
	\end{Remark}

	\subsection{The ari bracket} \label{subsec:ari_lie}

	There is a one-to-one correspondence of the elements in the magma $M(V)$ and planar, rooted, binary trees with leaves labeled by elements of $V$. For example the tree
	
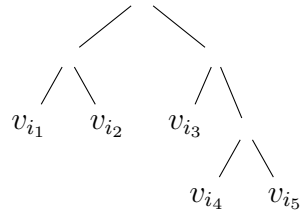
\begin{center}
		\begin{forest}
			for tree={l=0.2cm} 
			[ 
			[ [$v_{i_1}$][$v_{i_2}$]]
			[ [$v_{i_3}$][[$v_{i_4}$][$v_{i_5}$]]]]
			]
		\end{forest}
	\end{center}
	corresponds to  $(v_{i_1}\star v_{i_2})\star (v_{i_3}\star (v_{i_4}\star v_{i_5}))$.
	If we order the leaves of such a tree, e.g. from
	left to the right, then we may consider a tree $t$ with $d$ leaves as a map  
	$V^d \to M(V)$. 
	
	\begin{Notation}\label{nota:t_func}
		Let $t\in M(V)$ be a product of the letters $v_{k_1},\ldots,v_{k_d}$, such that $v_{k_i}$ occurs before $v_{k_{i+1}}$ for $i=1,\ldots,d-1$. Then, $t$ defines a function
		\begin{align*}
			t:V^d&\to M(V), \\
			(v_{l_1},\ldots,v_{l_d}) &\mapsto t(v_{l_1},\ldots,v_{l_d}),
		\end{align*}
		where $t(v_{l_1},\ldots,v_{l_d})$ is obtained from $t$ by replacing each letter $v_{k_i}$ by $v_{l_i}$ for $i=1,\ldots,d$. Even more simplified, $t$ also defines a function
		\begin{align*}
			t:\ZZnonneg^d&\to M(V), \\
			\lvec &\mapsto t(\lvec),
		\end{align*}
		where $t(\lvec)$, $\lvec=(l_1,\ldots,l_d)$ is obtained from $t$ by replacing each index $k_i$ by $l_i$ for $i=1,\ldots,d$.
		
		In particular, we have $t=t(v_{k_1},\ldots,v_{k_d})=t(k_1,\ldots,k_d)$. In the following, we will identify $t\in M(V)$ with its induced functions, and thus we will often also write $t(\kvec)$ instead of $t$.
	\end{Notation}
	
	For example, if $t= v_1\star(( v_0\star v_1 ) \star (v_4 \star v_0))$, 
	then 
	\begin{align*} t&=t(v_1,v_0,v_1,v_4,v_0)=t(1,0,1,4,0), \\
		t(7,0,8,9,10)&=t(v_7,v_0,v_8,v_9,v_{10})=v_7\star(( v_0\star v_8 ) \star (v_9 \star v_{10})).
	\end{align*}
	
We are now prepared to give our second example for the magmatic approach towards a post-Lie structure on $\Lie(V)$ as 
described in Subsection \ref{subsec:magma_post_lie}.	
	\begin{Definition}\label{def:ari_triangle}	
		Define a map
		\begin{align} \eqlabel{eq:tra_def}
			\tr_a:  M(V) \times V & \to M(V)_\Q,\\
			(t(\kvec),v_s) &\mapsto \begin{cases}    0  & s=0, \\  \sum\limits_{\lvec\in \ZZnonneg^{\ell(\kvec)}}  \mua_{\kvec,\lvec} \,\,  v_{s + |\kvec -\lvec|} \star t(\lvec)  \quad & \text{else}.  \\ \end{cases} \nonumber
		\end{align}
		\end{Definition}

		Recall from \eqref{eq:ari_mult_nonzero}, that we have $\mua_{\kvec,\lvec}=0$ for $\lvec>\kvec$. Thus, the sum in the definition \eqref{eq:tra_def} of $\tr_a$ is finite. Note that by Definition \ref{def:ari_mult}, we get non-zero multiplicities if  $k_i=0$ only for $l_i=0$ and if $k_i \ge 1$ only for $l_i \ge 1$.
        
        If it is clear from the context, we will often omit that $\lvec$ is an element in $\ZZnonneg^{\ell(\kvec)}$.
		
		We extend $\tr_a$ to a derivation in the right factor $\tr_a:M(V)\times M(V)_\Q\to M(V)_\Q$, i.e., we require
		\[t\tr_a (u\star v)=(t \tr_a u)\star v + u\star (t\tr_a v)\]
		for all $t\in M(V)$, $u,v\in M(V)$. Furthermore, we extend $\tr_a$ to a $\Q$-linear map \[\tr_a:M(V)_\Q\otimes M(V)_\Q\to M(V)_\Q.\]

	\begin{Lemma}\label{lem:ari_descends}
		The map $\tr_a$ on $M(V)_\Q$ descends to $\tr_a:\Lie(V)\otimes\Lie(V)\to \Lie(V)$.
	\end{Lemma}   
	
	\begin{proof}  By Lemma \ref{lem:lie_descend} it suffices to check $I_{\Lie} \tr_a V \subseteq I_{\Lie}$. Let $x,y,z\in M(V)_\Q$. Then, we have for $s\geq1$
		\begin{align*}
			\big(& x(\kvec)\star x(\kvec) \big) \tr_a v_s =
			\sum_{\lvec}  \mua_{\kvec,\lvec}^2 v_{s +2 |\kvec - \lvec|} \star \big( x(\lvec)\star x(\lvec) \big) \\
			&
			+ \sum_{\lvec_1 <  \lvec_2 }  \mua_{\kvec,\lvec_1}\,\mua_{\kvec,\lvec_2}\,
			v_{s + |\kvec - \lvec_1|+ |\kvec - \lvec_2|} \star \big( x(\lvec_1)\star x(\lvec_2) +x(\lvec_2)\star x(\lvec_1) \big) \in I_{\Lie}.
		\end{align*}
		Next, set $J(\kvec_1,\kvec_2,\kvec_3) = J\big(x(\kvec_1),y(\kvec_2),z(\kvec_3) \big) $, then 
		\begin{align*}
			&J(\kvec_1,\kvec_2,\kvec_3)  \tr_a v_s \\
			&\,= 
			\sum_{\lvec_1,\lvec_2,\lvec_3} \mua_{\kvec_1,\lvec_1}\,\mua_{\kvec_2,\lvec_2}\,\mua_{\kvec_3,\lvec_3}\,
			v_{s + |\kvec_1 - \lvec_1|+ |\kvec_2 - \lvec_2|+ |\kvec_3 - \lvec_3|} \star
			J(\lvec_1,\lvec_2,\lvec_3)  \in I_{\Lie}. \qedhere
		\end{align*}
	\end{proof}
	
	Extend the map $\tr_a:M(V)_\Q\otimes M(V)_\Q\to M(V)_\Q$ to a $\Q$-linear map
	\[\tr_a:\Q\langle M(V)\rangle \otimes \Q\langle M(V)\rangle\to \Q\langle M(V)\rangle\]
	via the conditions \eqref{eq:nonass_tr_1}, \eqref{eq:nonass_tr_2}, \eqref{eq:nonass_tr_3}, and \eqref{eq:nonass_tr_4}.

	\begin{Lemma} \label{lem:ari_magma_boxed} We have for $s\in \mathbb{Z}_{>0}$ and
		$t_1(\kvec_1) , \ldots, t_n(\kvec_n) \in M(V)$
		\[
		(t_1(\kvec_1) \cdots t_n(\kvec_n) ) \tr_a v_s = \sum_{\lvec_1,\ldots,\lvec_n  } 
		 \prod_{i=1}^n \mua_{\kvec_i,\lvec_i} 
		\, \big(...\big(\big( v_{s + \sum\limits_{i=1}^n ( |\kvec_i  - \lvec_i|)}  \star t_1(\lvec_1)\big) \star t_2(\lvec_2)\big) ... \big)\star t_n(\lvec_n)\,.
		\]
	\end{Lemma}
	\begin{proof}	We give here a proof that relies on an identity given in Lemma \ref{lem:key_ari}, which we prove  below.
		We prove the statement by induction on $n$. We have
		\begin{align*}
			&\big( t_1(\kvec_1)   \cdot t_2(\kvec_2)\cdots t_n(\kvec_n) \big) \tr_a v_s \\
			&= t_1(\kvec_1) \tr_a \Big( \big( t_2(\kvec_2)   \cdot t_3(\kvec_3)\cdots t_n(\kvec_n) \big) \tr_a v_s \Big) \\
			&\hspace{6cm}- \Big(  t_1(\kvec_1) \tr_a\big( t_2(\kvec_2)   \cdot   t_3(\kvec_3)  \cdots t_n(\kvec_n) \big)\Big) \tr_a v_s\\	
			&\overset{\text{I.H.}}{=} 
			\sum_{ \lvec_2,...,\lvec_n } 
			\Big( \prod_{i=2}^n \mua_{\kvec_i,\lvec_i} \Big)
			\, t_1(\kvec_1) \tr_a
			\Big(\big(...\big(\big( v_{s + \sum\limits_{i=2}^n ( |\kvec_i  - \lvec_i|)}  \star t_2(\lvec_2)\big) \star t_3(\lvec_3)\big) ...\big)\star t_n(\lvec_n)\Big)
			\\		 
			& \hspace{0,4cm} - \sum_{i=2}^n \Big( t_2(\kvec_2)   \cdots \big( t_1(\kvec_1) \tr_a t_i(\kvec_i) \big) \cdots t_n(\kvec_n) \Big) \tr_a v_s\\
			&=  \sum_{ \lvec_1,...,\lvec_n  }
			\Big( \prod_{i=1}^n \mua_{\kvec_i,\lvec_i} \Big)
			\, \big(...\big(\big( v_{s + \sum\limits_{i=1}^n ( |\kvec_i  - \lvec_i|)}  \star t_1(\lvec_1)\big) \star t_2(\lvec_2)\big) ...\big)\star t_n(\lvec_n)\\
			&\hspace{0,4cm} +\sum_{i=2}^n 
			\sum_{ \lvec_2,...,\lvec_n  }
			\Big( \prod_{i=2}^n \mua_{\kvec_i,\lvec_i} \Big)
			\, 
			\big(... \big(\big(...\big( v_{s + \sum\limits_{i=2}^n ( |\kvec_i  - \lvec_i|)}  \star t_2(\lvec_2)\big)...\big) \\
			&\hspace{8,3cm}\star \big( t_1(\kvec_1)\tr_a t_i(\lvec_i)\big)\big)...\big)\star t_n(\lvec_n)
			\\		 
			&\hspace{0,4cm} - \sum_{i=2}^n \Big( t_2(\kvec_2)   \cdots \big( t_1(\kvec_1) \tr_a t_i(\kvec_i) \big) \cdots t_n(\kvec_n) \Big) \tr_a v_s
		\end{align*}
		It remains to show that the two sums over the index $i$ cancel, this is the content of Lemma \ref{lem:key_ari} below. 
	\end{proof}
	
	\begin{Notation} \label{not:ari_a} 
		Let $t_1(\kvec), t_2({\bf n}) \in M(V)$ with $\kvec\in\ZZnonneg^{d_1}$ and ${\bf n}=(n_1,\ldots,n_{d_2})\in\ZZnonneg^{d_2}$.
		For $j\in\{1,\ldots,d_2\}$, we define
		\[
		a^{(t_1,t_2)}_j( \kvec , {\bf n}  ) =  \begin{cases} t_2( v_{n_1}, \ldots, v_{n_{j-1}}, 
		v_{n_j } \star t_1(  \kvec ), v_{n_{j+1}}, \ldots, v_{n_{d_2}} ), &\quad n_j>0, \\
		0, &\quad n_j=0. \end{cases}
		\] 
		So for $n_j>0$, the element $a^{(t_1,t_2)}_j( \kvec , {\bf n}  )\in M(V)$ is obtained from $t_2({\bf n})$ by replacing the letter $v_{n_j}$ by the product $v_{n_j}\star t_1(\kvec)$. In terms of trees, $a_j^{(t_1,t_2)}(\kvec,{\bf n})$ is obtained from $t_2$ by replacing the leaf labeled by $v_{n_j}$ with the tree $t_1$.
		
		This notation should be seen as an extension of Notation \ref{nota:t_func}, since we allow now that letters get not only replaced by letters, but also by products in $M(V)$. 
	\end{Notation}
	
	With Notation \ref{not:ari_a}, we get the following simplified expression of $\tr_a$
	\begin{align*}
		t_1(\kvec) \tr_a  t_2( {\bf n} ) &=  
		\sum_{j=1}^{d_2}  \sum_{\lvec}  \mua_{\kvec,\lvec}  \, t_2( v_{n_1}, \ldots, v_{n_{j-1}}, 
		v_{n_j +|\kvec - \lvec| } \star t_1(\lvec), v_{n_{j+1}}, \ldots, v_{n_{d_2}} )  \\
		&= \sum_{j=1}^{d_2}  \sum_{\lvec}  \mua_{\kvec,\lvec}  \, a_j^{(t_1,t_2)}(\lvec ,  {\bf n} + |\kvec-\lvec| \evec_j ).  
	\end{align*}

	\begin{Lemma} \label{lem:key_ari} Assume Lemma \ref{lem:ari_magma_boxed} holds for $n-1$ factors, then
	we have for all $i=2,...,n$
		\begin{align*}
			&   \big(  t_2(\kvec_2)\cdots  
			\big( t_1(\kvec_1) \tr_a t_i(\kvec_i) \big)  \cdots t_n(\kvec_n)\big) \tr_a v_s\\
			&=\sum_{ \lvec_2,\ldots,\lvec_n  }
			\Big( \prod_{i=2}^n \mua_{\kvec_i,\lvec_i} \Big)
			\, 
			\big(... \big(\big(...\big( v_{s + \sum\limits_{i=2}^n ( |\kvec_i  - \lvec_i|)}  \star t_2(\lvec_2)\big)...\big) \star \big( t_1(\kvec_1)\tr_a t_i(\lvec_i)\big)\big)...\big)\star t_n(\lvec_n).
		\end{align*}

	\end{Lemma}
	\begin{proof}  
		Let $\kvec_i\in \N^{d_i}$ for $i=1,\ldots,n$. With Notation \ref{not:ari_a} and $[n]_i=\{2,\ldots,n\}\backslash\{i\}$ we need to show 
		\begin{align*}		
			&\sum_{j= 1}^{d_i} 
			\sum_{\substack{\lvec_t, \,\, t\in [n]_i    \\ {\bf n}_1, {\bf n}_i, {\bf r} }}  
			\Big( \prod_{\substack{ t \in [n]_i }}   \mua_{\kvec_t,  \lvec_t } \Big)
			\mua_{ \kvec_{i} +|\kvec_1-{\bf r}|\evec_j  , {\bf n}_i }   \mua_{\kvec_1,{\bf r}}  \mua_{{\bf r},{\bf n}_1} \\
			& \hspace{0,4cm}  
			\big(...\big(\big(...\big( v_{s+ \sum\limits_{t \in [n]_i} |\kvec_t -\lvec_t| + |(\kvec_{i} +|\kvec_1-{\bf r}|\evec_j)  - {\bf n}_i | 
				+  |{\bf r}-{\bf n}_1|   }  \star t_2(\lvec_2)\big)... \big) 
			\star   a_j^{(t_1,t_i)}( {\bf n}_1, {\bf n}_i)\big) ...\big)   \star t_n(\lvec_n)\\
			&=  \sum_{j= 1}^{d_i}  \sum_{\lvec_1,\ldots,\lvec_n}  \Big( \prod_{i=1}^n \mua_{\kvec_i,\lvec_i} \Big)  
			\, \big(...\big(\big(...\big( v_{s+\sum\limits_{i=2}^n |\kvec_i  - \lvec_i |} \star t_2(\lvec_2) \big)...\big)  
			\\
			&\hspace{7,8cm} \star a_j^{(t_1,t_i)} (\lvec_1, \lvec_{i} +|\kvec_1 - \lvec_1| \evec_j)\big) ...\big)\star t_n(\lvec_n), \\
		\end{align*}
		We compute the coefficients of 
		\begin{align*}
			&  \big(...\big(\big(...\big( v_{s+  \sigma  }  \star t_2(\lvec_2)\big)...\big) 
			\star   a_j^{(t_1,t_i)}( {\bf n}_1, {\bf n}_i)    \big) ...\big)   \star t_n(\lvec_n)
		\end{align*}
		in both sides of the previous equation.
		For the right-hand side we must have 
		${\bf n}_1= \lvec_1 $ and ${\bf n}_i = \lvec_{i } +|\kvec_1 - \lvec_1| \evec_j$, 
		and thus the coefficient is
		\begin{align*}
			C_r &= \Big( \prod_{\substack{t \in [n]_i }}   \mua_{\kvec_t,  \lvec_t} \Big)   
			\mua_{\kvec_1,{\bf n}_1} \mua_{ \kvec_i, {\bf n}_i  - |\kvec_1 - {\bf n}_1| \evec_j  }\\
			&=\Big( \prod_{\substack{t \in [n]_i }}   \mua_{\kvec_t,  \lvec_t} \Big)   
			\mua_{\kvec_1,{\bf n}_1}   \mua_{ \widehat{\kvec_i}^j, \widehat{{\bf n}_i}^j}   \mua_{ \kvec_{i,j} , {\bf n}_{i,j}  - |\kvec_1 - {\bf n}_1|    }.
		\end{align*}
		In the last equality we employed the multiplicativity of the $\mua_{\kvec,\lvec}$ given in \eqref{eq:ari_mult_mult}.	
		For the left-hand side the coefficient is given by 
		\begin{align*}
			C_l&=  \Big( \prod_{\substack{t \in [n]_i }}   \mua_{\kvec_t,  \lvec_t} \Big) \cdot\, \sum_{{\bf r}}
			\mua_{ \kvec_{i } +|\kvec_1-{\bf r}|\evec_j  , {\bf n}_i }   \mua_{\kvec_1,{\bf r}}  \mua_{{\bf r},{\bf n}_1} \\
			&=  \Big( \prod_{\substack{t \in [n]_i }}   \mua_{\kvec_t,  \lvec_t} \Big) \cdot  \mua_{\kvec_1,{\bf n}_1} \cdot
			\mua_{ \widehat{\kvec_i}^j, \widehat{{\bf n}_i}^j} \cdot   \, \sum_{{\bf r}}   (-1)^{\abs{\kvec_1} + \abs{\bf r}} 
			\mua_{\kvec_1 - {\bf n}_1+\one,{\bf r}-{\bf n}_1+\one}  \mua_{ \kvec_{i,j } +|\kvec_1-{\bf r}|  , {\bf n}_{i,j} }    
		\end{align*}
		In the last step  we used again \eqref{eq:ari_mult_mult} and in addition  Lemma \ref{lem:ari_mult_klr}. Therefore, we have $C_l=C_r$ if and only if
		\begin{align*}
			\sum_{{\bf r}}   (-1)^{\abs{\kvec_1} + \abs{\bf r}} 
			\mua_{\kvec_1 - {\bf n}_1+\one,{\bf r}-{\bf n}_1+\one}  \mua_{ \kvec_{i,j } +|\kvec_1-{\bf r}|  , {\bf n}_{i,j} }    
			&= \mua_{ \kvec_{i,j} , {\bf n}_{i,j}  - |\kvec_1 - {\bf n}_1|}.
		\end{align*}
		Substitute ${\bf r}-{\bf n}_1={\bf r}$, then this is exactly the formula from Lemma \ref{lem:ari_mult_identities} with $a=\kvec_{i,j},\ b=n_{i,j}$, and ${\bf k}=\kvec_1-{\bf n}_1$. 
		 \end{proof}

	\begin{Theorem} \label{thm_ari_is_post_lie} The tuple $(\Lie(V), [\cdot,\cdot], \tr_a)$ is a post-Lie algebra.
	\end{Theorem}
	\begin{proof}  By Lemma \ref{lem:ari_descends} the map $\tr_a$ descends to a map on $\Lie(V)$. Therefore, it remains to check the second condition in Proposition \ref{prop:tr_descends}. 
		
		Let $t_1(\kvec_1), t_2(\kvec_2)\in M(V)$. By definition of $\tr_a$ we have for $s\geq1$
		\[
		( t_1(\kvec_1) \star t_2 (\kvec_2) ) \tr_a v_s = 
		\sum_{\lvec_1,\lvec_2 } \mua_{\kvec_1 ,\lvec_1 }\mua_{\kvec_2,\lvec_2} \,
		v_{s +  |\kvec   - \lvec |  } \star (t_1(\lvec_1) \star t_2(\lvec_2) ), 
		\]
		where we abbreviate $|\kvec-\lvec|=|\kvec_1-\lvec_1|+|\kvec_2-\lvec_2|$. On the other hand, we have by Lemma \ref{lem:ari_magma_boxed}
		\begin{align*}
			( t_1(\kvec_1) \cdot& t_2(\kvec_2)  - t_2(\kvec_2) \cdot t_1( \kvec_1) ) \tr_a v_s \\
			&= 
			\sum_{ \lvec_1,\lvec_2 } 
			\mua_{\kvec_1 ,\lvec_1 }\mua_{\kvec_2,\lvec_2}\,   \Big( ( v_{s + |\kvec  - \lvec |  } 
			\star t_1(\lvec_1)) \star t_2(\lvec_2)  - ( v_{s +  |\kvec   - \lvec |  } \star t_2(\lvec_2)) \star t_1(\lvec_1) \Big).
		\end{align*}
		
		Applying the Lie-map gives by means of the Jacobi relation for all $(\lvec_1 , \lvec_2)$
		\begin{align*}
			&\Lie\Big( v_{s +   |\kvec  -  \lvec |} \star (t_1(\lvec_1) \star  t_2(\lvec_2) )   -\big(  ( v_{s +   |\kvec    -  \lvec |} 
			\star t_1(\lvec_1)) \star  t_2(\lvec_2)  - ( v_{s +  |\kvec   -  \lvec | } \star  t_2(\lvec_2)) \star t_1(\lvec_1)    \big)   \Big)\\
			&\quad = 
			[ v_{s +   |\kvec  -  \lvec |} ,[t_1(\lvec_1) ,  t_2(\lvec_2) ]]   - [[ v_{s + |\kvec   -  \lvec |} , t_1(\lvec_1)],   t_2(\lvec_2)]  + [[ v_{s + |\kvec   -  \lvec |},  t_2(\lvec_2)],  t_1(\lvec_1)] \\  
			&\quad = 0. \qedhere
		\end{align*}
	\end{proof}
	
	An immediate consequence of Proposition \ref{prop:post-lie_relates_lie-algebras} is the following.
	
	\begin{Corollary}
		The pair $( \Lie(V),\emptybrace_a)$ is a Lie algebra, where
		\[\{x,y\}_a=x\tr_a y- y\tr_a x + [x,y].\] 
	\end{Corollary}
	
	\begin{Remark}
		We call the induced post-Lie bracket $\emptybrace_a$ following Ecalle \cite{Ecalle02} the ari bracket, as it translates into the ari bracket of bimoulds, cf Section \ref{sec:ecalle}.
		In fact, it also follows from this non-trivial translation of our setup into the language of alternal bi-moulds that $\emptybrace_a$ is a Lie bracket.
	\end{Remark}

	\subsection{The Grossman-Larson product for the ari bracket}

	The universal enveloping algebra
	$\mathcal{U}(\Lie(V))$ is isomorphic to the free non-commutative algebra $\QV$ and the extension in Definition \ref{def:ext_post-lie} yields a linear map
	\begin{align}\eqlabel{eq:tra_on_QV}
		\tr_a: \QV \otimes \QV  \to \QV  .
	\end{align}
	Because of  Remark \ref{rem:post_lie_recursion} it suffices to understand $A \tr_a v$ for $A\in \QV$ and $v \in V$.
	
	\begin{Remark} \label{exm:ari} 
		Evidently, the map $\tr_a$ in Definition \ref{def:ari_triangle} can be easily extended to a derivation on $\Q\langle V\rangle$ by allowing instead of Lie brackets also arbitrary $\Q$-linear combinations of words. As for the Ihara bracket explained in Example \ref{exm:post-lie1}, this naive extension does not agree with the map \eqref{eq:tra_on_QV} obtained from Definition \ref{def:ext_post-lie}.
	\end{Remark}

	\begin{Notation} \label{nota:t_func_on_QV} We want to use the Notation \ref{nota:t_func} also for $\QV$. Any word $A = v_{k_1}\cdots v_{k_d}\in \QV$, can also be seen as a function
		\begin{align*}
			A:\ZZnonneg^d&\to \QV, \\
			(l_1,\ldots,l_d)&\mapsto A(l_1,\ldots,l_d)=v_{l_1}\cdots v_{l_d}.
		\end{align*}
		As before, we have $A=A(k_1,\ldots,k_d)$ and from now on we will use both notations.
		
		Furthermore, we extend this to the Sweedler notation
		\begin{align} \label{eq:sweed_ex}
		\Delta A(\kvec) = A(\kvec)_{(1)} \otimes A(\kvec)_{(2)} =  A_{(1)}( \kvec_1)  \otimes A_{(2)}(\kvec_2).
		\end{align}
	\end{Notation}
	
	\begin{Example}
		Let $A(\kvec)=v_{k_1}v_{k_2}$. Then, we compute
		\begin{align*}
			\Delta A(\kvec)&= A(\kvec)\otimes\one + v_{k_1} \otimes v_{k_2}+v_{k_2}\otimes v_{k_1}+ \one\otimes A(\kvec).
		\end{align*}
    So for example in the third term, we have $A(\kvec)_{(1)}=v_{k_2}$ and $\kvec_1=(k_2)$.
	\end{Example}

	As indicated in Remark \ref{rem:post_lie_recursion} we have the following efficient way to calculate $\tr_a$. 
	
	\begin{Proposition}\label{prop:mult_tr_ari} \graphref{rem:post_lie_recursion}
		For all $A(\kvec) \in \QV$, we have
		\begin{equation} \eqlabel{eq:tr_iterated_brackets_ari}
			A(\kvec)  \tr_a v_s =
			\begin{cases}
			 	(A(\kvec)|\one) v_0 & \quad \text{if} \quad s=0,\,   \\
				\sum\limits_{\lvec }  \, \mua_{\kvec,\lvec}  \,       S\big(A(\lvec)_{(1)}\big) v_{s+|\kvec-\lvec|} A(\lvec)_{(2)} & \quad  \text{else}. 
			\end{cases}
		\end{equation}
	\end{Proposition} 

By the convention in Definition \ref{def:ari_mult}, we have $\one \tr_a v=v$ for all $v\in V$.

	\begin{proof} Let $\kvec=(\kvec_1,\ldots,\kvec_n)$ and $A(\kvec)=a_1(\kvec_1)\cdots a_n(\kvec_n)\in \QV$ with $a_i(\kvec_i)\in \Lie(V)$. With Lemma \ref{lem:ari_magma_boxed} and Lemma \ref{lem:boxed_brackets_antipode} we obtain for $s\geq1$
		\begin{align*}
			A(\kvec)\tr_a v_s&=\sum_{\lvec=(\lvec_1,\ldots,\lvec_n)} \mua_{\kvec,\lvec} [\ldots [[v_{s+|\kvec-\lvec|},a_1(\lvec_1)],a_2(\lvec_2)],\ldots,a_n(\lvec_n)] \\
			&= \sum_{\lvec} \mua_{\kvec,\lvec}       S\big(A(\lvec)_{(1)}\big) v_{s+|\kvec-\lvec|} A(\lvec)_{(2)}. \qedhere
		\end{align*} 
	\end{proof}

	\begin{Example}
		One computes that
		\begin{align*}
			v_2v_3 \tr_a v_3 &= -v_6v_1v_1 - v_1v_1v_6 + v_5v_2v_1 + 2\, v_2v_1v_5 - v_3v_1v_4 + v_3v_4v_1 - v_4v_1v_3 \\
			&\quad  + 2\, v_5v_1v_2 + v_1v_2v_5 + v_1v_4v_3 - 2\, v_4v_2v_2 - 2\, v_2v_2v_4 - v_3v_3v_2 \\
			&\quad + 2\, v_3v_2v_3 - v_2v_3v_3 + 2\, v_1v_6v_1 - 3\, v_2v_5v_1 - 3\, v_1v_5v_2 + 4\, v_2v_4v_2.
		\end{align*}
	\end{Example}

	From Theorem \ref{thm:glp_hopf_algebra}, we deduce that 
	\[(\QV,\glp_a,\co)\]
	is a Hopf algebra with the product $\glp_a$ induced by $\tr_a$ as in \eqref{eq:glp_def}. The next proposition gives an explicit formula for the product $\glp_a$.
	
	\begin{Proposition}\label{prop:glp_ari_formula} 
		Let $A(\kvec)\in\QV$ and $w = w_1v_{i_1}\cdots w_{d}v_{i_d}w_{d+1}\in V^*$ with  $w_1,\dots,w_{d+1} \in \{v_0^m\mid m\geq0\}$ and $v_{i_1},\dots,v_{i_d}\in V \setminus \{v_0\}$. Then,
		\begin{align*}
			A(\kvec) \glp_a w = \sum_{\lvec_2,\ldots \lvec_{2d+1} }  
			& \prod_{j=2}^{2d+1} \mua_{\kvec_j,\lvec_j}\, A_{(1)}(\kvec_1) w_1 S\big(A_{(2)}(\lvec_2)\big) 
			v_{i_1+\abs{ \kvec_2+\kvec_3}- \abs{\lvec_2+\lvec_3  }} A_{(3)}(\lvec_3) w_2 \\
			&\quad  \cdots w_{d} S\big(A_{(2d)}(\lvec_{2d})\big) v_{i_d +\abs{\kvec_{2d}+\kvec_{2d+1}} -\abs{\lvec_{2d}+\lvec_{2d+1}} } A_{(2d+1)}(\lvec_{2d+1} )w_{d+1}.
		\end{align*}
	\end{Proposition}
	Here, we use the iterated extended Sweedler notation as in Notation \ref{not:sweedler} and \eqref{eq:sweed_ex}.
	
	\begin{proof} We first observe that Lemma \ref{lem:tr_on_I-complement} also holds for $\tr_a$.
		Then by Proposition \ref{prop:glp_free_lie} we obtain
		\begin{align*}
			A(\kvec) \glp_a w = A_{(1)}(\kvec_1) w_1 (A_{(2)}(\kvec_2)\tr_a v_{i_1})\cdots w_d(A_{(d+1)}(\kvec_{d+1})\tr_a v_{i_d})w_{d+1}.
		\end{align*}
		Applying Proposition \ref{prop:mult_tr_ari}, we get the claim.
	\end{proof}
	
	\begin{Example}
		One computes
		\begin{align*}
			v_2v_3 \glp_a v_1v_0 &= 3\, v_1v_2v_3v_0 - 2\, v_2v_3v_1v_0 - 2\, v_3v_1v_2v_0 - 3\, v_1v_3v_2v_0 \\
			&\hspace{0.4cm}+ 5\, v_3v_2v_1v_0 - v_4v_1v_1v_0 + 2\, v_1v_4v_1v_0 - v_1v_1v_4v_0.
		\end{align*}
	\end{Example}

	\subsection{The coproduct for the ari bracket}

    As in Subsection \ref{sec:coprod_ihara}, we define a weight grading on $\QV$ by setting 
    \begin{align} \label{eq:wt} \operatorname{wt}(v_0)=1,\qquad \operatorname{wt}(v_s)=s \text{ for } s\geq1.\end{align}
    We also define a \emph{depth grading} on $\QV$ by
    \begin{align} \label{eq:dep} \operatorname{dep}(v_i)=1 \text{ for } i\geq0.\end{align}
    
    \begin{Lemma} \label{lem:ari_bigrad}
    The post-Lie algebra $(\Lie(V),\emptybrackets,\tr_a)$ from Theorem \ref{thm_ari_is_post_lie} is graded with respect to the weight and depth.
    
\end{Lemma}
\begin{proof}
This is an immediate consequence from the the formula in Definition \ref{def:ari_triangle}.
\end{proof}

We derive from Lemma \ref{lem:ari_bigrad} that $(\Lie(V),\emptybrace_a)$ is a Lie algebra graded by weight and depth, and hence also the universal enveloping algebra $(\QV,\glp_a,\co)$ is bigraded with respect to the weight and depth. We want to describe the Hopf algebra
    \[(\QV,\shuffle,\Delta_a)\]
	graded dual to $(\QV,\glp_a,\co)$ as explicit as possible.
      First, we determine the reduced triangle map \eqref{eq:reduced_cotr} by using Notation \ref{nota:indices}.

	\begin{Proposition} \label{prop:red_tr_ari} 
   For $\kvec=(k_1,\ldots,k_d)\in \ZZnonneg^d$ with $d\ge 1$, we have  
	\begin{align*}
	\cotr_a^{\operatorname{irr}} ( v_{k_1} \ldots v_{k_d})
			= &\sum_{i=1}^d \sum_{\substack{\lvec=(l_1,\ldots,\widehat{l_i},\ldots,l_d)\in \ZZnonneg^{d-1} \\ \lvec=\emptyset \text{ or }|\lvec|\leq k_i-1}}  \\
        &\mua_{\widehat\kvec^i +\lvec,\widehat\kvec^i } \big(S( v_{k_1 +l_1} \ldots v_{k_{i-1}+l_{i-1}}) \shuffle (v_{k_{i+1} + l_{i+1}}  \ldots v_{k_n +l_n} )\big) \otimes v_{k_i-|\lvec|}.
	\end{align*} 
    \end{Proposition}  
Observe the case $\lvec = \emptyset$ only occurs if $d=1$ and then 
\[
\cotr_a^{\operatorname{irr}}(v)=\one\otimes v
\]
for any    letter $v\in V$.
If $d\geq2$ and we have $k_i=0$ for some $i$, the summation index of the sum $\sum_{\substack{\lvec=(l_1,\ldots,\widehat{l_i},\ldots,l_d)\in \ZZnonneg^{d-1}, \ |\lvec|\leq k_i-1}}$ is empty and hence the sum vanishes by convention. So for any word $w$ with at least two letters, the letter $v_0$ will never occur as a right factor in $\cotr_a^{\operatorname{irr}}(w)$. In particular, $\cotr_a^{\operatorname{irr}}(v_0^m)=0$ for all $m\geq2$. 

    \begin{proof} For $d=1$, the formula $\cotr_a^{\operatorname{irr}}(v)=\one\otimes v$ evidently holds. Thus we assume $d\geq2$ and compute with Proposition \ref{prop:mult_tr_ari} for $v_{k_1},\ldots,v_{k_d}\in V$
    \begin{align*}
    &\cotr_a^{\operatorname{irr}}(v_{k_1}\cdots v_{k_d})=\sum_{w\in V^*,\ a\in V} (w\tr_a a\mid v_{k_1}\cdots v_{k_d})\ w\otimes a \\
    &=\sum_{w(\mathbf{r})\in V^*, s\geq1} \sum_{\lvec} \mua_{\mathbf{r},\lvec} \Big(S(w(\lvec)_{(1)})v_{s+|\mathbf{r}-\lvec|}w(\lvec)_{(2)}\mid v_{k_1}\cdots v_{k_d}\Big)\ w(\mathbf{r})\otimes v_s \\
	&=\sum_{i=1}^d \sum_{w(\mathbf{r})\in V^*} \sum_{\substack{\lvec, \\ |\mathbf{r}-\lvec|<k_i}} \mua_{\mathbf{r},\lvec} \Big(S(w(\lvec)_{(1)})v_{k_i}w(\lvec)_{(2)}\mid v_{k_1}\cdots v_{k_d}\Big)\ w(\mathbf{r})\otimes v_{k_i-|\mathbf{r}-\lvec|} \\
	&=\sum_{i=1}^d \sum_{w(\mathbf{r})\in V^*} \sum_{\substack{\lvec, \\ |\mathbf{r}-\lvec|<k_i}} \mua_{\mathbf{r},\lvec} \Big(w(\lvec)\mid S(v_{k_1}\cdots v_{k_{i-1}})\shuffle v_{k_{i+1}}\cdots v_{k_d}\Big)\ w(\mathbf{r})\otimes v_{k_i-|\mathbf{r}-\lvec|} \\
	&= \sum_{i=1}^d \sum_{w\in V^*} \sum_{\substack{\lvec=(l_1,\ldots,l_{i-1},l_{i+1},\ldots,l_d)\in \ZZnonneg^{d-1} \\ |\lvec|\leq k_i-1}} \mua_{\widehat\kvec^i+\lvec,\widehat\kvec^i} \\
	&\hspace{3cm}\cdot\Big(w\mid S(v_{k_1+l_1}\cdots v_{k_{i-1}+l_{i-1}})\shuffle v_{k_{i+1}+l_{i+1}}\cdots v_{k_d+l_d}\Big)\ w\otimes v_{k_i-|\lvec|},
    \end{align*}
	which immediately implies the claimed formula.
    \end{proof} 

	Having a formula for the reduced triangle map, we can use Theorem \ref{thm:free_coprod2} to give an explicit formula for the coproduct $\Delta_a$ dual to the Grossman-Larson product $\glp_a$.
	
	\begin{Theorem} \label{thm:coprod_ari} For $\kvec=(k_1,\ldots,k_d)\in \ZZnonneg^d$, we have
	\begin{align*}
	\Delta_a(v_{k_1}\cdots v_{k_d})&= \sum_{\substack{0\leq n\leq d \\ 0\leq i_1<j_1\leq i_2<j_2 \leq \cdots \leq i_n<j_n\leq d}} \sum_{\substack{\lvec=(\lvec_1,\ldots,\lvec_n) \\ \lvec_s=(l_{i_s+1},\ldots,\widehat{l_{j_s}},\ldots,l_{i_{s+1}}) \text{ for } s=1,\ldots,d \\
    \lvec_s=\emptyset \text{ or } |\lvec_s|\leq k_{j_s}-1}} \hspace{-0,1cm}\mua_{\widehat\kvec^{j_1,\ldots,j_n}+\lvec,\widehat\kvec^{j_1,\ldots,j_n}} \\
    &\Big(v_{k_1}\cdots v_{k_{i_1}}\shuffle S(v_{k_{i_1+1}+l_{i_1+1}}\cdots v_{k_{j_1-1}+l_{j_1-1}})\shuffle v_{k_{j_1+1}+l_{j_1+1}}\cdots v_{k_{i_2}+l_{i_2}}\shuffle\\
    &\hspace{3,5cm}\cdots \shuffle v_{k_{j_n+1}+l_{j_n+1}}\cdots v_{k_d+l_d}\Big)\otimes v_{k_{j_1}-|\lvec_1|}\cdots v_{k_{j_n}-|\lvec_n|},
	\end{align*}
    where we formally set $i_{n+1}=d$.
	\end{Theorem}

    In the formula, we used Notation \ref{nota:indices}. By the discussion after Proposition \ref{prop:red_tr_ari}, a $v_0$ at the $s$-th position in the right factor corresponds to $i_s+1=i_{s+1}$ and $\lvec_s=\emptyset$.

    \begin{proof}
    From Theorem \ref{thm:free_coprod2} and Proposition \ref{prop:red_tr_ari}, we deduce
    \begin{align*}
    &\Delta_a(v_{k_1}\cdots v_{k_d})\\
    &=\sum_{\substack{0\leq n\leq d \\ 0\leq i_1<\cdots <i_n<d}} \big(v_{k_1}\cdots v_{k_{i_1}}\big)\shuffle_\bullet \cotr_a^{\operatorname{irr}}\big(v_{k_{i_1+1}}\cdots v_{k_{i_2}}\big)\shuffle_\bullet\cdots \shuffle_\bullet \cotr_a^{\operatorname{irr}}\big(v_{k_{i_n+1}}\cdots v_{k_d}\big) \\
    &=\sum_{\substack{0\leq n\leq d \\ 0\leq i_1<\cdots <i_n<d}} \sum_{j_1=i_1+1}^{i_2}\cdots \sum_{j_n=i_n+1}^d \sum_{\substack{\lvec_s=(l_{i_s+1},\ldots,\widehat l_{j_s},\ldots,l_{i_{s+1}}),\ s=1,\ldots,n \\
    |\lvec_s|\leq k_{j_s}-1}} \ \ \prod_{s=1}^n \mua_{\widehat\kvec_s^{j_s}+\lvec_s,\widehat\kvec_s^{j_s}} \\
    & \Big(v_{k_1}\cdots v_{k_{i_1}}\otimes \one\Big)\shuffle_\bullet
    \Big(\big(S(v_{k_{i_1+1}+l_{i_1+1}}\cdots v_{k_{j_1-1}+l_{j_1-1}})\shuffle v_{k_{j_1+1}+l_{j_1+1}}\cdots v_{k_{i_2}+l_{i_2}}\big)\otimes v_{k_{j_1}-|\lvec_1|}\Big)\\
    &\shuffle_\bullet\cdots \shuffle_\bullet \Big(\big(S(v_{k_{i_n+1}+l_{i_n+1}}\cdots v_{k_{j_n-1}+l_{j_n-1}})\shuffle v_{k_{j_n+1}+l_{j_n+1}}\cdots v_{k_d+l_d}\big)\otimes v_{k_{j_n}-|\lvec_n|}\Big),
    \end{align*}
    where we denote $\kvec_s=(k_{i_s+1},\ldots,k_{i_{s+1}})$ for $s=1,\ldots,n$ and any decomposition $0\leq i_1<\cdots<i_n<d$.
    By applying the definition of $\shuffle_\bullet$ given in \eqref{eq:sh_bull_free}, we obtain the claimed formula.
    \end{proof}

\begin{Remark}
There seems not to be a formula for the coproduct $\Delta_a$ in the shape of Theorem $\ref{thm:coprod_similar_Gon}$ using the $I$-notation. Applying the $I$-notation means to eliminate the sums running over $i_1,\ldots i_n$ in Theorem \ref{thm:coprod_ari}. But here, the decompositions $\lvec=(\lvec_1,\ldots,\lvec_n)$ depend on the choice $i_1,\ldots,i_n$ and therefore we cannot apply the antipode relation.
\end{Remark}

Though the coproduct $\Delta_a$ corresponding to the ari bracket is not in the shape of the Goncharov coproduct, it still possesses some similar properties.

For any word $w=v_{k_1}...v_{k_d} \in V^*$, set $\indmax(w)=\max\limits_{1 \le s \le d}(k_s)$. This notion defines an increasing as well as a decreasing filtration on $\QV$ via
    \begin{align} \label{eq:fil_max}
    \operatorname{Fil}^n\QV&=\operatorname{span}\{w\in V^*\mid \indmax(w)\leq n\}, \\
    \operatorname{Fil}_n\QV&=\operatorname{span}\{w\in V^*\mid \indmax(w)\geq n\}.
    \end{align}
	
	\begin{Proposition} \label{prop:glpa_123}
	(i)	 For all $A\in \QV$ and $w\in V*$, we have 
		\[ A \glp_a w \in \operatorname{Fil}_{\indmax(w)} \QV.\]
       (ii)  For $w\in \QV$ write $\Delta_a(w)=w_{(1)}\otimes w_{(2)}$, then we always have  
        \[
        w_{(2)}\in \operatorname{Fil}^{\indmax(w)}\QV.
     \] 
	\end{Proposition} 
	
	\begin{proof} 
     Since $\mua_{\kvec,\lvec}$ is non-zero only for $|\kvec| - |\lvec|\ge 0$ the claim (i) follows from the explicit formula for $A \glp_a w$ given in   Proposition \ref{prop:glp_ari_formula}. Now (ii) follows by duality. 
	\end{proof}


    Brown \cite{br} proved that multiple zeta values $\zeta(k_1,\ldots,k_d)$ with $k_i\in\{2,3\}$ span the algebra $\Z$ of all multiple zeta values. A key ingredient in the proof is that one factor of the Goncharov coproduct preserves the subalgebra $\Q\langle x_0x_1,x_0x_0x_1\rangle$.
    Propositon \ref{prop:glpa_123} implies that if $w\in \Q\langle v_1,v_2,v_3\rangle$, then the right tensor product factor of $\Delta_a(w)$ is also contained in $w\in \Q\langle v_1,v_2,v_3\rangle$. This is a first indication that there might be a generalization of Brown's result. 
	
\begin{Example}  
We have
\begin{align*}
\Delta_a(v_2v_2v_1v_3)
&= -v_1v_3 \otimes v_2v_2 - v_2v_1v_3 \otimes  v_2 + v_2v_2 \otimes v_2v_2 + v_2v_2v_2 \otimes v_2\\ 
&\hspace{0.4cm} - v_2v_2 \otimes v_1v_3 + 2 v_1v_3 \otimes v_1v_3 + v_2 \otimes v_2v_2v_2 - v_2 \otimes v_2v_1v_3.
\end{align*}
\end{Example}
By observation we found the following conjecture, which we checked then numerically for a large number of words.
\begin{Conjecture} \label{conj:glpa_123} For $k\ge 0$ we define the $k$-level of a word  $v_{k_1}\cdots v_{k_d}$ by
\[
\operatorname{k-level}(v_{k_1}\cdots v_{k_d})= d- \# \{j\mid k_j = k\},
\]
and extend this to an increasing filtration on $\QV$. For $w\in \QV$, write as usual $\Delta_a(w)=w_{(1)}\otimes w_{(2)}$. Then, we always have  
\[
  \operatorname{k-level}(w_{(2)}) \le
  \operatorname{k-level}(w).
     \] 
\end{Conjecture} 
For words in a single letter $v_k$,  this conjecture predicts that the right tensor product factors are again words only in $v_k$. Indeed, this consequence follows directly from the explicit  description  $\tr_a$ in Proposition \ref{prop:mult_tr_ari}.  Moreover, if $w\in \Q\langle v_1,v_2,v_3\rangle$, then conjecturally the right tensor product factors of $\Delta_a(w)$ have not more 
indicies distinct to $2$ than  $w$. This is a another indication that there might be a generalization of Brown's result \cite{br} on a generating set for multiple zeta values.

By Theorem \ref{thm_ari_is_post_lie} 
the ari bracket belongs to a post-Lie algebra,  this yields the following 
an explicit version of \ref{thm:dual_pair_free_lie}.	

\begin{Theorem}\label{thm:dual_pair_ari}
Let $V=\{v_0,v_1,\ldots\}$ be equipped with the grading $\operatorname{wt}(v_0)=1$ and $\operatorname{wt}(v_i)=i$ for $i>0$. 
Consider the graded post-Lie algebra $\big(\Lie(V),\emptybrackets,\tr_a \big)$  associated with the ari bracket. Then, we get a dual pair of graded Hopf algebras 
    \[
    \xymatrix{ 
     (\Q\langle V\rangle,\shuffle,\Delta_a ) 
   \ar@{<~>}[rrr]^{\text{graded dual}}&&&    
   (\Q\langle V \rangle,\glp_a,\co). 
    }
   \]
An explicit and effective formula for the Grossman-Larson product $\glp_a$ is given by Proposition \ref{prop:glp_ari_formula}, and for its dual coproduct $\Delta_a$ it is given by Theorem \ref{thm:coprod_ari}.   
\graphref{thm_ari_is_post_lie} 
\graphref{thm:dual_pair_free_lie}  
\end{Theorem}
\noproof{Theorem}
We further have the following functorialities.

\begin{Lemma}
\label{lem:sublie_ortho_ari}  Let $i \ge 0$, then the space 
\[
\mathfrak{o}_a(v_i)  = \big\{ w \in \Lie(V\setminus \{v_0\} ) \, \big| \, w \tr_a v_i =0 \big\}
\]
is a Lie subalgebra  of $( \Lie(V),\emptybrace_a)$.
\end{Lemma}

\begin{proof}  One of the identities in Remark \ref{rem:post_lie_structures}  implies that for any $w,w' \in \mathfrak{o}_a(v_i) $ we get
\[
\{ w, w'\}_a \tr_a v_i = w \tr_a ( w' \tr_a v_i ) - w' \tr_a ( w \tr_a v_i )=0.
\] 
From Proposition \ref{prop:glp_ari_formula} it follows $\{ w, w'\}_a \in 
\Lie(V\setminus \{v_0\} )$.
\end{proof}  

\begin{Proposition} \label{prop:igeq1_sub_glpa}
	(i)	The tuple $\big( \Q\langle v_1,v_2,v_3,... \rangle , \glp_a, \co \big)$ is a Hopf subalgebra of \\
	$\big( \QV, \glp_a, \co \big) $.
    \vspace{0,3cm} \\
    (ii) The tuple $(\Q\langle v_1,v_2,v_3,\ldots\rangle,\shuffle,\Delta_a)$ is a Hopf subalgebra of $(\QV,\shuffle,\Delta_a)$.   
\end{Proposition}
	
\begin{proof}
	The result in (i) follows from Lemma \ref{lem:sublie_ortho_ari} applied with $v_i=v_0$. The second claim follows by duality.       
\graphref{thm:dual_pair_ari}
\end{proof}
	
\begin{Proposition} \label{prop:glpI_into_glpa}
	(i)	We have an embedding of Hopf algebras 
		\[
		\big( \Q\langle v_0,v_1 \rangle , \glp_I, \co \big)  \hookrightarrow  \big( \QV, \glp_a, \co \big).
		\] 
		(ii) The canonical projection $\QV\rightarrow \Q\langle v_0,v_1\rangle$ induces a surjective Hopf algebra morphism
		\[
		\big( \QV, \shuffle, \Delta_a \big)  \twoheadrightarrow \big( \Q\langle v_0,v_1 \rangle , \shuffle, \Delta_I \big).
		\] 
\end{Proposition}
	
\begin{proof}
	Let $A(\kvec),w\in \Q\langle v_0,v_1\rangle$. Then the only summand in Proposition \ref{prop:glp_ari_formula} corresponds to $\lvec_j=\kvec_j$ for $j=2,\ldots,2d+1$. Thus, we get for $w=w_1 v_{i_1}\cdots w_d v_{i_d} w_{d+1}$ as above 
	\[ A(\kvec)\glp_a w= A(\kvec)_{(1)} w_1S(A(\kvec)_{(2)}) v_{i_1} A(\kvec)_{(3)} w_2 \cdots w_d S(A(\kvec)_{(2d)}) v_{i_d} A(\kvec)_{(2d+1)}.\]
	This agrees with the formula for $\glp_I$ given in Proposition \ref{prop:glp_ihara_grouplike-formula}. \graphref{thm:dual_pair_ari}
\end{proof}

\begin{Proposition}\label{prop:v1_ortho_ari}
For  all  $w \in \mathfrak{o}_a(v_1)$  and   for all $b  \in \Lie(v_0,v_1)$ 
 we have 
\[
\{w,b\}_a =0.
\]
\end{Proposition}

\begin{proof}  
Assume $w  \in \mathfrak{o}_a(v_1)$, then for all $b=b_1\cdots b_n$, $b_i\in\{ v_0, v_1 \}$, we get  by the assumption $ w \tr_a b_i=0$ in both cases and therefore
\begin{align*}
w \glp_a b &= wb + w \tr_a b \\
&= wb + \sum_i  b_1 \ldots ( w \tr_a b_i ) \ldots b_n \\
&= wb.
\end{align*}
If $b \in \Lie(v_0,v_1)$, then for all $w=w_1\cdots w_m$, $w_i \in V\backslash \{v_0\}$, we get by 
Definition \ref{def:ari_triangle} of the triangle map $\tr_a$
\begin{align*}
b \glp_a w &= bw + b \tr_a w \\
&= bw + \sum_i  w_1 \ldots ( b \tr_a w_i ) \ldots w_m \\
&= bw + \sum_i  w_1 \ldots (w_i b-  b w_i ) \ldots w_m \\
&= bw- bw + w_1 b w_2 \ldots w_m - w_1 b w_2 \ldots w_m  + \ldots + wb\\
&=wb.
\end{align*}
Therefore,
\[
\{w,b\}_a = w \glp_a b - b \glp_a w  = 0.
\]
 \end{proof}
Using the isomorphism from Theorem \ref{thm:glp_hopf_algebra}, we obtain the following consequence.

\begin{Corollary}
Consider $ \mathcal{U}( \mathfrak{o}_a(v_1)) $ as a Hopf subalgebra of $\big( \QV, \glp_a, \co \big)$, then
for any $w \in \mathcal{U}( \mathfrak{o}_a(v_1)) $ and any $b \in   \Q\langle v_0,v_1 \rangle $, we
have
\[
w \glp_a b = b \glp_a w.
\]
\end{Corollary}

\begin{Remark}\label{rem:period_relations}
A simple calculations shows that $v_{2i+1} \in \mathfrak{o}_a(v_1)$ 
and $v_{2i} \notin \mathfrak{o}_a(v_1)$ for all $i \ge 1$.  

For $k,l\in\mathbb{N}$, we have
the explicit formula for the ari bracket\footnote{The same coefficients occur if we consider $\{ \operatorname{ad}(v_0)^{k_1-1}(v_1) ,
\operatorname{ad}(v_0)^{k_2-1}(v_1)\}_I$, cf \cite[eq. (28)]{LNT}} 
\begin{align} \label{eq:ari_two_lett}
\{v_{k_1},v_{k_2} \}_a &= [v_{k_1},v_{k_2}] + v_{k_1} \tr_a v_{k_2} - v_{k_2} \tr_a v_{k_1} \nonumber\\
&=  [v_{k_1},v_{k_2}]  +
\sum_{l_1=1}^{k_1} (-1)^{k_1+l_1} \binom{k_1-1}{l_1-1} [v_{k_2+k_1-l_1},v_{l_1}] \\
&\hspace{2.1cm}- \sum_{l_2=1}^{k_2} (-1)^{k_2+l_2} \binom{k_2-1}{l_2-1} [v_{k_1+k_2-l_2},v_{l_2}]. \nonumber
\end{align}
Consider the map $\rho_{C^{\operatorname{bi}}}$ from \eqref{eq:def_rho_C} restricted to $V\backslash\{v_0\}$,
\begin{align*}
\rho_{C^{\operatorname{bi}}}|_{V\backslash\{v_0\}}:\Q \langle V \setminus \{v_0\} \rangle & \to  \Q [X_1,X_2, \ldots] \\
v_{k_1}\ldots v_{k_d} &\mapsto X_1^{k_1-1}\cdots X_d^{k_d-1}.  
\end{align*}
Then due to \eqref{eq:ari_two_lett}, the assignment $[v_{k_1},v_{k_2}] \mapsto \{v_{k_1},v_{k_2}\}_a$ corresponds under $\rho_{C^{\operatorname{bi}}}|_{V\backslash\{v_0\}}$ to 
\[p(X_1,X_2) \mapsto p(X_1, X_2 ) + p(X_2-X_1,X_1) - p(X_1,X_2-X_1).\]
Restricting to anti-symmetric polynomials, the kernel of this map equals the space of even period polynomials, see for example \cite{brown-depth}.
Therefore, in $\mathfrak{o}_a(v_1)$ the well-established period polynomial relations in depth $2$ hold. For example  
in weight $12$, we get a relation that was first observed by Ihara
\[\{v_5,v_7\}_a - \frac{1}{3} \{ v_3,v_9\}_a =0.\]
\end{Remark}

\section{The uri bracket} \label{sec:uri}
	
 Schneps and the second author \cite{K-montreal}  
 proposed with the uri bracket a conjectural Lie structure on the space of alternil bimoulds.	An alternative description via non-commutative polynomials is discussed in the thesis of the first author \cite{BPhD}.  In this section, we rephrase the second approach in the alphabet $V=\{v_0,v_1,v_2,\ldots\}$ and show that we get a post-Lie algebra under the assumption of some combinatorial identities.

We expect a Hopf algebra structure on the algebra $\Z_q$ of multiple $q$-zeta values modulo quasi-modular forms \cite{BK20}. In particular the Lie bracket on the dual of the space of indecomposables 
should be the uri bracket.

\subsection{The uri multiplicities} \label{subsec:uri_mult}

In the following, we use Notation \ref{nota:indices} and Definition \ref{def:ari_mult}. For $n,s\in \mathbb{N}$, we denote the set of compositions of $n$ in $s$ parts by
\begin{align*}
 \mathcal{C}^s(n) &= \{ \boldsymbol{\alpha} \in \mathbb{N}^s \mid |\boldsymbol{\alpha}| = n \},
\end{align*}
and write $\mathcal{C}= \bigcup\limits_{s,n\in \mathbb{N}} \mathcal{C}^s(n)$.
    
\begin{Definition}
For $a\in \mathbb{N}$, $\boldsymbol{\alpha}=(\alpha_1,\ldots,\alpha_s)\in \mathbb{N}^s$, the \emph{threshold indicator function} is given by
\begin{align*}
 \ind{a}{\boldsymbol{\alpha}}&=\min_{1\leq j\leq s} \left\{ j \mid \alpha_1 + \dots + \alpha_j \geq a \right\}.
\end{align*}      
By convention we set $\ind{a}{\boldsymbol{\alpha}}=0$, if either $\boldsymbol{\alpha} = \emptyset$, $a > |\boldsymbol{\alpha}|$.
\end{Definition}

Consider the power series
\begin{align*}
	\mathbf{B}_1(x,y)&= \frac{y^2 x e^{x}}{(1-y)(e^{xy}-1)}.
\end{align*}
Write $\mathbf{B}_1(x,y)= \sum_{m,n\geq0 } B_1 (m,n) x^m y^n$, then we have
\begin{align} \label{eq:def_B1}
B_1(m,n)= \frac{1}{m!} \sum_{k=0}^{n-1} \binom{m}{k} B_ k,\qquad m,n\geq0,
\end{align}
where $B_k$ in the $k$-th Bernoulli-number.

We have for all $r\geq 1$ that $B_1(r,1) =  \frac{1}{r!}$, and for all $r\geq2$ that $B_1(r,r)= 0$. Moreover, for $m\ge 2$ we get the identity
	\begin{align*}
		B_1(m,n)&= -  \sum_{k=n}^{m-1} \frac{B_k}{k! (m-k)!}.
	\end{align*}

\begin{Remark}	With Maple we got the following table for $B_1(m,n)$, $m,n\leq 10$,
	
	\[
	\left[\begin{matrix} 
		1 &  
		\\
		\frac{1}{2} & 0  
		\\
		\frac{1}{6} & -\frac{1}{12} & 0 & 
		\\
		\frac{1}{24} & -\frac{1}{24} & 0 & 0 &  \\
		\frac{1}{120} & -\frac{1}{80} & \frac{1}{720} & \frac{1}{720} & 0 & 
		\\
		\frac{1}{720} & -\frac{1}{360} & \frac{1}{1440} & \frac{1}{1440} & 0 & 0 &   
		\\
		\frac{1}{5040} & -\frac{1}{2016} & \frac{1}{5040} & \frac{1}{5040} & -\frac{1}{30240} & -\frac{1}{30240} & 0 &  
		\\
		\frac{1}{40320} & -\frac{1}{13440} & \frac{1}{24192} & \frac{1}{24192} & -\frac{1}{60480} & -\frac{1}{60480} & 0 & 0 &   
		\\
		\frac{1}{362880} & -\frac{1}{103680} & \frac{1}{145152} & \frac{1}{145152} & -\frac{17}{3628800} & -\frac{17}{3628800} & \frac{1}{1209600} & \frac{1}{1209600} & 0  
		\\
		\frac{1}{3628800} & -\frac{1}{907200} & \frac{1}{1036800} & \frac{1}{1036800} & -\frac{1}{1036800} & -\frac{1}{1036800} & \frac{1}{2419200} & \frac{1}{2419200} & 0 & 0 
	\end{matrix}\right].
	\]	
\end{Remark}

\begin{Definition} \label{def:uri_mult} The \emph{uri multiplicities} are given by
\begin{align*}
\mu_{a, \boldsymbol{\alpha}} = B_1\left(\ell(\boldsymbol{\alpha}), \ind{a}{\boldsymbol{\alpha}}\right), \qquad a\in \mathbb{N}, \boldsymbol{\alpha}\in \mathcal{C}.  
\end{align*}
Here, we set $\mu_{a, \boldsymbol{\alpha}}=0$ if $ \ind{a}{\boldsymbol{\alpha}}=0$.
\end{Definition}

\begin{Proposition} \label{prop:uri_descend_identity} For all
 $\mathbf{d} \in \mathbb{N}^d$, $a \in \mathbb{N}$, and $\mathbf{b} \in \mathbb{N}^s$, we have
\[
\sum_{\mathbf{0} \leq \mathbf{r} \leq \mathbf{d}} (-1)^{|\mathbf{d}| + |\mathbf{r}|} \mathbf{m}_{\mathbf{d}  + 1, \mathbf{r} + 1} \sum_{\boldsymbol{\alpha} \in \mathcal{C}^s(a + |\mathbf{d} - \mathbf{r}|)} \mathbf{m}_{\boldsymbol{\alpha}, \mathbf{b}} \, \mu_{a, \boldsymbol{\alpha}} = \mathbf{m}_{a, |\mathbf{b}| - |\mathbf{d} |} \,\mu_{|\mathbf{b}| - |\mathbf{d} |, \mathbf{b}}.
\]
\end{Proposition}

\begin{proof}
Let $1\leq j \leq s$. The coefficient of $B_1(s,j)$ defined in \eqref{eq:def_B1} in the left hand side of the claimed formula is given by
\begin{align*}
c(s,j) &= \sum_{\mathbf{0} \leq \mathbf{r} \leq \mathbf{d}} (-1)^{|\mathbf{d}| + |\mathbf{r}|} \, \mathbf{m}_{\mathbf{d} + \mathbf{1}, \, \mathbf{r}  + \mathbf{1}} \sum_{\substack{\boldsymbol{\alpha} \in \mathcal{C}^s(a + |\mathbf{d} - \mathbf{r}|) \\ \ind{a}{\boldsymbol{\alpha}}=j}} \mathbf{m}_{\boldsymbol{\alpha}, \mathbf{b}}\\
&=\sum_{\mathbf{0} \leq \mathbf{r} \leq \mathbf{d}} \binom{\mathbf{d}}{\mathbf{r}} \sum_{\substack{\boldsymbol{\alpha} \in \mathcal{C}^s(a + |\mathbf{d} - \mathbf{r}|) \\ \ind{a}{\boldsymbol{\alpha}}=j}} (-1)^{|\boldsymbol{\alpha}|+|\mathbf{b}|} \binom{\boldsymbol{\alpha}-1}{\mathbf{b}-1}.
\end{align*}
Thus, the statement is equivalent to
\begin{align*}
c(s,j) &= \begin{cases} \mathbf{m}_{a, |\mathbf{b}| - |\mathbf{d}|} ,& j= \ind{|\mathbf{b}| - |\mathbf{d}|}{\mathbf{b}}\\
 0, &\text{else}\\
 \end{cases}\\
 & = \begin{cases} (-1)^{a+|\mathbf{b}| - |\mathbf{d}|}\binom{a-1}{|\mathbf{b}| - |\mathbf{d}|-1} ,& j= \ind{|\mathbf{b}| - |\mathbf{d}|}{\mathbf{b}}\\
 0, &\text{else},
 \end{cases}
\end{align*}
which is exactly the statement of Theorem \ref{thm:conj_uri_prop}.
\end{proof}

\subsection{The uri bracket}

As in the previous section, we consider the alphabet $V=\{v_0,v_1,v_2,\ldots\}$.
We present here our third example for the 
magmatic approach towards a post-Lie structure on $\Lie(V)$ as 
described in Subsection \ref{subsec:magma_post_lie}.

\begin{Definition} \label{def:uri_triangle}
We define a map
\begin{align*}
&\tr_u:  M(V) \times V  \to M(V)_\Q\\
&(t(\kvec),v_a) \mapsto \begin{cases}    0,  & \hspace{-0.2cm} a=0, \\  \sum\limits_{\lvec\in \ZZnonneg^{\ell(\kvec)}} \, \mua_{\kvec,\lvec} \,
\sum\limits_{\boldsymbol{\alpha} \in \mathcal{C}(a+|\kvec-\lvec|)}  \,\muu_{a,\boldsymbol{\alpha} }\,\, 
v_{\alpha_1}\star(v_{\alpha_2}\star( \cdots \star(v_{\alpha_{\ell({\boldsymbol{\alpha}})}} \star t(\lvec))\cdots))  \quad & \hspace{-0.2cm}\text{else}.   \\ \end{cases}
		\end{align*}

Recall from \eqref{eq:ari_mult_nonzero}, that we have $\mua_{\kvec,\lvec}=0$ for $\lvec>\kvec$. Thus, the sum in the above definition of $\tr_u$ is finite. By Definition \ref{def:ari_mult}, we get non-zero multiplicities if  $k_i=0$ only for $l_i=0$ and if $k_i \ge 1$ only for $l_i \ge 1$.
        
If it is clear from the context, we will omit that $\lvec$ is an element in $\ZZnonneg^{\ell(\kvec)}$.
        
We denote by $\tr_u$ also its extensions to a derivation on $M(V)_\Q$ as in \eqref{eq:tr_ext_M}, and to $\Q\langle M(V)\rangle$ as in \eqref{eq:tr_ext_AM}.
\end{Definition}
	
\begin{Lemma}\label{lem:uri_descends}
The map $\tr_u$ on $M(V)_\Q$ descends to $\tr_u:\Lie(V)\otimes\Lie(V)\to \Lie(V)$.
\end{Lemma}   
	
\begin{proof} By Lemma \ref{lem:lie_descend}, it suffices to show that $I_{\Lie} \tr_u V   \subseteq I_{\Lie}$.  For $x,y,z\in M(V)_\Q$, and $a\geq1$, we have
\begin{align*}
&\big( x(\kvec)\star x(\kvec) \big) \tr_u v_a \\
&=\sum_{\lvec}  \mua_{\kvec,\lvec}^2 \sum_{\boldsymbol{\alpha}\in \mathcal{C}(a+2(\abs{\kvec-\lvec}))} \muu_{a,\boldsymbol{\alpha}} v_{\alpha_1}\star(v_{\alpha_2}\star(\cdots\star(v_{\alpha_{\ell(\boldsymbol{\alpha})}}\star (x(\lvec)\star x(\lvec)))\cdots ))\\
&+ \sum_{\lvec_1 <  \lvec_2 }  \mua_{\kvec,\lvec_1}\,\mua_{\kvec,\lvec_2}\, \sum_{\boldsymbol{\alpha}\in\mathcal{C}(s + |\kvec - \lvec_1|+ |\kvec - \lvec_2|)} \muu_{a,\boldsymbol{\alpha}} \\
&\hspace{3cm} \cdot v_{\alpha_1}\star(v_{\alpha_2}\star(\cdots\star(v_{\alpha_{\ell(\boldsymbol{\alpha})}}\star (x(\lvec_1)\star x(\lvec_2)+x(\lvec_2)\star x(\lvec_1)))\cdots )) \\ 
&\in I_{\Lie}.
\end{align*}
Next, set $J(\kvec_1,\kvec_2,\kvec_3) = J\big(x(\kvec_1),y(\kvec_2),z(\kvec_3) \big) $, then 
\begin{align*}
J(\kvec_1,\kvec_2,\kvec_3)  \tr_a v_s &=
\sum_{\lvec_1,\lvec_2,\lvec_3} \mua_{\kvec_1,\lvec_1}\,\mua_{\kvec_2,\lvec_2}\,\mua_{\kvec_3,\lvec_3}\,
\sum_{\boldsymbol{\alpha}\in \mathcal{C}(s + |\kvec_1 - \lvec_1|+ |\kvec_2 - \lvec_2|+ |\kvec_3 - \lvec_3|)}  \muu_{a,\boldsymbol{\alpha}} \\
& \hspace{2cm}\cdot v_{\alpha_1}\star(v_{\alpha_2}\star(\cdots\star(v_{\alpha_{\ell(\boldsymbol{\alpha})}}\star J(\lvec_1,\lvec_2,\lvec_3))\cdots ))\\
&\in I_{\Lie}. \qedhere
\end{align*}
\end{proof}

Recall from \eqref{eq:wt}, \eqref{eq:dep} that we have a weight and depth grading on the Lie algebra $\Lie(V)$. In particular, this also induces a descending filtration with respect to the depth.

\begin{Lemma}  \label{lem:ari_depthgraded_of_uri}
(i) The map $\tr_u$ is graded with respect to the weight and filtered with respect to the depth.

(ii) The associated depth-graded of 
$\tr_u$ equals $\tr_a$.
\end{Lemma}

\begin{proof}
The first claim follows directly from the Definition \ref{def:uri_triangle} of $\tr_u$. To obtain the associated depth-graded map of $\tr_u$, we have to omit all terms in $t(\kvec)\tr_u v_a$ which are of depth $\dep(\kvec)+1$. These are all terms, where we choose the composition $\boldsymbol{\alpha}$ in Definition \ref{def:uri_triangle} to be of length $>1$. The only composition of length $1$ is $\boldsymbol{\alpha}=(a+\abs{\kvec-\lvec})$ and in this case we get $\mu_{\boldsymbol{\alpha},a}=1$. Thus, the associated depth-graded of $\tr_u$ equals the formula for $\tr_a$ given in Definition \ref{def:ari_triangle}. 

\end{proof}

As indicated in Remark \ref{rem:post_lie_recursion} we have the following efficient way to calculate $\tr_u$. 

\begin{Lemma}\label{lem:uri_magma_boxed} \graphref{rem:post_lie_recursion} 
We have for $a\geq1$
\begin{align*}
	\big(&t_1(\kvec_1)\cdot t_2(\kvec_2)\big)  \tr_u v_a =  
	\sum_{\lvec_1,\lvec_2} \,m_{\kvec_1,\lvec_1}m_{\kvec_2,\lvec_2}\,\,
	\sum\limits_{\boldsymbol{\alpha} \in \mathcal{C}(a+|\kvec_2|-|\lvec_2|)}   \,\muu_{a,\boldsymbol{\alpha}}  \sum_{i=1}^{\ell(\boldsymbol{\alpha})} 
	\\
	&\sum_{\boldsymbol{\beta} \in \mathcal{C}( \alpha_i +|\kvec_1-\lvec_1|) }  \,\muu_{\alpha_i,\boldsymbol{\beta}} v_{\alpha_1}\star \big(\cdots \star\big( v_{\alpha_{i-1}}\star \big(\,
	\underbrace{ (v_{\beta_1}\star(v_{\beta_2}\star( \cdots \star(v_{\beta_{\ell(\boldsymbol{\beta})}} \star t_1(\lvec_1)) \cdots)))}_{\text{$i$-th position}} \\
	&\hspace{7.3cm}\star \big(v_{\alpha_{i+1}}
	\cdots\star\big(v_{\alpha_{\ell(\boldsymbol{\alpha})}} \star t_2(\lvec_2)\big)\cdots\big)\big)\big)\cdots \big)
\end{align*}
\end{Lemma}	
	
\begin{proof}  
We have by definition of $\tr_u$  
\begin{align*}
&\big( t_1(\kvec_1)   \cdot t_2(\kvec_2)  \big) \tr_u v_a \\
&= t_1(\kvec_1) \tr_u \Big(   t_2(\kvec_2) \tr_u v_a  \Big) - \Big(  t_1(\kvec_1) \tr_u  t_2(\kvec_2)\Big) \tr_u v_a\\	
& =  \sum\limits_{\lvec_2  } \, \mua_{\kvec_2,\lvec_2} \, \sum\limits_{\boldsymbol{\alpha} \in \mathcal{C}(a+|\kvec_2-\lvec_2|)}  \,\muu_{a,\boldsymbol{\alpha} }\,\, 
t_1(\kvec_1) \tr_u \Big( v_{\alpha_1}\star(v_{\alpha_2}\star( \cdots\star(v_{\alpha_{\ell(\boldsymbol{\alpha})}} \star t_2(\lvec_2)) \cdots)) \Big)\\
& \hspace{0.4cm}- \Big(  t_1(\kvec_1) \tr_u  t_2(\kvec_2) \Big) \tr_u v_a\\	
& =  \sum\limits_{\lvec_2} \, \mua_{\kvec_2,\lvec_2} \,
\sum\limits_{\boldsymbol{\alpha} \in \mathcal{C}(a+|\kvec_2-\lvec_2|)}  
\,\muu_{a,\boldsymbol{\alpha}} \sum_{i=1}^{\ell(\boldsymbol{\alpha})}\,\, v_{\alpha_1}\star(v_{\alpha_2}\star(\cdots  \star((t_1(\kvec_1) \tr_u v_{\alpha_i}) \\
&\hspace{8cm}\star( \cdots \star(v_{\alpha_{\ell(\boldsymbol{\alpha})}} \star t_2(\lvec_2))\cdots )) \cdots ))\\
&\hspace{0.4cm}+\sum\limits_{\lvec_2 } \, \mua_{\kvec_2,\lvec_2} \,
\sum\limits_{\boldsymbol{\alpha} \in \mathcal{C}(a+|\kvec_2-\lvec_2|)}  \,\muu_{a,\boldsymbol{\alpha} }\,\, 
 v_{\alpha_1}\star(v_{\alpha_2}\star( \cdots \star(v_{\alpha_{\ell(\boldsymbol{\alpha})}} \star (t_1(\kvec_1) \tr_u t_2(\lvec_2)))\cdots) )		\\
&\hspace{0.4cm} - \Big(  t_1(\kvec_1) \tr_u  t_2(\kvec_2) \Big) \tr_u v_a 
\end{align*}
Now by  Lemma \ref{lem:uri_lemma_2} below the  claimed identity follows.
\end{proof}

\begin{Remark}\label{rem:uri_tr_not_efective} We were not successful in deriving a general formula. However, we expect that Lemma \ref{lem:uri_magma_boxed} generalizes to $n$ factors \[(t_1(\kvec_1)\cdots t_n(\kvec_n))\tr_u v_a.\]  The right hand side of such a formula is then made in an iterative manner extending Lemma \ref{lem:ari_magma_boxed}:  In  the $r$-th step we replace each $v_{\alpha_i}$ in $(t_{n-r+2}(\kvec_{n-r+2})\cdots t_n(\kvec_n))\tr_u v_a$ coming from a composition corresponding to some $t_{n-s+1}(\kvec_{n-s+1})$, $s<r$, by $ t_{n-r+1}(\kvec_{n-r+1}) \tr_u v_{\alpha_i}$. If we picture this in trees as explained in Subsection \ref{subsec:ari_lie}, this means we replace only every left leave $v_{\alpha_i}$ of $(t_{n-r+2}(\kvec_{n-r+2})\cdots t_n(\kvec_n))\tr_u v_a$ by $t_{n-r+1}(\kvec_{n-r+1}) \tr_u v_{\alpha_i}$. We checked this generalized formula up to degree $15$.
	\end{Remark}
     
\makeinvisible{keep   ... , then we have a conjectural expression, which had been proven experimentally upto weight $15$. The conjectural right hand side is made in an iterative manner:  In  the $r$-th step we have to sum over all terms, where in the formula of the previous step one of the "$v_i$'s" has been replaced by $ t_{n-r+1}(\lvec_{n-r+1}) \tr_u v_i$, with the appropriate multiplicity. Here we mean by the  "$v_i$'s" exactly those factors who are left to a "$\star$" or in the first step it is just $v_a$. }

	\begin{Lemma}\label{lem:uri_lemma_2} 
   We have for $a\geq1$ the identity
		\begin{align} \eqlabel{eq:uri_double_2}
			&\Big(  t_1(\kvec_1) \tr_u  t_2(\kvec_2)     \Big) \tr_u v_a  \\
			&=   \sum\limits_{\lvec_2 } \, \mua_{\kvec_2,\lvec_2} \,
			\sum\limits_{\boldsymbol{\alpha} \in \mathcal{C}(a+|\kvec_2-\lvec_2|)}  \,\muu_{a,\boldsymbol{\alpha} }\,\, 
			v_{\alpha_1}\star(v_{\alpha_2}\star( \cdots \star(v_{\alpha_{\ell(\boldsymbol{\alpha})}} \star (t_1(\kvec_1) \tr_u t_2(\lvec_2)))\cdots) ).	 \nonumber
		\end{align}
	\end{Lemma}
	\begin{proof} Similar to Notation \ref{not:ari_a}, we define for $t_1(\kvec),t_2(\kvec)\in M(V)$ with $\kvec\in \ZZnonneg^{d_1}$, $\mathbf{n}=(n_1,\ldots,n_{d_2})\in \ZZnonneg^{d_2}$, $j\in \{1,\ldots,d_2\}$, and $\boldsymbol{\alpha}\in\mathcal{C}^s(n_j)$
		\[	u^{t_1,t_2}_j ( \lvec ,\mathbf{n},\boldsymbol{\alpha}) =
			t_2( v_{n_1}, \ldots, v_{n_{j-1}}, 
		( v_{\alpha_1}\star(v_{\alpha_2}\star( \cdots\star(v_{\alpha_s} \star   t_1(\lvec))\cdots )) )		 
		, v_{n_{j+1}},\ldots, v_{n_{d_2}} )\] 
	if $n_j>0$ and $u^{t_1,t_2}_j ( \lvec,\mathbf{n},\boldsymbol{\alpha} )=0$ else.
		Then, we observe that  
		\begin{align*}
			t_1(\kvec) \tr_u  t_2( {\bf n} ) &=  
			\sum_{j=1}^{d_2}  \sum_{\lvec}  \sum\limits_{\boldsymbol{\alpha} \in \mathcal{C}(n_j+|\kvec-\lvec|)}  
			\,  \mua_{\kvec,\lvec}  \, \muu_{n_j,\boldsymbol{\alpha}}\,\,   
			\, u^{t_1,t_2}_j  ( \lvec, \mathbf{n},\boldsymbol{\alpha}).  
		\end{align*}
		Using this observation we get for the left hand side of \eqref{eq:uri_double_2} 
		\begin{align*}
			&\Big(  t_1(\kvec_1) \tr_u  t_2(\kvec_2)  \Big) \tr_u v_a =  
			\sum_{j=1}^{\ell(\kvec_2)}  \sum_{\lvec_1 }  \sum\limits_{\boldsymbol{\beta}\in \mathcal{C}(k_{2,j}+|\kvec_1-\lvec_1|)}  
			\,  \mua_{\kvec_1,\lvec_1}  \, \muu_{k_{2,j},\boldsymbol{\beta} }\,\,   
			\, u^{t_1,t_2}_j  ( \lvec_1, \kvec_2,\boldsymbol{\beta})    \tr_u v_a \\
			&=  
			\sum_{j=1}^{\ell(\kvec_2)}  \sum_{\lvec_1 }  \sum\limits_{\boldsymbol{\beta} \in \mathcal{C}(k_{2,j}+|\kvec_1-\lvec_1|)}  
			\,  \mua_{\kvec_1,\lvec_1}  \, \muu_{k_{2,j},\boldsymbol{\beta}}\, \sum_{{\bf n_1},{\bf n_2},\boldsymbol{\gamma}} 
			\mua_{\lvec_1, {\bf n}_1}  \mua_{\widehat{\kvec_2}^{ j}, {\bf n}_2} \mua_{ \boldsymbol{\beta},  \boldsymbol{\gamma}}  \\
			&\sum\limits_{\boldsymbol{\alpha} \in \mathcal{C}(a +|\lvec_1-{\bf n}_1 | +| \widehat{\kvec_2}^{ j} - {\bf n}_2| + |\boldsymbol{\beta} - \boldsymbol{\gamma}|)} 
			\muu_{a,\boldsymbol{\alpha}} v_{\alpha_1}\star(v_{\alpha_2}\star( \cdots \star(v_{\alpha_{\ell(\boldsymbol{\alpha}) }} \star   u^{t_1,t_2}_j({\bf n}_1, {\bf n}_2,\boldsymbol{\gamma} ))\cdots )),
		\end{align*}
		and for the right hand side of \eqref{eq:uri_double_2}  we obtain
		\begin{align*}
			&   \sum\limits_{\lvec_2  } 
			\sum\limits_{\boldsymbol{\alpha }\in \mathcal{C}(a+|\kvec_2-\lvec_2|)}  \, \mua_{\kvec_2,\lvec_2} \,\muu_{a,\boldsymbol{\alpha}}\, v_{\alpha_1}\star(v_{\alpha_2}\star( \cdots\star(v_{\alpha_{\ell(\boldsymbol{\alpha})}} \star \big( t_1(\kvec_1) \tr_u t_2(\lvec_2))) \cdots))	\\	 
			&=\sum\limits_{\lvec_2} 
			\sum\limits_{\boldsymbol{\alpha} \in \mathcal{C}(a+|\kvec_2-\lvec_2|)}  \, \mua_{\kvec_2,\lvec_2} \,\muu_{a,\boldsymbol{\alpha}}\,\, 
			\sum_{j=1}^{\ell(\kvec_2)}  \sum_{\lvec_1 }  \sum\limits_{\boldsymbol{\beta} \in \mathcal{C}(l_{2,j}+|\kvec_1-\lvec_1|)}  
			\,  \mua_{\kvec_1,\lvec_1}  \, \muu_{l_{2,j},\boldsymbol{\beta}}\,\,   
			\\ 
			&\hspace{1cm}  v_{\alpha_1}\star(v_{\alpha_2}\star( \cdots\star(v_{\alpha_{\ell(\boldsymbol{\alpha})}} \star u^{t_1,t_2}_j  ( \lvec_1, \lvec_2,\boldsymbol{\beta}))\cdots )).
		\end{align*}
		
		We compare both sides by computing the coefficient of 
		\[v_{\alpha_1}\star(v_{\alpha_2}\star( \cdots\star(v_{\alpha_{\ell(\boldsymbol{\alpha})}} \star u^{t_1,t_2}_j  ( {\bf n}_1, {\bf n}_2,\boldsymbol{\gamma}))\cdots )).\]
		For the left hand side we get with Lemma \ref{lem:ari_mult_klr} 
		\begin{align*}
			C_l&=   \muu_{a,\boldsymbol{\alpha}}\, \mua_{\widehat{\kvec_2}^{ j}, {\bf n}_2}\,\,   \sum_{\lvec_1 } \,  \mua_{\kvec_1,\lvec_1}\,\mua_{\lvec_1, {\bf n}_1}\,
			\sum\limits_{\boldsymbol{\beta} \in \mathcal{C}(k_{2,j}+|\kvec_1-\lvec_1|)}  
			\,   \mua_{\boldsymbol{\beta},\boldsymbol{\gamma}}  \, \muu_{k_{2,j},\boldsymbol{\beta} }  \\
			&=   \muu_{a,\boldsymbol{\alpha}}\, \mua_{\widehat{\kvec_2}^{ j}, {\bf n}_2}\,\,   
			\mua_{\kvec_1,{\bf n}_1}\, \sum_{\lvec_1 }  
			\, (-1)^{\abs{\kvec_1}+ \abs{\lvec_1} }   \mua_{\kvec_1-{\bf n}_1+\one, \lvec_1 - {\bf n}_1+\one }\,
			\sum\limits_{\boldsymbol{\beta} \in \mathcal{C}(k_{2,j}+|\kvec_1-\lvec_1|)}  
			\,   \mua_{\boldsymbol{\beta},\boldsymbol{\gamma}}  \, \muu_{k_{2,j},\boldsymbol{\beta} }.   
		\end{align*}
		For the right hand side we get with $\lvec_1 =  {\bf n}_1,\ \widehat{\lvec_2}^{ j} =  {\bf n}_2,\ l_{2,j} = |\boldsymbol{\gamma}| - |\kvec_1- {\bf n}_1 |,$ and $\boldsymbol{\beta}=\boldsymbol{\gamma}$ the coefficient
		\begin{align*}
			C_r&= \muu_{a,\boldsymbol{\alpha}} \, \mua_{\widehat{\kvec_2}^j, {\bf n}_2}\,
			\mua_{\kvec_1,{\bf n}_1}  \,
			\mua_{k_{2,j}, |\boldsymbol{\gamma}| - |\kvec_1- {\bf n}_1 |}\,
			\muu_{|\boldsymbol{\gamma}| - |\kvec_1- {\bf n}_1 |,\boldsymbol{\gamma}} .
		\end{align*}
		By Proposition \ref{prop:uri_descend_identity}, the multiplicities $C_l,C_r$ agree. Hence, the identity \eqref{eq:uri_double_2} holds.
	\end{proof}

	For the next result, we have to assume Conjecture \ref{conj:uri_mult_threshold_shuffle}, which states that the uri multiplicities satisfy the threshold shuffle identities from Definition \ref{def:threshold_shuffle}. In the following proof we will make use of the notations introduced in Appendix \ref{app:threshold_shuffle} in the context of the threshold shuffle identities.

	\begin{Theorem}\label{thm:uri_is_post_lie} 
   If Conjecture \ref{conj:uri_mult_threshold_shuffle} holds, then the triple $\big( \Lie(V), \emptybrackets, \tr_u\big) $ is a post-Lie algebra.
	\end{Theorem}

	\begin{proof} Because of Lemma \ref{lem:uri_descends} and Proposition \ref{prop:tr_descends} it suffices 
	 to show that for any $a\geq 1$, the following expression vanishes
		\begin{align*}\eqlabel{eq:ari_lie_identity_old}
			&\Lie\Big(( t_1(\kvec_1) \star t_2(\kvec_2)  ) \tr_u v_a   -( t_1(\kvec_1) \cdot t_2(\kvec_2)  - t_2(\kvec_2) \cdot t_1(\kvec_1) ) \tr_u v_a   \Big) \\
			&= \Lie\Big(( t_1(\kvec_1) \star t_2(\kvec_2)  ) \tr_u v_a  \Big)  -\Lie\Big( ( t_1(\kvec_1) \cdot t_2(\kvec_2))  \tr_u v_a \Big) \\
			&\hspace{9cm} +\Lie\Big( (t_2(\kvec_2) \cdot t_1(\kvec_1) ) \tr_u v_a   \Big).
		\end{align*}
		We  set
		$ \Lie(t_i(\lvec_i))= T_i(\lvec_i)$ for $i=1,2$. We will show that each of the above three terms is of the following type 
		\begin{align*}
			&\sum\limits_{\lvec_1,\lvec_2} \Bigg( \sum\limits_{\bsigma, \btheta} \Big(
			m_{\bsigma, \btheta}(\lvec_1,\lvec_2) \,\,v_\bsigma T_1(\lvec_1)T_2(\lvec_2) v_{S(\btheta)}  
			- m_{\bsigma, \btheta}(\lvec_2,\lvec_1) \,\,v_\bsigma T_2(\lvec_2)T_1(\lvec_1) v_{S(\btheta)} \Big)\\
			&+ \sum\limits_{\bsigma, \btheta} \Big(
			m_{\bsigma, \btau, \btheta}(\lvec_1,\lvec_2) \,\,v_\bsigma T_1(\lvec_1)v_\btau T_2(\lvec_2) v_{S(\btheta)}  
			- m_{\bsigma, \btau,\btheta}(\lvec_2,\lvec_1) \,\,v_\bsigma T_2(\lvec_2)v_\btau T_1(\lvec_1) v_{S(\btheta)} \Big) 
			\Bigg)
		\end{align*}
        for some indices $\bsigma,\btheta,\btau$. Here, for $\bsigma=(\sigma_1,\ldots,\sigma_n)$ we abbreviate 
        \[ v_\bsigma=v_{\sigma_1}\cdots v_{\sigma_n},\qquad v_{S(\bsigma)}=S(v_{\bsigma}).\]
		It suffices to study these multiplicities $m_*(*,*)$ for fixed  $\lvec_1$ and $\lvec_2$. We set $d_1= \abs{ \kvec_1 - \lvec_1}$ and $d_2=\abs{\kvec_2-\lvec_2}$.

		At first, we study $\Lie\Big(( t_1(\kvec_1) \star t_2(\kvec_2)  ) \tr_u v_a  \Big)$. 
		We have by definition of $\tr_u$ that
		\begin{align*}
		&( t_1(\kvec_1) \star t_2(\kvec_2)  ) \tr_u v_a \\
		&=  \sum\limits_{\lvec_1,\lvec_2} \,m_{\kvec_1,\lvec_1}m_{\kvec_2,\lvec_2}\,\,
		\sum\limits_{\boldsymbol{\alpha} \in \mathcal{C}(a+d_1+d_2)}  \,\muu_{a,\boldsymbol{\alpha}}\,\, 
		v_{\alpha_1}\star(v_{\alpha_2}\star( \cdots \star(v_{\alpha_{\ell(\boldsymbol{\alpha})}} \star ( t_1(\lvec_1) \star t_2(\lvec_2)) )\cdots)). 	
		\end{align*}
		Applying $\Lie$ to each summand yields with Lemma \ref{lem:boxed_brackets_antipode}
		\begin{align*}
		&\Lie \Big(  v_{\alpha_1}\star(v_{\alpha_2}\star( \cdots\star(v_{\alpha_{\ell(\alpha)}} \star ( t_1(\lvec_1) \star t_2(\lvec_2)))\cdots))  \Big) \\
		&= \big[  v_{\alpha_1},\big[ v_{\alpha_2},\big[ \cdots,
			\big[v_{\alpha_{\ell(\alpha)}} ,  T_1(\lvec_1) T_2(\lvec_2)  -  T_2(\lvec_2) T_1(\lvec_1) \big ]\cdots\big]\big] \big] \\
			&= v_{\balpha_{(1)} } T_1(\lvec_1)  T_2(\lvec_2) v_{S(\balpha_{(2)})}
			-   v_{\balpha_{(1)} } T_2(\lvec_2)  T_1(\lvec_1) v_{S(\balpha_{(2)})}.
		\end{align*}
        Here, we used the shorthand notation $\co(v_{\balpha})=v_{\balpha_{(1)}}\otimes v_{\balpha_{(2)}}$.
		We obtain for the first of our sought multiplicities 
		\begin{align*}
		m_{\bsigma, \btheta}(\lvec_1,\lvec_2)\big(  ( t_1 \star t_2  ) \tr_u v_a \big)&:=\Big( \Lie(( t_1(\kvec_1) \star t_2(\kvec_2)  ) \tr_u v_a) \big | v_{\bsigma} T_1(\lvec_1) T_2(\lvec_2) v_{\btheta} \Big)\\ 
		&=  m_{\kvec_1,\lvec_1}m_{\kvec_2,\lvec_2}\,\, \sum\limits_{ \boldsymbol{\alpha} \in \mathcal{C}(\bsigma, \btheta)} \,\muu_{a,\boldsymbol{\alpha}},
		\end{align*}
		where the sum runs over the set
		\begin{align*}
			\mathcal{C}( \bsigma,\btheta) 
			&=\Big \{ \balpha \in \mathcal{C}(a+d_1+d_2)\, : \,  \big( v_{\balpha_{(1)}} T_1 T_2 v_{S(\balpha_{(2)})} \, \big| 
			\, v_{\bsigma }   T_1 T_2 v_{S(\btheta)}  \big ) \neq 0 \,\Big\} \\
			& = 
			\Big\{ \balpha \in \mathcal{C}(a+d_1+d_2)\, : \,  \,  \big( \co( \balpha)\, \big | \, ( \sigma_1, \ldots,  \sigma_{\ell(\bsigma)} )\otimes 
			(  \vartheta_1, \ldots, \vartheta_{\ell(\btheta)} ) \big ) \neq 0 \,\Big\} \\
			&= 
			\Big\{ \balpha \in \mathcal{C}(a+d_1+d_2)\, : \,  \,  \big(  \balpha \, \big|  ( \sigma_1, \ldots,  \sigma_{\ell(\bsigma)} ) \shuffle
			( \vartheta_1, \ldots, \vartheta_{\ell(\btheta)} ) \big ) \neq 0 \,\Big\}.
		\end{align*}
		We observe that
		\[ m_{\bsigma, \btheta}(\lvec_2,\lvec_1)\big(  ( t_1 \star t_2  ) \tr_u v_a \big) = m_{\bsigma, \btheta}(\lvec_1,\lvec_2)\big(  ( t_1 \star t_2  ) \tr_u v_a \big). \]
		
		Now, we consider $\Lie\Big( ( t_1(\kvec_1) \cdot t_2(\kvec_2))  \tr_u v_a \Big)$. We get by Lemma \ref{lem:uri_magma_boxed}
	\begin{align*}
		\big(&t_1(\kvec_1)\cdot t_2(\kvec_2)\big)  \tr_u v_a =  
		\sum_{\lvec_1,\lvec_2} \,m_{\kvec_1,\lvec_1}m_{\kvec_2,\lvec_2}\,\,
		\sum\limits_{\boldsymbol{\alpha} \in \mathcal{C}(a+|\kvec_2|-|\lvec_2|)}   \,\muu_{a,\boldsymbol{\alpha}}  \sum_{i=1}^{\ell(\boldsymbol{\alpha})} 
		\\
		&\sum_{\boldsymbol{\beta} \in \mathcal{C}( \alpha_i +|\kvec_1-\lvec_1|) }  \,\muu_{\alpha_i,\boldsymbol{\beta}} v_{\alpha_1}\star \big(\cdots \star\big( v_{\alpha_{i-1}}\star \big(\,
		\underbrace{ (v_{\beta_1}\star(v_{\beta_2}\star( \cdots \star(v_{\beta_{\ell(\boldsymbol{\beta})}} \star t_1(\lvec_1)) \cdots)))}_{\text{$i$-th position}} \\
		&\hspace{7.3cm}\star \big(v_{\alpha_{i+1}}
		\cdots\star\big(v_{\alpha_{\ell(\boldsymbol{\alpha})}} \star t_2(\lvec_2)\big)\cdots\big)\big)\big)\cdots \big).
	\end{align*}
		
	We rewrite the last expression in the new indices $\widehat{\balpha}^{i} =( \bmu, \beeta) $ as follows
		\begin{align*}
			&\Lie\Bigg(
			v_{\mu_1}\star \big(\cdots \star\big( v_{\mu_{{\ell({\bmu})}}}\star \big((v_{\beta_1}\star(v_{\beta_2}\star( \cdots \star(v_{\beta_{\ell( \bbeta)}} \star t_1(\lvec_1)) \cdots )))\\
			&\hspace{7.5cm} \star \big(v_{\eta_{1}}\cdots \star\big(v_{\eta_{\ell(\beeta)}} \star t_2(\lvec_2)\big)\cdots \big)\big)\big)\cdots \big)\Bigg)\\
			&= \big[ v_{\mu_1} ,\big[\cdots,\big[ v_{\mu_{{\ell({\bmu})}}} ,
			\big[ [v_{\beta_1},[ v_{\beta_2},[ \cdots ,[ v_{\beta_{\ell({\bbeta})}} , T_1(\lvec_1)]\cdots] ]]] , \\
			&\hspace{7.8cm}   \big[v_{\eta_{1}},\big[\cdots ,\big[ v_{\eta_{\ell({\beeta})}}, T_2(\lvec_2) \big] \cdots \big]\big]\big]\big]\cdots \big]\big]\\
			&= \big[ v_{\mu_1} ,\big[\cdots ,\big[ v_{\mu_{{\ell({\bmu})}}} ,\big[\,
			v_{ \bbeta_{(1)}} T_1(\lvec_1) v_{ S(\bbeta_{(2)})}   ,   
			v_{ \beeta_{(1)}}  T_2(\lvec_2)v_{ S(\beeta_{(2)})}
			\big]\big]\cdots \big]\big]\\
			&=  v_{\bmu_{(1)}} v_{ \bbeta_{(1)}} T_1(\lvec_1) v_{ S(\bbeta_{(2)})} v_{ \beeta_{(1)}}  T_2(\lvec_2)v_{ S(\beeta_{(2)})} v_{ S(\bmu_{(2)})} \\
			&\hspace{5cm} -v_{\bmu_{(1)}} v_{ \beeta_{(1)}} T_2(\lvec_2) v_{ S(\beeta_{(2)})} v_{ \bbeta_{(1)}}  T_1(\lvec_1)v_{ S(\bbeta_{(2)})} v_{ S(\bmu_{(2)})} .
		\end{align*}
		
		Thus, we need to understand which triple $( \bmu, \bbeta, \beeta)$ contribute to the coefficient of 
		\[ v_\bsigma T_1(\lvec_1)v_\btau T_2(\lvec_2) v_{S(\btheta)}.\]
		We define 
		\begin{align*}
		C\Big(\begin{smallmatrix} \bsigma & \btau &\btheta \\ i_\bsigma & i_\btau & i_\btheta \end{smallmatrix}\Big)&=
			\Big \{ \bmu \in \N^{i_\bsigma + i_\btheta} , \bbeta \in \N^{\ell(\bsigma)-i_\bsigma + i_\btau}, 
			\beeta \in \N^{\ell(\btau)-i_\btau + \ell(\btheta)-i_\btheta} \,: \\
			&\hspace{-0.5cm}  \big(   v_{\bmu_{(1)}} v_{ \bbeta_{(1)}} T_1(\lvec_1) v_{ S(\bbeta_{(2)})} v_{ \beeta_{(1)}}  T_2(\lvec_2)v_{ S(\beeta_{(2)})} v_{ S(\bmu_{(2)})} \,
			\big| \,v_\bsigma T_1(\lvec_1)v_\btau T_2(\lvec_2) v_{S(\btheta)} 
			\big ) \neq 0 \,\Big\} \\
			& =
			\Big\{ \bmu \in \N^{i_\bsigma + i_\btheta} , \bbeta \in \N^{\ell(\bsigma)-i_\bsigma + i_\btau}, 
			\beeta \in \N^{\ell(\btau)-i_\btau + \ell(\btheta)-i_\btheta} \,:   \\
			& \qquad \quad  \big( \co( \bmu)\,   \big| \,( \sigma_1, \ldots, \sigma_{ i_\sigma  } )\otimes 
			(   \vartheta_1,\ldots, \vartheta_{i_\vartheta} )  \big ) \neq 0, \\
			& \qquad \quad  \big( \co(\bbeta)\,   \big |\,  (\sigma_{ i_\sigma  +1 },\ldots,\sigma_{\ell(\sigma)} )  \otimes   
			(  \tau_{ i_\tau},\ldots, \tau_1 ) \big) \neq 0 ,  \\
			& \qquad \quad \big(   \co(\beeta) \, \big |\,  (\tau_{ i_\tau+1}, \ldots ,\tau_{\ell(\tau)} )\otimes 
			( \vartheta_{ i_\vartheta +1 },\ldots,\vartheta_{\ell(\vartheta)}  )
			\big) \neq 0  \,\Big\} \\
			& =
			\Big\{ \bmu \in \N^{i_\bsigma + i_\btheta} , \bbeta \in \N^{\ell(\bsigma)-i_\bsigma + i_\btau}, 
			\beeta \in \N^{\ell(\btau)-i_\btau + \ell(\btheta)-i_\btheta} \,: \\
			& \qquad \quad  \big( \bmu\,   \big|\, ( \sigma_1, \ldots, \sigma_{ i_\sigma  } )\shuffle 
			(   \vartheta_1,\ldots, \vartheta_{i_\vartheta} ) \big ) \neq 0 , \\
			& \qquad \quad  \big( \bbeta\,   \big |\,  (\sigma_{ i_\sigma  +1 },\ldots,\sigma_{\ell(\sigma)} )  \shuffle   
			(  \tau_{ i_\tau},\ldots, \tau_1 ) \big) \neq 0 ,  \\
			& \qquad \quad \big(    \beeta \, \big  |\,  (\tau_{ i_\tau+1}, \ldots,\tau_{\ell(\tau)} )\shuffle 
			( \vartheta_{ i_\vartheta +1 },\ldots, \vartheta_{\ell(\vartheta)}   )
			\big) \neq 0 \,\Big\} ,
		\end{align*}
		and 
		\begin{align*}
		C(\bsigma , \btau ,\btheta) &= 
			\bigcup_{\substack{0 \le i_\bsigma \le \ell(\bsigma) \\ 
					0 \le i_\btau \le \ell(\btau)\\ 0 \le i_\btheta \le \ell(\btheta)
			}}
		C\Big(\begin{smallmatrix} \bsigma & \btau &\btheta \\ i_\bsigma & i_\btau & i_\btheta \end{smallmatrix}\Big).
		\end{align*} 
		
		Then, we get
		\begin{align*}
			m_{\bsigma, \btau, \btheta}(\lvec_1,\lvec_2)\big( ( t_1  \cdot t_2   ) \tr_u v_a\big)   &:=   
			\Big( \Lie(( t_1(\kvec_1) \cdot t_2(\kvec_2)  ) \tr_u v_a) \, \big |\,  v_\bsigma T_1(\lvec_1)v_\btau T_2(\lvec_2) v_{S(\btheta)} 
			\Big) \\
			&=m_{\kvec_1,\lvec_1}m_{\kvec_2,\lvec_2}\,\,
			\sum\limits_{ (\bmu,\bbeta,\beeta)  \in \mathcal{C}(\bsigma, \btau, \btheta)} 
			\muu_{a,(\bmu, \abs\bbeta -d_1 , \beeta  )} 
			\muu_{\abs\bbeta -d_1,  \bbeta },
		\end{align*} 
		and 
		\begin{align*}
			m_{\bsigma, \btau, \btheta}(\lvec_2,\lvec_1)\big( ( t_1  \cdot t_2   ) \tr_u v_a\big)    &=m_{\kvec_1,\lvec_1}m_{\kvec_2,\lvec_2}\,\,
			\sum\limits_{ (\bmu,\beeta,\bbeta)  \in \mathcal{C}(\bsigma, \btau, \btheta)} 
			\muu_{a,(\bmu, \abs\bbeta -d_1 , \beeta  )} 
			\muu_{\abs\bbeta -d_1,  \bbeta }.
		\end{align*}

		The other previously defined multiplicities occur only if $\btau = \emptyset$. Thus,
		\begin{align*}
			m_{\bsigma,  \btheta}(\lvec_1,\lvec_2)\big( ( t_1  \cdot t_2   ) \tr_u v_a\big)  
			&=   
			\Big( \Lie(( t_1(\kvec_1) \cdot t_2(\kvec_2)  ) \tr_u v_a) \, \big |\,  v_\bsigma T_1(\lvec_1) T_2(\lvec_2) v_{S(\btheta)} 
			\Big) \\
			&=   m_{\bsigma, \emptyset, \btheta}(\lvec_1,\lvec_2)\big( ( t_1  \cdot t_2   ) \tr_u v_a\big) \\
			&=m_{\kvec_1,\lvec_1}m_{\kvec_2,\lvec_2}\,\,
			\sum\limits_{ (\bmu,\bbeta,\beeta)  \in \mathcal{C}(\bsigma, \emptyset, \btheta)} 
			\muu_{a,(\bmu, \abs\bbeta -d_1 , \beeta  )} 
			\muu_{\abs\bbeta -d_1,  \bbeta },
		\end{align*} 
		and
		\begin{align*}
			m_{\bsigma,  \btheta}(\lvec_2,\lvec_1)\big( ( t_1  \cdot t_2   ) \tr_u v_a\big)  
			&=m_{\kvec_1,\lvec_1}m_{\kvec_2,\lvec_2}\,\,
			\sum\limits_{ (\bmu,\beeta,\bbeta)  \in \mathcal{C}(\bsigma, \emptyset, \btheta)} 
			\muu_{a,(\bmu, \abs\bbeta -d_1 , \beeta  )} 
			\muu_{\abs\bbeta -d_1,  \bbeta }.
		\end{align*}
		
		Finally, for $\Lie\big( ( t_2  \cdot t_1 ) \tr_u v_a\big)$ we get by the same considerations the multiplicities
		\begin{align*}
			m_{\bsigma, \btau, \btheta}(\lvec_1,\lvec_2)\big( ( t_2  \cdot t_1   ) \tr_u v_a\big)    &=m_{\kvec_1,\lvec_1}m_{\kvec_2,\lvec_2}\,\,
			\sum\limits_{ (\bmu,\beeta,\bbeta)  \in \mathcal{C}(\bsigma, \btau, \btheta)} 
			\muu_{a,(\bmu, \abs\bbeta -d_2 , \beeta  )} 
			\muu_{\abs\bbeta -d_2,  \bbeta },\\
			m_{\bsigma, \btau, \btheta}(\lvec_2,\lvec_1)\big( ( t_2  \cdot t_1   ) \tr_u v_a\big)    &=m_{\kvec_1,\lvec_1}m_{\kvec_2,\lvec_2}\,\,
			\sum\limits_{ (\bmu,\bbeta,\beeta)  \in \mathcal{C}(\bsigma, \btau, \btheta)} 
			\muu_{a,(\bmu, \abs\bbeta -d_2 , \beeta  )} 
			\muu_{\abs\bbeta -d_2,  \bbeta },\\
			\intertext{and}
			m_{\bsigma,  \btheta}(\lvec_1,\lvec_2)\big( ( t_2  \cdot t_1   ) \tr_u v_a\big)  
			&=m_{\kvec_1,\lvec_1}m_{\kvec_2,\lvec_2}\,\,
			\sum\limits_{ (\bmu,\beeta,\bbeta)  \in \mathcal{C}(\bsigma, \emptyset, \btheta)} 
			\muu_{a,(\bmu, \abs\bbeta -d_2 , \beeta  )} 
			\muu_{\abs\bbeta -d_2,  \bbeta },\\
			m_{\bsigma,  \btheta}(\lvec_2,\lvec_1)\big( ( t_2  \cdot t_1   ) \tr_u v_a\big)  
			&=m_{\kvec_1,\lvec_1}m_{\kvec_2,\lvec_2}\,\,
			\sum\limits_{ (\bmu,\bbeta,\beeta)  \in \mathcal{C}(\bsigma, \emptyset, \btheta)} 
			\muu_{a,(\bmu, \abs\bbeta -d_2 , \beeta  )} 
			\muu_{\abs\bbeta -d_2,  \bbeta }.
		\end{align*}
		
		So we derive that the desired vanishing is equivalent to the identities
		\begin{align*}
			&\sum\limits_{ \balpha \in \mathcal{C}(\bsigma, \btheta)} 
			\,\muu_{a,\balpha }
			=
			\sum\limits_{ (\bmu,\beeta,\bbeta)  \in \mathcal{C}(\bsigma, \emptyset, \btheta)} 
			\muu_{a,(\bmu, \abs\bbeta -d_1 , \beeta  )} 
			\muu_{\abs\bbeta -d_1,  \bbeta } \\
			&\hspace{8cm}
			- \sum\limits_{ (\bmu,\bbeta,\beeta)  \in \mathcal{C}(\bsigma, \emptyset, \btheta)} 
			\muu_{a,(\bmu, \abs\bbeta -d_2 , \beeta  )} 
			\muu_{\abs\bbeta -d_2,  \bbeta }, \\
			&0=
			\sum\limits_{ (\bmu,\beeta,\bbeta)  \in \mathcal{C}(\bsigma, \btau, \btheta)} 
			\muu_{a,(\bmu, \abs\bbeta -d_1 , \beeta  )} 
			\muu_{\abs\bbeta -d_1,  \bbeta }
			- \sum\limits_{ (\bmu,\bbeta,\beeta)  \in \mathcal{C}(\bsigma, \btau, \btheta)} 
			\muu_{a,(\bmu, \abs\bbeta -d_2 , \beeta  )} 
			\muu_{\abs\bbeta -d_2,  \bbeta }.
		\end{align*}
    Note that it is enough to consider the case $(\lvec_2,\lvec_1)$, since the case $(\lvec_1,\lvec_2)$ gives the same equations with $d_1,d_2$ exchanged. These equations are exactly the threshold shuffle identities from Definition \ref{def:threshold_shuffle}, which by Conjecture \ref{conj:uri_mult_threshold_shuffle} should hold for the uri multiplicities.
	\end{proof}

	\begin{Remark}
	The map $\tr_u$ in Definition \ref{def:uri_triangle} can be extended to a derivation on $\Q\langle V\rangle$ by allowing instead of Lie brackets also arbitrary $\Q$-linear combinations of words. As for the Ihara bracket explained in Example \ref{exm:post-lie1}, this naive extension does not agree with the map obtained from Definition \ref{def:ext_post-lie}.
\end{Remark}

\begin{Remark}	For words $w$ of weight less than $9$, there are algorithms that calculate $\Delta_u(w)$ by brute force \cite{Confurius:phd_thesis}. But in contrast to the previous examples, we did not even conjecturally found a closed formula for either the coproduct $\Delta_u$ or the cotriangle map $\cotr_u$ yet.
\end{Remark}

We have the following explicit version of Theorem \ref{thm:dual_pair_free_lie}.

\begin{Theorem}\label{thm:dual_pair_uri}
Let $V=\{v_0,v_1,\ldots\}$ be equipped with the grading $\operatorname{wt}(v_0)=1$ and $\operatorname{wt}(v_i)=i$ for $i>0$. If Conjecture \ref{conj:uri_mult_threshold_shuffle} holds, then we get the graded post-Lie algebra $\big(\Lie(V),\emptybrackets,\tr_u \big)$  associated with the uri bracket and 
  we get a dual pair of graded Hopf algebras 
    \[
    \xymatrix{ 
     (\Q\langle V\rangle,\shuffle,\Delta_u ) 
   \ar@{<~>}[rrr]^{\text{graded dual}}&&&    
   (\Q\langle V \rangle,\glp_u,\co). 
    }
   \]
\graphref{thm:uri_is_post_lie} 
\graphref{thm:dual_pair_free_lie}  
\end{Theorem}
\noproof{Theorem}
 
For the rest of this chapter we assume Conjecture \ref{conj:uri_mult_threshold_shuffle}, i.e., that the uri multiplicities satisfy the threshold shuffle identities from Definition \ref{def:threshold_shuffle}. In this case, the Hopf algebras in Theorem \ref{thm:dual_pair_uri} are well-defined.

\begin{Proposition} \label{prop:glpI_into_glpu}
	(i)	We have an embedding of Hopf algebras 
		\[
		\big( \Q\langle v_0,v_1 \rangle , \glp_I, \co \big)  \hookrightarrow  \big( \QV, \glp_u, \co \big).
		\] 
		(ii) The canonical projection $\QV\rightarrow \Q\langle v_0,v_1\rangle$ induces a surjective Hopf algebra morphism
		\[
		\big( \QV, \shuffle, \Delta_u \big)  \twoheadrightarrow \big( \Q\langle v_0,v_1 \rangle , \shuffle, \Delta_I \big).
		\] 
\end{Proposition}

\begin{proof} If $ a\in \Lie(v_0,v_1)$, then in Definition  \ref{def:uri_triangle}, we must have $\lvec=\kvec$. Hence, we get for $v_i$, $i=0,1$, that 
\[ 
a \tr_u v_i = [v_i,a] = a\tr_I v_i.
\]
This means, the products $\tr_I$ and $\tr_u$, agree on $\Lie (v_0,v_1)$. By Definition \ref{prop:post-lie_relates_lie-algebras} of the post-Lie bracket, we get 
\[
\{a,b\}_I= \{a,b\}_u.
\]
and thus 
the claimed equalities.
\end{proof}

 \begin{Lemma}
\label{lem:sublie_ortho_uri}  Let $i \ge 0$, then the space 
\[
\mathfrak{o}_u(v_i)  = \big\{ w \in \Lie(V\setminus \{v_0\} ) \, \big| \, w \tr_u v_i =0 \big\}
\]
is a Lie subalgebra  of  $(\Lie(V),\emptybrace_u)$.
\end{Lemma}
\begin{proof}
 The proof is the same as for Lemma \ref{lem:sublie_ortho_ari}   
\end{proof}

\begin{Proposition} \label{prop:igeq1_sub_glpu}
	(i)	The tuple $\big( \Q\langle v_1,v_2,v_3,... \rangle , \glp_u, \co \big)$ is a Hopf subalgebra of \\
	$\big( \QV, \glp_u, \co \big) $.
    \vspace{0,3cm} \\
    (ii) The tuple $(\Q\langle v_1,v_2,v_3,\ldots\rangle,\shuffle,\Delta_u)$ is a Hopf subalgebra of $(\QV,\shuffle,\Delta_u)$.   
\end{Proposition}
	
\begin{proof} 
	The result in (i) follows from Lemma \ref{lem:sublie_ortho_uri} applied with $v_i=v_0$. The second claim follows by duality.      \graphref{thm:dual_pair_uri}
\end{proof}

By $\gr_D$ we denote associated graded of any algebraic structure on (a subspace of) $\QV$ with respect to the depth \eqref{eq:dep}. From Lemma \ref{lem:ari_depthgraded_of_uri},  we deduce the following.
\begin{Corollary}\label{cor:ari=depth_graded_uri} There is a natural isomorphism of post-Lie structures
\[
\gr_D  \big( \Lie(V), \emptybrackets, \tr_u \big)
  \cong \big( \Lie(V), \emptybrackets, \tr_a \big).
\]
In particular we have an isomorphism of Hopf algebras
\[
\gr_D \big( \QV, \glp_u, \co \big)
  \cong \big( \QV, \glp_a, \co \big). 
 \] 
\end{Corollary}
The ari bracket $\emptybrace_a$ is therefore the depth-graded version of the uri bracket $\emptybrace_u$.

By duality, the Hopf algebra $\big( \QV, \shuffle, \Delta_a \big)$ is the associated depth-graded of the Hopf algebra $\big( \QV, \shuffle, \Delta_u \big)$. There are relations in the Lie algebra $(\Lie(V),\emptybrace_a)$, see for example Remark \ref{rem:period_relations}, thus we cannot expect a natural morphism between these two Hopf algebras.

\begin{Proposition}\label{prop:v1_ortho_uri}
For  all  $b  \in \Lie(v_0,v_1)$  and   for all 
$w \in \mathfrak{o}_u(v_1)$ we have 
\[
\{w,b\}_u =0.
\]
\end{Proposition}
 \begin{proof}
The proof is the same as for  Proposition \ref{prop:v1_ortho_ari}.  
 \end{proof}

\section{Connection to Ecalle's bimoulds } \label{sec:ecalle}

In \cite{Ecalle02}, Ecallé introduced the Lie algebra BARI of bimoulds, and studied various Lie subalgebras defined with respect to certain symmetries. We rephrase the post-Lie structures studied in Section \ref{sec:ihara}, \ref{sec:ari}, and \ref{sec:uri} in Ecalle's language of bimoulds.

\subsection{Lazard elimination and
 bimoulds }\label{subsec:bimoulds}
Define 
\[C_{k,m}=\operatorname{ad}(v_0)^m(v_k), \qquad k\geq1,\ m\geq0,\]
and let $C^{\operatorname{bi}}$ be the alphabet consisting of all the letters $C_{k,m}$.
By Lazard elimination (cf \cite[Theorem 0.6]{Reut}), the Lie algebra $\Lie(C^{\operatorname{bi}})$ is also free and we have
\[
\Lie(V)=\Lie(C^{\operatorname{bi}})\oplus \Q v_0.
\]
The triangle maps $\tr_I$, $\tr_a$ and $\tr_u$ from 
Definition \ref{def:magma_tr_ihara}, \ref{def:ari_triangle}, and \ref{def:uri_triangle} preserve the space $\Lie(C^{\operatorname{bi}})\subset \Lie(V)$.

Similar to Notation \ref{nota:t_func} and \ref{nota:t_func_on_QV}, we view a Lie element $t\in \Lie(C^{\operatorname{bi}})$ as a function
\begin{align*} 
t:\mathbb{N}^d\times\ZZnonneg^d&\to \Lie(C^{\operatorname{bi}}), \\(\lvec,\mathbf{n})&\mapsto t(\lvec,\mathbf{n}).\end{align*}
So, if $t$ is a multiple Lie bracket of the letters $C_{k_1,m_1}\cdots C_{k_d,m_d}$, then we identify $t=t(\kvec,\mathbf{m})$. Moreover, we extend the multinomial coefficient to $n,n_1\ldots,n_d\in \mathbb{Z}$ by
\[\binom{n}{n_1,\ldots,n_d}=\begin{cases}\binom{n}{n_1,\ldots,n_d}, & \quad n_1+\cdots+n_d=n,\ n_1,\ldots,n_d\in \ZZnonneg, \\0 & \quad \text{else}.\end{cases}\]
 
\begin{Lemma}\label{lem:lazard_tr}
Let $t(\kvec,\mathbf{m})\in \Lie(C^{\operatorname{bi}})$ be a Lie element.

(i) We have for $\tr_I$ with respect to $V_0=\{v_0\}$ 
\begin{align*}
	t(\kvec,\mathbf{m})\tr_I C_{s,n}
	&= \sum_{\mathbf{n} \in \ZZnonneg^d }  \binom{n}{ n-\abs{\mathbf{n}}, \mathbf{n}}[C_{s ,n-\abs{\mathbf{n}}},t(\kvec,\mathbf{m}+\mathbf{n})], 
 \end{align*}
(ii) We have for $\tr_a$ 
\begin{align*}   
t(\kvec,\mathbf{m})\tr_a C_{s,n}
&=\sum_{\lvec,\mathbf{n}\in \ZZnonneg^d} \mua_{\kvec,\lvec} \binom{n}{n-\abs{\mathbf{n}},\mathbf{n}}[C_{s+\abs{\kvec-\lvec},n-\abs{\mathbf{n}}},t(\lvec,\mathbf{m}+\mathbf{n})], 
\end{align*}
(iii) we have for $\tr_u$
\begin{align*}
t(\kvec,\mathbf{m})\tr_u C_{s,n}
&=\sum_{\lvec,\mathbf{n}\in \ZZnonneg^d} \mua_{\kvec,\lvec} \sum_{\substack{\boldsymbol{\alpha}\in \mathcal{C}(s+\abs{\kvec-\lvec}) \\ \boldsymbol{\eta}\in\mathcal{C}^{\ell(\boldsymbol{\alpha})}(n-\abs{\mathbf{n}}) }} \muu_{s,\boldsymbol{\alpha}} \binom{n}{\mathbf{n},\boldsymbol{\eta}}\\
& \hspace{0.4cm} [C_{\alpha_1,\eta_1},[C_{\alpha_2,\eta_2},[\cdots [C_{\alpha_{\ell(\boldsymbol{\alpha})},\eta_{\ell(\boldsymbol{\alpha})}},t(\lvec,\mathbf{m}+\mathbf{n})]\cdots ]]].
\end{align*}
\end{Lemma}

\begin{proof} We can consider the first formula as an identity of elements in $\QV$. By Definition $t(\kvec,\mathbf{m}) \tr_I$ is a derivation for the 
Lie bracket $\emptybrackets$ which vanishes on  $v_0$. Thus, using that also $\ad(v_0)$ is a derivation we get
\begin{align*}
t(\kvec,\mathbf{m})\tr_I C_{s,n} &=    	
t(\kvec,\mathbf{m})\tr_I \ad(v_0)^n (v_s) \\
&= \ad(v_0)^n (t(\kvec,\mathbf{m})\tr_I v_s )\\
&= \ad(v_0)^n (  [  v_s, t(\kvec,\mathbf{m})  ] ) \\
&= \sum_{\mathbf{n} \in \ZZnonneg^d} \binom{n}{ n-\abs{\mathbf{n}}, \mathbf{n}}[C_{s,n-\abs{\mathbf{n}} },t(\kvec,\mathbf{m}+\mathbf{n})].
\end{align*}
The other two cases follow by the same considerations.
\end{proof}
In particular, $(\Lie(C^{\operatorname{bi}}),\emptybrackets,\tr_I)$ and $(\Lie(C^{\operatorname{bi}}),\emptybrackets,\tr_a)$ are post-Lie algebras, and $(\Lie(C^{\operatorname{bi}}),\emptybrackets,\tr_u)$ is a post-Lie algebra under the assumption of Conjecture \ref{conj:uri_mult_threshold_shuffle}.

In the following, a \emph{bimould} means a sequence $A=(A_d)_{d\geq0}$ of polynomials 
$A_d\in \QQ[X_1,Y_1,\ldots,X_d,Y_d]$, such that only finitely many polynomials $A_d$ are nonzero. The space of bimoulds $\operatorname{BIMU}$ equipped with the multiplication $\operatorname{mu}$ given by
\[\operatorname{mu}(A,B)_d\binom{X_1,\ldots,X_d}{Y_1,\ldots,Y_d}=
\sum_{i=0}^d A_i\binom{X_1,\ldots,X_i}{Y_1,\ldots,Y_i}B_{d-i}\binom{X_{i+1},\ldots,X_d}{Y_{i+1},\ldots,Y_d}
\]
is an associative, unital $\Q$-algebra.

We explain how to understand $\Lie(C^{\operatorname{bi}})$ equipped with the different structures from Section \ref{sec:ihara}, \ref{sec:ari}, and \ref{sec:uri} in the context of bimoulds, for details we refer to \cite{BPhD}.

By the \emph{depth} of a word in $\Q\langle C^{\operatorname{bi}}\rangle$ we mean the number of letters. We write $\Q\langle C^{\operatorname{bi}}\rangle^{(d)}$ for the homogeneous subspace of $\Q\langle C^{\operatorname{bi}}\rangle$ spanned by the words of depth $d$, and for an element $f\in \Q\langle C^{\operatorname{bi}}\rangle$ we write $f^{(d)}$ for the homogeneous part of $f$ of depth $d$. We have a $\Q$-linear map
\begin{align}\label{eq:def_rho_C}
\rho_{C^{\operatorname{bi}}}: \Q\langle C^{\operatorname{bi}}\rangle &\to \QQ[X_1,Y_1,X_2,Y_2,\ldots], \\
C_{k_1,m_1}\cdots C_{k_d,m_d}&\mapsto X_1^{k_1-1}Y_1^{m_1}\cdots X_d^{k_d-1}Y_d^{m_d}, \nonumber
\end{align}
which restricts to isomorphisms $\rho_{C^{\operatorname{bi}}}:\Q\langle C^{\operatorname{bi}}\rangle^{(d)}\to \Q[X_1,Y_1,\ldots,X_d,Y_d]$, $d\geq0$.

Observe that the map $\rho_{C^{\operatorname{bi}}}$ maps a word of weight $w$ and depth $d$ to a homogeneous polynomial of degree $w-d$ in $2d$ variables.  

\begin{Definition}
For each $f\in \Q\langle C^{\operatorname{bi}}\rangle$, we associate a bimould $\rho_{C^{\operatorname{bi}}}(f)=(\rho_{C^{\operatorname{bi}}}(f)_d)_{d\geq0}$ by 
\[\rho_{C^{\operatorname{bi}}}(f)_d\binom{X_1,\ldots,X_d}{Y_1,\ldots,Y_d}=\rho_{C^{\operatorname{bi}}}(f^{(d)}), \qquad d\geq0.\]
\end{Definition}

\begin{Proposition}
The map  
\[
\rho_{C^{\operatorname{bi}}}:(\Q\langle C^{\operatorname{bi}}\rangle , \conc ) \to (\operatorname{BIMU}, \operatorname{mu})	
\]
is an algebra isomorphism. 	
\end{Proposition}

\subsection{The ari bracket on bimoulds}

Because of \cite[Lem. 5.61, Prop. 5.63]{BPhD}, we  can define the Lie algebra
\begin{align} \label{eq:def_bari_al}
\big(\operatorname{BARI}_{\operatorname{al}},\emptybrackets_{\operatorname{mu}}\big)  = \rho_{C^{\operatorname{bi}}} \big(\Lie(C^{\operatorname{bi}}),\emptybrackets\big).
\end{align}
Of course, $\emptybrackets_{\operatorname{mu}}$ equals the commutator bracket induced by the multiplication $\operatorname{mu}$.
This particular subspace $\operatorname{BARI}_{\operatorname{al}}$ of 
bimoulds has been called by Ecalle the space of \emph{alternal} 
bimoulds. By \cite[Thm 1.4]{Reut} and duality, an element 
$f\in \Q\langle C^{\operatorname{bi}}\rangle$ is a Lie element if and only if the 
coefficients of $f$ multiply to $0$ under the shuffle product. 
Alternality just means that the coefficients of a bimould $A$ multiply to $0$ under the 
shuffle product, for the precise definition we refer to \cite[Section 2.3]{Sch15}.
 
The definition of $(\operatorname{BARI}_{\operatorname{al}},
\emptybrackets_{\operatorname{mu}})$ 
in \eqref{eq:def_bari_al} directly implies that it is a free Lie algebra generated 
by bimoulds $p_{k,m}$, whose depth $1$ component is homogeneous and all others are $0$. Precisely, we have  
\[p_{k,m}\binom{X_1}{Y_1}= \rho_{C^{\operatorname{bi}}}(C_{k,m})=X_1^{k-1}Y_1^m, \qquad k\geq1,\ m\geq0.\]

Any post-Lie structure $(\Lie(C^{\operatorname{bi}}),\emptybrackets, \tr)$ induces a post-Lie structure $(\operatorname{BARI}_{\operatorname{al}},
\emptybrackets_{\operatorname{mu}}, \tr)$ via 
\[
A \tr B = \rho_{C^{\operatorname{bi}}} \big( \rho_{C^{\operatorname{bi}}}^{-1} (A)   \tr \rho_{C^{\operatorname{bi}}}^{-1}(B) \big) .
\]
for all $A,B \in \operatorname{BARI}_{\operatorname{al}}$.
For the maps $\tr_I$ and $\tr_a$ on $\operatorname{BARI}_{\operatorname{al}}$ we have the following.

 \begin{Lemma}\label{lem:tr_on_bari_al} For $k,s\geq1$ and $m,n\geq0$, we have
\begin{align*}
\text{(i)} \qquad p_{k,m} \tr_I p_{s,n} \binom{X_1,X_2}{Y_1,Y_2}&= [ p_{k,m} , p_{s,n} ]_{\operatorname{mu}} \binom{X_1,X_2}{Y_1+Y_2,Y_2}, \\
\text{(ii)} \qquad 
p_{k,m} \tr_a p_{s,n} \binom{X_1,X_2}{Y_1,Y_2} &= [ p_{k,m} , p_{s,n} ]_{\operatorname{mu}} \binom{X_1,X_2-X_1}{Y_1+Y_2,Y_2}.
\end{align*}
\end{Lemma}

For the map $\tr_u$ it is much more difficult to find a closed formula, since it also has higher depth components and the additional uri multiplicities $\muu_{a,\boldsymbol{\alpha}}$.

\begin{proof} The cases (i) and (ii) are similar. Thus, we only consider (ii) 
\begin{align*}
&p_{k,m} \tr_a p_{s,n} \binom{X_1,X_2}{Y_1,Y_2} =  \rho_{C^{\operatorname{bi}}}  \big( C_{k,m} \tr_a C_{s,n} \big)\binom{X_1,X_2}{Y_1,Y_2}  \\
&=  \sum_{l,\eta \in \ZZnonneg}  \mua_{k,l} \binom{n}{ \eta } \rho_{C^{\operatorname{bi}}}  \Big( [C_{s+ k-l,n- \eta},C_{l,m+\eta} ]  \Big) \binom{X_1,X_2}{Y_1,Y_2}\\
&= \sum_{l,\eta \in \ZZnonneg}  (-1)^{k+l} \binom{k-1}{l-1} \binom{n}{ \eta }
\big( X_1^{s+k-l-1} Y_1^ {n-\eta} X_2^{l-1} Y_2^{m+\eta}   
-  X_1^{l-1} Y_1^{m+\eta} X_2^{s+k-l-1} Y_2^ {n-\eta} \big)\\
&=  X_1^{s-1} (Y_1+Y_2)^n (X_2-X_1)^{k-1} Y_2^m - 
(X_2-X_1)^{s-1} Y_1^n   X_2^{k-1} (Y_1+Y_2)^m \\
&= [ p_{k,m} , p_{s,n} ]_{\operatorname{mu}} \binom{X_1,X_2-X_1}{Y_1+Y_2,Y_2}. \qedhere
\end{align*} 
\end{proof}
 
Since we determined $\tr_a$ on the generators, we could use the general formalism to calculate the Grossman-Larson product in $\operatorname{BIMU}$. Indeed, by the recursive Definition \ref{def:ext_post-lie} any of these products is expressible in the terms we get from  Lemma \ref{lem:tr_on_bari_al}. 
In Remark \ref{rem:glp_post-lie_bracket_algorithm} we express the post-Lie bracket in terms of Grossman-Larson products, thus we get some expression for the induced post-Lie bracket $\emptybrace_a$ on   $\operatorname{BARI}_{\operatorname{al}}$. 
Actually, a better and direct way to calculate the ari bracket $\emptybrace_a$  had been known for long before. 
As in \cite[2.2]{Ecalle11}, we define for two bimoulds $A,B$ the operation $\operatorname{arit}$ by
\begin{align*}&
\operatorname{arit}(B)(A)_d\binom{X_1,\ldots,X_d}{Y_1,\ldots,Y_d}\\
&=\sum_{j=1}^d\sum_{i=1}^{d-j} A_{d-j}\binom{X_1,\ldots,X_{i-1},X_{i+j},X_{i+j+1},\ldots,X_d}{Y_1,\ldots,Y_{i-1},Y_i+\cdots+Y_{i+j},Y_{i+j+1},\ldots,Y_d} \\
&\hspace{8cm} \cdot B_j\binom{X_i-X_{i+j},\ldots, X_{i+j-1}-X_{i+j}}{Y_i,\ldots,Y_{i+j-1}} \\
&-\sum_{j=1}^d\sum_{i=1}^{d-j} A_{d-j}\binom{X_1,\ldots,X_{i-1},X_i,X_{i+j+1},\ldots,X_d}{Y_1,\ldots,Y_{i-1},Y_i+\cdots+Y_{i+j},Y_{i+j+1},\ldots,Y_d} \\
&\hspace{8cm} \cdot B_j\binom{X_{i+1}-X_i,\ldots, X_{i+j}-X_i}{Y_{i+1},\ldots,Y_{i+j}}.
\end{align*}

An immediate consequence of \cite[Thm 5.52.]{BPhD} is the following.

\begin{Theorem} \label{thm:iso_ari_bari} We have an isomorphism of post-Lie algebras
\[\big(\Lie(C^{\operatorname{bi}}),\emptybrackets,\tr_a\big)\to \big(\operatorname{BARI}_{\operatorname{al}},\emptybrackets_{\operatorname{mu}},\operatorname{arit}\big),\quad f \mapsto \rho_{C^{\operatorname{bi}}}(f).\]
\end{Theorem}

\begin{proof}
Define the derivation
\begin{align*}   
d^A_{C_{k_1,m_1}\cdots C_{k_d,m_d}}(C_{s,n})
&=\sum_{\lvec,\mathbf{n}\in \ZZnonneg^d} \mua_{\kvec,\lvec} \binom{n}{n-\abs{\mathbf{n}},\mathbf{n}}[C_{s+\abs{\kvec-\lvec},n-\abs{\mathbf{n}}},C_{l_1,m_1+n_1}\cdots C_{l_d,m_d+n_d}].
\end{align*}
Then the $\Q$-linearity directly implies that $\tr_a$ and $d^A$ agree on the Lie algebra $\Lie(C^{\operatorname{bi}})$. Straight-forward calculations, that extend those of the proof of Lemma \ref{lem:tr_on_bari_al},  show now (see \cite[Proposition 5.70]{BPhD}) that
\[\rho_{C^{\operatorname{bi}}}(d^A_f(g))=\operatorname{arit}(\rho_{C^{\operatorname{bi}}}(f))(\rho_{C^{\operatorname{bi}}}(g)),\qquad f,g\in \Q\langle C^{\operatorname{bi}}\rangle.\]
So, restricting to $\Lie(C^{\operatorname{bi}})$, we deduce
\[\rho_{C^{\operatorname{bi}}}(f\tr_a g)=\operatorname{arit}(\rho_{C^{\operatorname{bi}}}(f))(\rho_{C^{\operatorname{bi}}}(g)),\qquad f,g\in \Lie(C^{\operatorname{bi}}).\]
This implies Theorem \ref{thm:iso_ari_bari}.
\end{proof}

The associated post-Lie bracket $\emptybrace_a$ corresponds to the Ecalle's ari bracket, see \cite[2.2]{Ecalle11}, under the map $\rho_{C^{\operatorname{bi}}}$.

\begin{Remark}\label{ref:linear_shuffle_embeds}
The subspace of alternal and swap invariant bimoulds $\operatorname{BARI}_{\underline{\operatorname{al}},\operatorname{swap}}$, where in addition the depth 1 component is required to be even, is a Lie subalgebra of $(\operatorname{BARI}_{\operatorname{al}},\operatorname{ari})$. The proof is with minor modifications an application of \cite[Lem 2.4.1, Prop 2.5.2, and Lem 2.5.5]{Sch15} and is revisited  in  \cite{Bu_extended}. This Lie algebra should contain a Lie subalgebra corresponding to the associated depth-graded of the algebra $\Z_q$ of multiple $q$-zeta values, its conjectural structure has been presented in \cite{K-montreal}.

Recently, the first author studied in \cite{Bu_extended} a Lie subalgebra $\mathfrak{lq}$ of $(\Lie(V),\emptybrace_a)$, which is related to the associated depth-graded of the algebra $\Z_q$ of multiple $q$-zeta values in a similar way. In op. cit. it is shown that 
we have a Lie algebra isomorphism
\[(\mathfrak{lq},\emptybrace_a) \to (\operatorname{BARI}_{\underline{\operatorname{al}},\operatorname{swap}}, \operatorname{ari}), \] and an embedding 
\[
(\mathfrak{ls}, \emptybrace_I ) \hookrightarrow
(\mathfrak{lq},\emptybrace_a),   \]
where the linearised double shuffle Lie algebra $\mathfrak{ls}$ corresponds to depth-graded multiple zeta values, see e.g. \cite{Ra}, \cite{brown-depth}, and \cite{Sch15}.

\end{Remark}

\subsection{The uri bracket on bimoulds} 

The map $\tr_u$ from Section \ref{sec:uri} seems to be hard to translate to the space $\operatorname{BARI}_{\operatorname{al}}$. Though, we want to explain a connection to other Lie structures expected on $\operatorname{BIMU}$.

Consider the algebra morphism $\log_*:\Q\langle C^{\operatorname{bi}}\rangle \to \Q\langle C^{\operatorname{bi}}\rangle$ given by
\begin{align*}
\log_*(C_{k,m})&=
\frac{1}{m!} \sum\limits_{\substack{l_1+\cdots+l_d=k \\ l_j\geq1}} \frac{(-1)^{d+1}}{d} \ad(v_0)^m(v_{l_1}\cdots v_{l_d}) \\
&=\sum\limits_{\substack{l_1+\cdots+l_d=k \\ n_1+\cdots+n_d=m \\ l_j\geq1,\ n_j\geq0}} \frac{(-1)^{d+1}}{d} \frac{1}{n_1!\cdots n_d!} C_{l_1,n_1}\cdots C_{l_d,n_d}.
\end{align*}

This map is dual to the Hoffman logarithm, see \cite{HI17}, which describes the isomorphism of the (bi-)stuffle Hopf algebra with the shuffle Hopf algebra  \cite{BPhD}.

\begin{Proposition}
We have a homomorphism
 \begin{align*}
\rho_{D^{\operatorname{bi}}}=\rho_{C^{\operatorname{bi}}}\circ\log_*: (\Q\langle C^{\operatorname{bi}}\rangle ,\conc)\to ( \operatorname{BIMU} , \operatorname{mu}).
\end{align*}   
\end{Proposition}

We define the Lie algebra
\begin{align}
(\operatorname{BARI}_{\operatorname{il}},\emptybrackets_{\operatorname{mu}})=\rho_{D^{\operatorname{bi}}}\big(\Lie(C^{\operatorname{bi}}),\emptybrackets\big).
\end{align}
The subspace $\operatorname{BARI}_{\operatorname{il}}$ equals the space of \emph{alternil} bimoulds defined by Ecalle, this follows by a combination of results in \cite[Section 5.4]{BPhD}. The elements in $\log_*(\Lie(C^{\operatorname{bi}}))$ are the primitive elements under a certain coproduct $\Delta_\ast$. Thus, the coefficients of $f\in \log_*(\Lie(C^{\operatorname{bi}}))$ multiply to $0$ under the dual product $\ast$, usually called the (bi-) stuffle product. Alternility means that the coefficients of a bimould $A$ multiply to $0$ under the stuffle product. The precise definition is given in \cite[Section 2.3]{Sch15}.

As in \cite[Definition 5.12]{BPhD}, we define for bimoulds $A,B$ the operation $\operatorname{urit}$ by
\begin{align*}
\operatorname{urit}(B)(A)_d\binom{X_1,\ldots,X_d}{Y_1,\ldots,Y_d} &= \operatorname{arit}(B)(A)_d\binom{X_1,\ldots,X_d}{Y_1,\ldots,Y_d} \\
&\hspace{-4.3cm}+\sum_{l=2}^d\sum_{i=1}^{d-l} A_{d-l}\binom{X_1,\ldots,X_{i-1},X_i,X_{i+l+1},\ldots,X_d}{Y_1,\ldots,Y_{i-1},Y_i+\cdots+Y_{i+l},Y_{i+l+1},\ldots,Y_d} \\
&\hspace{-3.8cm}\cdot\Bigg(\sum_{r=1}^{l-1}\sum_{s=0}^r \prod_{\substack{t=0 \\ t\neq s}}^r \frac{1}{X_{i+t}-X_{i+s}}B_{l-r}\binom{X_{i+r+1}-X_{i+s},\ldots,X_{i+l}-X_{i+s}}{Y_{i+r+1},\ldots,Y_{i+l}} \\
& \hspace{-3.2cm} -\sum_{r=1}^{l-1}\sum_{\substack{s=0 \\ \text{or } s=l}}^{r-1} \prod_{\substack{t=0 \\ \text{or } t=l \\ t\neq s}}^{r-1} \frac{1}{X_{i+t}-X_{i+s}}B_{l-r}\binom{X_{i+r}-X_{i+s},\ldots,X_{i+l+1}-X_{i+s}}{Y_{i+r},\ldots,Y_{i+l-1}}\Bigg).
\end{align*}
By \cite[Prop 5.14]{BPhD}, the operation $\operatorname{urit}$ preserves the space $\operatorname{BIMU}$.

\begin{Conjecture} \label{conj:uri_bari_iso} There is an isomorphism of post-Lie algebras
\[ (\Lie(C^{\operatorname{bi}}),\emptybrackets,\tr_u)\to (\operatorname{BARI}_{\operatorname{il}},\emptybrackets_{\operatorname{mu}},\operatorname{urit}),\quad f\mapsto \rho_{D^{\operatorname{bi}}}(f). \]
\end{Conjecture}
In particular, the associated post-Lie bracket $\emptybrace_u$ should endow $\operatorname{BARI}_{\operatorname{il}}$ with a Lie algebra structure. 

\begin{Example} We have
\begin{align*}
\rho_{D^{\operatorname{bi}}} (C_{3,0}) &= \big(X_1^2, -\frac{1}{2}(X_1+X_2),\frac{1}{3},0,\ldots\big), \\
\rho_{D^{\operatorname{bi}}} (C_{2,1}) &= \big(X_1Y_1,-\frac{1}{2}(Y_1+Y_2),0,\ldots\big).
\end{align*}
Thus, we compute  
\begin{align*}
\operatorname{urit}(\rho_{D^{\operatorname{bi}}}(C_{3,0}))(\rho_{D^{\operatorname{bi}}}(C_{2,1}))&=\big(0,(X_1^3-X_2^3+3X_1X_2^2-3X_1^2X_2)(Y_1+Y_2), \\
&(X_1^2-X_3^2-\frac{1}{2}X_1X_2-3X_1X_3+\frac{7}{2}X_2X_3)(Y_1+Y_2+Y_3),\\ 
&(\frac{13}{12}X_1-\frac{3}{4}X_2-\frac{3}{4}X_3+\frac{5}{12}X_4)(Y_1+Y_2+Y_3+Y_4),0,\ldots\big)
\end{align*}
On the other hand, we compute with Lemma \ref{lem:lazard_tr} 
\begin{align*}
C_{3,0} \tr_u C_{2,1}=\ & 3[C_{2,1},C_{3,0}]+3[C_{2,0},C_{3,1}]+[C_{4,1},C_{1,0}]+[C_{4,0},C_{1,1}]\\ &+\frac{3}{2}[C_{2,1},[C_{2,0},C_{1,0}]]+\frac{3}{2}[C_{2,0},[C_{2,1},C_{1,0}]]+\frac{3}{2}[C_{2,0},[C_{2,0},C_{1,1}]]\\
&-\frac{1}{4}[C_{1,1},[C_{1,0},[C_{2,0},C_{1,0}]]]-\frac{1}{4}[C_{1,0},[C_{1,1},[C_{2,0},C_{1,0}]]]\\
&-\frac{1}{4}[C_{1,0},[C_{1,0},[C_{2,1},C_{1,0}]]]-\frac{1}{4}[C_{1,0},[C_{1,0},[C_{2,0},C_{1,1}]]].
\end{align*}
One easily verifies that these two calculations coincide under the map $\rho_{D^{\operatorname{bi}}}$.
\end{Example}

\begin{Remark}
There is another approach by the second author and Schneps \cite{K-montreal} to a Lie bracket on $\operatorname{BARI}_{\operatorname{il}}$ that transforms the ari bracket on $\operatorname{BARI}_{\operatorname{al}}$ into the so-called uri bracket by composing it with maps $\operatorname{ganit}_{\operatorname{pic}}$ and $\operatorname{ganit}_{\operatorname{poc}}$. The maps $\operatorname{ganit}_{\operatorname{pic}}$, $\operatorname{ganit}_{\operatorname{poc}}$ map into the space of bimoulds, whose components are rational functions instead of polynomials. So it is not known yet, whether the uri bracket respects the space $\operatorname{BIMU}$, as we defined it above. We expect that this construction coincides with the post-Lie bracket induced by $\operatorname{urit}$. The first author showed in her thesis \cite{BPhD} that this is in fact the case for depth $d\le 6$.    
\end{Remark}

\begin{Remark} 
Recall by Corollary \ref{cor:ari=depth_graded_uri} we have a natural isomorphism of post-Lie algebras
\[\gr_D(\Lie(V),\emptybrackets,\tr_u)\simeq (\Lie(V),\emptybrackets,\tr_a).\]

In the context of bimoulds this corresponds to 
\[
\big(\gr_D\operatorname{BARI}_{\operatorname{il}},\emptybrackets_{\operatorname{mu}},\gr_D \operatorname{urit}\big)\simeq \big(\operatorname{BARI}_{\operatorname{al}},\emptybrackets_{\operatorname{mu}},\operatorname{arit}\big).
\] 
By construction we have  $\gr_D \operatorname{urit} =\operatorname{arit}$. 
Since $\gr_D\log_*=\operatorname{id}$, we have $\gr_D\rho_{D^{\operatorname{bi}}}=\rho_{C^{\operatorname{bi}}}$ and therefore
\[\gr_D \operatorname{BARI}_{\operatorname{il}}
=\gr_D ( \rho_{D^{\operatorname{bi}}} (\Lie(C^{\operatorname{bi}})) )  
= \rho_{C^{\operatorname{bi}}}(   \Lie(C^{\operatorname{bi}}) )
= \operatorname{BARI}_{\operatorname{al}}.
\]
We also observe that Ecalle's definition of alternil implies that forgetting all higher depth terms of an  alternil bimould yields  an alternal bimould. 
Note that above isomorphism holds even without assuming Conjecture \ref{conj:uri_bari_iso}. In other words,   although we do not know whether $\operatorname{urit}$ preserves $\operatorname{BARI}_{\operatorname{il}}$, the statement holds for the associated depth-graded structures.

\end{Remark}

\begin{Remark}
The first author studied in \cite{BPhD} a space $\mathfrak{bm}_0\subset \Lie(V)$, which should be isomorphic to the space of indecomposables of the algebra $\Z_q$ of multiple $q$-zeta values modulo quasi-modular forms. We expect that $\mathfrak{bm}_0$ is a Lie subalgebra of $(\Lie(V),\emptybrace_u)$, and that the vector space isomorphism $\mathfrak{bm}_0\to \operatorname{BARI}_{\underline{\operatorname{il}},\operatorname{swap}}$ in \cite[Cor 5.51]{BPhD} extends to a Lie algebra isomorphism
\[\big(\mathfrak{bm}_0,\emptybrace_u\big)\to \big(\operatorname{BARI}_{\underline{\operatorname{il}},\operatorname{swap}},\operatorname{uri}\big), \quad f\mapsto\rho_{D^{\operatorname{bi}}}(f).\]
Moreover, the embedding of vector spaces $\mathfrak{dm}_0\hookrightarrow \mathfrak{bm}_0$ in \cite[Theorem 4.28]{BPhD} should give rise to an injective Lie algebra morphism
\[\big( \mathfrak{dm}_0 , \emptybrace_I \big) \hookrightarrow  \big(\mathfrak{bm}_0,\emptybrace_u\big),\]
where $\mathfrak{dm}_0$ is the double shuffle Lie algebra \cite{Ra}. Observe that these expectations actually hold in the associated depth-graded case, c.f. Remark
\ref{ref:linear_shuffle_embeds}.
In \cite{K-montreal} a proof for the claim that the image of a similar map of $\mathfrak{dm}_0$  
into $\operatorname{BARI}_{\operatorname{il},\operatorname{swap}}$  is a Lie algebra with respect to the uri bracket $\emptybrace_u$ is announced. 
\end{Remark}

\appendix

\section{Identities 
via threshold functions 
on compositions} \label{app:comb_id}

For $\mathbf{a}=(a_1,\ldots,a_s) \in \mathbb{N}^s$, we set $|\mathbf{a}| =a_1+\dots + a_s$.
Thus we can view $\mathbf{a}$ as a composition of $|\mathbf{a}|$ with $s$ parts. We denote the set of such compositions  by
\begin{align*}
 \mathcal{C}^s(n) &= \{ \mathbf{a} \in \mathbb{N}^s \mid |\mathbf{a}| = n \}.
\end{align*}
We have for all $a\in \mathbb{N}$ a threshold indicator function\footnote{A similar construction was used by Andrews in \cite{andrews:compositions} in the proof of Theorem 3.}  
on the set of all compositions $\mathcal{C}= \bigcup_{s,n\geq0} \mathcal{C}^s(n)$, which is given by
\begin{align} \label{eq:app_def_threshold}
 \ind{a}{\mathbf{a}}&=\min_{1\leq j\leq \ell(\mathbf{a})} \left\{ j \mid a_1 + \dots + a_j \geq a \right\}.
\end{align}      
By convention we set $\ind{a}{\mathbf{a}}=0$, if either $\mathbf{a} = \emptyset$, $a > |\mathbf{a}|$.
In the following we want to study some properties of $\ind{a}{\mathbf{a}}$. 

\begin{Proposition}
We have for $n,j,r,s\in\mathbb{N}$ that
\[
\# \big\{  \mathbf{a} \in \mathcal{C}^s(n+r) \,| \,\ind{n}{\mathbf{a}}= j \,\big\} 
= \binom{n-1}{j-1} \binom{r}{s-j}.
\]
\end{Proposition}
\begin{proof}  We use the bijection
\begin{align*}
\big\{  \mathbf{a} \in \mathcal{C}^s(n+r) \,| \,\ind{n}{\mathbf{a}}= j \,\big\} 
& \to    \mathcal{C}^j(n)  \times  \bigcup_{t=0}^{r-1}  \mathcal{C}^{s-j}(r-t), \\  
\mathbf{ a} &\mapsto \Big( (  a_1, \ldots , a_j- t )  ,(a_{j+1},\ldots,a_s  ) \Big).
\end{align*}
Now using MacMahon's observation  
\[
\# \mathcal{C}^j(n) = \binom{n-1}{j-1},
\]
the claim follows since,
\begin{align*}
\# \big(  \bigcup_{t=0}^{r-1}  \mathcal{C}^{s-j}(r-t) \big) &= 
\sum_{t=0}^{r-1} \# \mathcal{C}^{s-j}(r-t) =\sum_{t=0}^{r-1}  \binom{r-t-1}{s-j-1} = \binom{ r}{s-j}. \qedhere
\end{align*}
\end{proof} 

\begin{Remark} Using the fact that the threshold indicator functions provide a disjoint decomposition of  $\mathcal{C}^s(n)$, we get as an immediate consequence the well known identity 
\[
\sum_{j=1}^s \binom{n-1}{j-1} \binom{ r}{s-j} = \binom{n+r-1}{s-1}, \qquad n,r,s\in\mathbb{N}.
\]
\end{Remark}

\subsection{Identities for binomial coefficients 
}

For $\mathbf{d}=(d_1,\ldots,d_l),\ \mathbf{r}=(r_1,\ldots,r_l) \in \ZZnonneg^l$, we write $\mathbf{r} \le \mathbf{d}$ if $d_i-r_i \ge 0$ for all $i=1,\ldots,l$, and we set
\begin{align*}
\binom{\mathbf{d}}{\mathbf{r}} &= \binom{d_1}{r_1} \cdots \binom{d_l}{r_l}.
\end{align*}

The following identity is applied in Section \ref{sec:uri}.

\begin{Theorem} \label{thm:conj_uri_prop}
For $\mathbf{d} \in \ZZnonneg^l$, $\mathbf{b} \in \mathbb{N}^s$, and $n \in \mathbb{N}$, we have for any $1\leq j\leq s$
    \begin{align*}
     \sum_{\mathbf{0} \leq \mathbf{r} \leq \mathbf{d}} (-1)^{|\mathbf{d}-\mathbf{r}|}
    \binom{\mathbf{d}}{\mathbf{r}}
    \sum_{\substack{\mathbf{a} \in \mathcal{C}^s(n + |\mathbf{r}|) \\ \ind{n}{\mathbf{a}}=j}} 
    \binom{\mathbf{a}-1}{\mathbf{b}-1}   = \begin{cases} \binom{n-1}{|\mathbf{b}| - |\mathbf{d}|-1} ,& j= \ind{|\mathbf{b}| - |\mathbf{d}|}{\mathbf{b}},\\
        0 &\text{else}.
    \end{cases}
    \end{align*}
\end{Theorem}

\begin{proof} By definition, we have for $r\in\ZZnonneg$, $j,n\in\mathbb{N}$
\begin{align} \label{eq:conj_uri_prop_form0}
&\{\mathbf{a}\in\mathcal{C}^s(n+r)\mid \ind{n}{\mathbf{a}}=j\} \nonumber \\
&=\{\mathbf{a}\in\mathcal{C}^s(n+r)\mid a_1+\cdots+a_j\geq n, a_1+\cdots a_{j-1}<n\} \nonumber\\
&=\{\mathbf{a}\in\mathcal{C}^s(n+r)\mid a_1+\cdots+a_{j-1}<n\}\backslash\{\mathbf{a}\in\mathcal{C}^s(n+r)\mid a_1+\cdots+a_j<n\}.
\end{align}
Let $\mathbf{b}=(b_1,\ldots,b_s)\in\mathbb{N}^s$. With the binomial identity
\begin{align} \label{eq:conj_uri_prop_binomial_id} 
\sum_{a_1+\cdots+a_j=n} \prod_{i=1}^s \binom{a_i-1}{b_i-1}=\binom{n-1}{b_1+\cdots+b_j-1},
\end{align}
we get for each $r\in \ZZnonneg$, $j,n\in\mathbb{N}$
\begin{align} \label{eq:conj_uri_prop_form2}
\sum_{\substack{a_1+\cdots+a_s=n+r \\ a_1+\cdots+a_j<n}} \prod_{i=1}^s \binom{a_i-1}{b_i-1}&=\sum_{m=1}^{n-1} \sum_{\substack{a_1+\cdots+a_j=n-m \\ a_{j+1}+\cdots+a_s=r+m}} \prod_{i=1}^s \binom{a_i-1}{b_i-1} \\
&=\sum_{m=1}^{n-1} \binom{n-m-1}{b_1+\cdots+b_j-1}\binom{r+m-1}{b_{j+1}+\cdots+b_s-1}.  \nonumber
\end{align}
Let $\mathbf{d}=(d_1,\ldots,d_l)$. With \eqref{eq:conj_uri_prop_form0} and \eqref{eq:conj_uri_prop_form2}, we get that the left hand side of the claimed formula is equal to
\begin{align} \label{eq:conj_uri_prop_form3}
&\sum_{r_1=0}^{d_1}\cdots \sum_{r_l=0}^{d_l} \prod_{i=1}^l (-1)^{d_i-r_i}\binom{d_i-1}{r_i-1}\Bigg(\sum_{m=1}^{n-1} \binom{n-m-1}{b_1+\cdots+b_{j-1}-1}\binom{r_1+\cdots+r_l+m-1}{b_j+\cdots+b_s-1} \nonumber \\
&\hspace{5.7cm}- \sum_{m=1}^{n-1} \binom{n-m-1}{b_1+\cdots+b_j-1}\binom{r_1+\cdots+r_l+m-1}{b_{j+1}+\cdots+b_s-1}\Bigg) \nonumber \\ 
&=\sum_{m=1}^{n-1} \binom{n-m-1}{b_1+\cdots+b_{j-1}-1} \sum_{r_1=0}^{d_1}\cdots \sum_{r_l=0}^{d_l} \prod_{i=1}^l (-1)^{d_i-r_i}\binom{d_i-1}{r_i-1} \binom{r_1+\cdots+r_l+m-1}{b_j+\cdots+b_s-1} \nonumber \\
&-\sum_{m=1}^{n-1} \binom{n-m-1}{b_1+\cdots+b_j-1} \sum_{r_1=0}^{d_1}\cdots \sum_{r_l=0}^{d_l} \prod_{i=1}^l (-1)^{d_i-r_i}\binom{d_i-1}{r_i-1} \binom{r_1+\cdots+r_l+m-1}{b_{j+1}+\cdots+b_s-1}.
\end{align}
For all $k,m\in\ZZnonneg, k\leq d_1+\cdots+d_l+m$, we have the identity
\[\sum_{r_1=0}^{d_1}\cdots \sum_{r_l=0}^{d_l} \prod_{i=1}^l (-1)^{r_i}\binom{d_i-1}{r_i-1} \binom{m+d_1-r_1+\cdots+d_l-r_l}{k-r_1-\cdots-r_l}=\binom{m}{m-k},\]
which can be checked via generating series \cite[Lemma 3.23]{BPhD}. This identity applied to the expression in \eqref{eq:conj_uri_prop_form3} yields (choose $k=m+d_1+\cdots+d_l-b_j-\cdots-b_s$ resp. $k=m+d_1+\cdots+d_l-b_{j+1}-\cdots-b_s$)
\begin{align} \label{eq:conj_uri_prop_form4}
&\sum_{m=1}^{n-1} \binom{n-m-1}{b_1+\cdots+b_{j-1}-1} \binom{m-1}{b_j+\cdots+b_s-d_1-\cdots-d_l-1}\\
&\hspace{3.5cm}-\sum_{m=1}^{n-1}\binom{n-m-1}{b_1+\cdots+b_j-1} \binom{m-1}{b_{j+1}+\cdots+b_s-d_1-\cdots-d_l-1}. \nonumber
\end{align}
Next, observe that
\begin{align*}
\ind{\abs{\mathbf{b}}-\abs{\mathbf{d}}}{\mathbf{b}}&=\min\{t\mid b_1+\cdots+b_t\geq b_1+\cdots+b_s-d_1-\cdots-d_l\} \\
&=\min\{t\mid b_{t+1}+\cdots+b_s\leq d_1+\cdots+d_l\}.
\end{align*}
We compute \eqref{eq:conj_uri_prop_form4} by considering the three cases $j>\ind{\abs{\mathbf{b}}-\abs{\mathbf{d}}}{\mathbf{b}}$, $j<\ind{\abs{\mathbf{b}}-\abs{\mathbf{d}}}{\mathbf{b}}$ and $j=\ind{\abs{\mathbf{b}}-\abs{\mathbf{d}}}{\mathbf{b}}$.

If $j>\ind{\abs{\mathbf{b}}-\abs{\mathbf{d}}}{\mathbf{b}}$, then we must have \[b_{j+1}+\cdots+b_s\leq d_1+\cdots+d_l,\quad b_j+\cdots+b_s\leq d_1+\cdots+d_l.\]
Therefore, both latter binomial coefficients in \eqref{eq:conj_uri_prop_form4} vanish, and thus the whole expression equals $0$.

If $j<\ind{\abs{\mathbf{b}}-\abs{\mathbf{d}}}{\mathbf{b}}$, then we must have
\[b_{j+1}+\cdots+b_s> d_1+\cdots+d_l,\quad b_j+\cdots+b_s> d_1+\cdots+d_l.\]
With \eqref{eq:conj_uri_prop_binomial_id} in the case $s=2$, we can reduce \eqref{eq:conj_uri_prop_form4} to
\begin{align*} 
&\binom{n-1}{(b_1+\cdots+b_{j-1})+(b_j+\cdots+b_s-d_1-\cdots-d_l)-1}\\
&\hspace{4cm}-\binom{n-1}{(b_1+\cdots+b_j)+(b_{j+1}+\cdots+b_s-d_1-\cdots-d_l)-1} \\
&=0.
\end{align*}
Finally, if $j=\ind{\abs{\mathbf{b}}-\abs{\mathbf{d}}}{\mathbf{b}}$, then we must have \[b_{j+1}+\cdots+b_s\leq d_1+\cdots+d_l,\quad b_j+\cdots+b_s> d_1+\cdots+d_l.\]
Therefore, the latter binomial coefficient in the second row of \eqref{eq:conj_uri_prop_form4} vanishes. So, with \eqref{eq:conj_uri_prop_binomial_id} the expression in \eqref{eq:conj_uri_prop_form4} reduces to
\[\binom{n-1}{(b_1+\cdots+b_{j-1})+(b_j+\cdots+b_s-d_1-\cdots-d_l)-1}=\binom{n-1}{\abs{\mathbf{b}}-\abs{\mathbf{d}}-1}.\]
\end{proof}

\makeinvisible{KEEP:
\begin{Proposition}\label{prop:s=1-case}
The conjecture \ref{conj:conj_uri_prop} holds for $s=1$, i.e. for all $b\geq 1$ we have
\begin{align*}
 \sum_{\mathbf{0} \leq \mathbf{r} \leq \mathbf{d}} (-1)^{|\mathbf{d}-\mathbf{r}|}
\binom{\mathbf{d}}{\mathbf{r}}
\binom{n+|\mathbf{r}|-1}{b-1}  =  \binom{n-1}{b - |\mathbf{d}|-1} .
\end{align*}
\end{Proposition}

\begin{proof} 
We consider the generating series over all $b\geq 1$ in the variable $T$. The right-hand side is given by 
\begin{align*}
    \sum_{b\geq 1} \binom{n-1}{b - |\mathbf{d}|-1} T^{b-1} = T^{|\mathbf{d}|} (1+T)^{n-1}.
\end{align*}
And the left-hand side is 
\begin{align*}
&\sum_{b\geq 1}  \sum_{\mathbf{0} \leq \mathbf{r} \leq \mathbf{d}} (-1)^{|\mathbf{d}-\mathbf{r}|}
\binom{\mathbf{d}}{\mathbf{r}}
\binom{n+|\mathbf{r}|-1}{b-1} T^{b-1} \\
&=       \sum_{\mathbf{0} \leq \mathbf{r} \leq \mathbf{d}} (-1)^{|\mathbf{d}-\mathbf{r}|}
\binom{\mathbf{d}}{\mathbf{r}}
\sum_{b\geq 1} \binom{n+|\mathbf{r}|-1}{b-1} T^{b-1} \\
&= \sum_{\mathbf{0} \leq \mathbf{r} \leq \mathbf{d}} (-1)^{|\mathbf{d}-\mathbf{r}|}
\binom{\mathbf{d}}{\mathbf{r}} (1+T)^{n+|\mathbf{r}|-1}\\
&= \prod_{j=1}^l \left( \sum_{r_j=0}^{d_j} (-1)^{d_j+r_j} \binom{d_j}{r_j} (1+T)^{r_j} \right) (1+T)^{n-1}
\\
&= \prod_{j=1}^l \left( T^{d_j} \right) (1+T)^{n-1} =T^{|\mathbf{d}|} (1+T)^{n-1},
\end{align*}
and therefore the statement is true for $s=1$.
\end{proof} 
}

As a special case of Theorem \ref{thm:conj_uri_prop}, we derive an identity used in Section \ref{sec:ari}.
\begin{Corollary} \label{cor:binomial_identities}
We have the binomial identity 
\[\sum_{r\in \ZZnonneg} (-1)^r \binom{d}{r}\binom{n+d-r-1}{b-1}=\binom{n-1}{b-d-1}, \qquad d\in \ZZnonneg,\ b,n\in\mathbb{N}.\]
\end{Corollary}  

\begin{proof} Set $s=l=1$ in Theorem \ref{thm:conj_uri_prop}. Then, the formula follows from the substitution $r\mapsto d-r$ and the observation that $\binom{d}{r}=0$ for $r>d$.
\end{proof}

\makeinvisible{KEEP: 
\begin{Proposition}  The conjecture \ref{conj:conj_uri_prop} holds for $l=1$ in average, i.e. we have for all $1\le j \le s$ the identity
\begin{align}\label{eq:average_composition_binomial_identity}
 \sum_{\mathbf{b} \in \mathbb{N}^s}    \sum_{0 \leq r \leq d} (-1)^{d-r }
    \binom{ d}{ r}
    \sum_{\substack{\mathbf{a} \in \mathcal{C}^s(n +  r ) \\ \ind{n}{\mathbf{a}}= j }} 
    \binom{\mathbf{a}-1}{\mathbf{b}-1} t^{|\mathbf{b}|-s}  = 
     \sum_{\substack{\mathbf{b} \in \mathbb{N}^s \\  \ind{|\mathbf{b}| -  d}{\mathbf{b}} =j  }}   \binom{n-1}{|\mathbf{b}| - d-1}  t ^{| \mathbf{b}|-s}.     
    \end{align}
 \end{Proposition}   
\begin{proof} We have 
\begin{align*} 
\sum_{\mathbf{b} \in \mathbb{N}^s}  \sum_{\substack{\mathbf{a} \in \mathcal{C}^s(n +  r ) \\ \ind{n}{\mathbf{a}}= j }} 
    \binom{\mathbf{a}-1}{\mathbf{b}-1} t^{|\mathbf{b}|-s}  &= 
\sum_{\mathbf{b}\in \mathbb{N}^s}  \sum_{\substack{\mathbf{a} \in \mathcal{C}^s(n +  r ) \\ \ind{n}{\mathbf{a}}= j }} 
\prod_{i=1}^s \binom{a_i-1}{b_i-1} t^{b_i-1} \\
&= \sum_{\substack{\mathbf{a} \in \mathcal{C}^s(n +  r ) \\ \ind{n}{\mathbf{a}}= j }}  ( 1+t)^{a_1-1} \cdots ( 1+t)^{a_s-1} \\
 &= \# \big\{  \mathbf{a} \in \mathcal{C}^s(n+r) \,| \,\ind{n}{\mathbf{a}}= j \,\big\} 
  (1+t)^{ n + r -s} \\
 &= \binom{n-1}{j-1} \binom{ r}{s-j} (1+t)^{ n + r -s} \\
 &= \binom{n-1}{j-1} \binom{ r}{s-j} \sum_{l} \binom{n+r-s}{l} t^{l}.
 \end{align*}  
This implies for the left hand side of \eqref{eq:average_composition_binomial_identity}  
 
\begin{align*}
lhs &= \sum_l \Big( \sum_{0 \le r \le d} (-1)^{d-r} \binom{d}{r} \binom{n-1}{j-1} \binom{ r}{s-j} \binom{n+r-s}{l}\Big)  t^{l}\\
&= \binom{n-1}{j-1}  \binom{d}{s-j} \sum_l \Big( \sum_{0 \le r \le d} (-1)^{d-r}   \binom{ d-s+j}{r-s+j} \binom{n+r-s}{l}\Big)  t^{l} \\
&= \binom{n-1}{j-1}  \binom{d}{s-j} \sum_l \Big( \sum_{0 \le r \le d} (-1)^{d-r}   \binom{ d-s+j}{d-r} \binom{n+r-s}{l}\Big)  t^{l}\\
&= \binom{n-1}{j-1}  \binom{d}{s-j} \sum_l  \binom{n-j}{l+s-d-j} t ^{l}.
\end{align*}
For the last step, we used the formula of Proposition \ref{prop:s=1-case} with $l=1$.
For the right hand side of \eqref{eq:average_composition_binomial_identity} we get
\begin{align*}
rhs&=\sum_{\substack{\mathbf{b} \in \mathcal{C}^s(l) \\  \ind{l - d}{\mathbf{b}} =j  }}
  \binom{n-1}{l - d-1}  t ^{l-s}   \\
  &= \sum_l 
  \# \big\{  \mathbf{b} \in \mathcal{C}^s(l) \,| \,  \ind{l - d}{\mathbf{b}} =j \,\big\} \,\binom{n-1}{l - d-1}  t ^{l-s}\\
 &= \sum_l  \binom{l-d-1}{j-1} \binom{d}{s-j} 
    \binom{n-1}{l - d-1}  t ^{l-s} \\
 &= \sum_l  \binom{l+s-d-1}{j-1} \binom{d}{s-j} 
    \binom{n-1}{l+s - d-1}  t ^{l}\\
  &=\binom{n-1}{j-1} \binom{d}{s-j} \sum_l  \binom{n-j}{l+s-d-j} t ^{l}. \qedhere
   \end{align*}
\end{proof}
}

\subsection{Identities for Bernoulli numbers
} \label{app:threshold_shuffle}
Recall that $\mathcal{C}$ is the set of all compositions. For $\boldsymbol{\alpha}=(\alpha_1,\ldots,\alpha_s)\in \mathcal{C}$, we denote by $\ell(\boldsymbol{\alpha})=s$ the length. We have a map 
\begin{align*}
 \mathbb{N} \times \mathcal{C} &\to  \mathbb{N}  \times \mathbb{N}\\
( a, \boldsymbol{\alpha})  &\mapsto   \big(\,\ell( \boldsymbol{\alpha}), \, \ind{a}{\boldsymbol{\alpha}} \,\big),
\end{align*} 
where the second component is the threshold indicator function from \eqref{eq:app_def_threshold}.
Given any formal power series 
\[ \mathbf{B}(x,y) = \sum_{m,n\geq0 } B (m,n) x^m y^n, \]
with the normalizing condition $B(1,1)=1$, we can define multiplicities by
\begin{align}\label{def:threshold_index} 
 \mu_{a, \boldsymbol{\alpha}}^B = B\left(\ell(\boldsymbol{\alpha}), \ind{a}{\boldsymbol{\alpha}}\right), \qquad a\in \mathbb{N}, \boldsymbol{\alpha}\in \mathcal{C}.  
\end{align}
Here, we set $\mu_{a, \boldsymbol{\alpha}}^B=0$ if $ \ind{a}{\boldsymbol{\alpha}}=0$.

We extend the multiplicity $\mu_{a,\_}^B$ linearly in the second argument, i.e., for a formal sum of indices $n_1 \boldsymbol{\alpha}_1 + n_2 \boldsymbol{\alpha}_2$ we set
\[   \mu_{a, n_1\boldsymbol{\alpha}_1+n_2\boldsymbol{\alpha}_2}^B = n_1\, \mu_{a,\boldsymbol{\alpha}_1}^B + n_2 \, \mu_{a,\boldsymbol{\alpha}_2}^B.\]
The concatenation of two indices $ \kvec,\lvec$ is  given by the index $( k_1,...,k_{\ell(\kvec)}, l_1,...,l_{\ell(\lvec)})$.
We denote it by $(\kvec,\lvec)$ and we extend this pairing bilinear.
We combine these convention with the shuffle product of indices, for example
\[
\mu_{2, (4, (3,2) \shuffle (2) )}^B =
\mu_{2, (4, 2 \,(3,2,2)+ (2,3,2))}^B = 2 \mu_{2,  (4,3,2,2)}^B +\mu_{2,   (4,2,3,2)}^B . \]
By a decomposition of an index
$\kvec$ we understand a pair of indices $\kvec_1,\kvec_2$ such that $\kvec= (\kvec_1,\kvec_2)$, where we allow either $\mathbf{k}_1$ or $\kvec_2$ to be the empty index, i.e.,
$\kvec= (\emptyset,\kvec) =(\kvec, \emptyset)$ is allowed in this context. 
Given an index $\kvec=(k_1,k_2,\ldots,k_d)$ its reversed index is given by $\overline{\kvec}= (-1)^d \, (k_d,\ldots,k_2,k_1)$.

\begin{Definition}\label{def:threshold_shuffle}
We say a multiplicity $\mu_{a,\boldsymbol{\alpha}}^B$ given by \eqref{def:threshold_index} satisfies the \emph{threshold shuffle identities} if the following equalities hold
 \begin{itemize}
\item[(i)]For  all $\bsigma \in \mathbb{N}^{\ell(\bsigma)}$,  $\btheta \in \mathbb{N}^{\ell(\btheta)}$ 
and all $d_1,d_2 \in \ZZnonneg$  we have 
that  
\begin{align*}
  \mu_{t, \bsigma \shuffle \btheta}^B  
 &= \sum\limits_{\substack{  ( \bsigma_1, \bsigma_2)=\bsigma  \\   ( \btheta_1, \btheta_2)=\btheta }} 
\mu_{\abs{ \bsigma_2}  -d_1, \bsigma_2  }^B \, \mu_{t,( \bsigma_1 \shuffle  \btheta_1   ,\abs{ \bsigma_2}  -d_1 , \btheta_2 ) }^B   +
\mu_{\abs{ \btheta_2}  -d_2, \btheta_2  }^B \, \mu_{t,( \bsigma_1 \shuffle  \btheta_1   ,\abs{ \btheta_2 } -d_2 , \bsigma_2 ) }^B,
\end{align*}
where  $
t= \abs{\bsigma} +\abs{ \btheta} -d_1-d_2$,
\item[(ii)]
For all $\bsigma \in \mathbb{N}^{\ell(\bsigma)}$, $\btau \in \mathbb{N}^{\ell(\btau)}$, $\btheta \in \mathbb{N}^{\ell(\btheta)}$ 
and all $d_1,d_2 \in \ZZnonneg$   we have 
that
\begin{align*}
0= \sum\limits_{\substack{  ( \bsigma_1, \bsigma_2)=\bsigma  \\  (\btau_1, \btau_2)=\btau  \\ ( \btheta_1, \btheta_2)=\btheta}} &
\mu_{\abs{ \bsigma_2}+\abs{\btau_1 }  -d_1, \bsigma_2 \shuffle \overline{\btau_1} }^B \, 
\mu_{t,( \bsigma_1 \shuffle \overline{ \btheta_2}   ,\abs{ \bsigma_2}+\abs{\btau_1 }-d_1 , \btau_2 \shuffle \overline{\btheta_1} ) }^B
\\
&\hspace{3cm}+\mu_{\abs{ \btau_2 }+\abs{ \btheta_1}  -d_2, \btau_2 \shuffle \overline{\btheta_1}  }^B  \,
\mu_{t,( \bsigma_1 \shuffle \overline{\btheta_2}   ,\abs{ \btau_2 }+\abs{ \btheta_1}  -d_2 , \bsigma_2 \shuffle \overline{\btau_1}) }^B ,
\end{align*}
where $t = \abs{ \bsigma} +\abs{ \btau}+\abs{ \btheta} -d_1-d_2$.
\end{itemize}
\end{Definition}

Observe that   the identities in (i) and (ii)  are homogeneous with respect to scaling of each $B(\ell(\boldsymbol{\alpha}),\ind{a}{\boldsymbol{\alpha}})$ by a factor
$t^{ \ell(\boldsymbol{\alpha}) -1}$.

\begin{Example} (i)
For $(d_1,d_2,\bsigma,\btheta)=\big( 2,4,[1,2],[4,3] \big)$, we get with PARI/GP   the relation
\[
B(4,3) + 2 B(4,2) + 3 B(4,1) =    \big( 2B(3,1) +3  B(3,2) + B(3,3)\big) B(2,1).
\]

\makeinvisible{With $(d_1,d_2,\bsigma,\btheta)=\big( 2,4,[1,2],[4,3] \big)$, we have $t= 1+2 +4+3 -2-4 = 4$, and  
\begin{align*}
\mu_{4,[1,2]\shuffle[4,3] }^B &= \mu_{4, [1, 2, 4, 3] }^B + 
\mu_{4, [1, 4, 2, 3]}^B+ \mu_{4, [1, 4, 3, 2]}^B+ \mu_{4, [4, 1, 2, 3]}^B+\mu_{4, [4, 1, 3, 2]}^B + \mu_{4, [4, 3, 1, 2] }^B \\
 &=   B(4,3)  + 2 B(4,2) +  3 B(4,1). 
 \end{align*}
The possible decompositions of $\bsigma,\btheta$ are
\begin{align*}
(\bsigma_1, \bsigma_2)& \in \big\{ (  \emptyset, [1, 2]), ([1], [2]) ,([4, 3],  \emptyset ) \big\},\\ 
( \btheta_1, \btheta_2) & \in \big\{ (    \emptyset, [4, 3]),  ([4], [3]), ([4, 3],  \emptyset) \big\}.
\end{align*}
We now determine all non zero summands individually: 

For the pair $( \emptyset, [1, 2]), (    \emptyset, [4, 3]) $, we get
\[
\mu_{1,[1,2]}^B \mu_{4,[1,4,3]}^B + \mu_{3,[4,3]}^B \mu_{4, [3, 1, 2]}^B = B(2,1) B(3,2) + B(2,1)B(3,2).
\]
%
For the pair $(  \emptyset, [1, 2]), (   [4], [3]) $, we get
\[
 \mu_{1, [1, 2]}^B \mu_{4, [4, 1, 3]}^B = B(2,1) B(3,1).
 \]

For the pair $(  \emptyset, [1, 2]), (    [4, 3], \emptyset) $, we get
\[
 \mu_{1, [1, 2]}^B \mu_{4, [4, 3, 1]}^B = B(2,1) B(3,1).
 \]

For the pair $(  [1], [2]), (    \emptyset, [4, 3]) $, we get
\[
 \mu_{3, [4, 3]}^B  \mu_{4,   [1, 3, 2]}^B = B(2,1) B(3,2).
 \]


For the pair $(  [1, 2],\emptyset  ), (    \emptyset, [4, 3]) $, we get
\[
 \mu_{3, [4, 3]}^B  \mu_{4,   [1,2,3]}^B = B(2,1) B(3,3).
 \]

Finally, we get the relation
\[
B(4,3) + 2 B(4,2) + 3 B(4,1) =    \big( 2B(3,1) +3  B(3,2) + B(3,3)\big) B(2,1).
\]
}
(ii)  For $(d_1,d_2,\bsigma,\btau,\btheta)=(1,2,[2,1],[4,3],[5,6])$, we get with PARI/GP
\[
0= 6 B(4,4) B(3,2) + 9 B(3,3) B(4,3) + 21 B(4,4) B(3,3) - 15 B(6,6) .
\]


%
\end{Example}

Recall that the uri multiplicities from Definition \ref{def:uri_mult} are given by
\begin{align*}
\mu_{a,\boldsymbol{\alpha}}=B_1(\ell(\boldsymbol{\alpha}),\ind{a}{\boldsymbol{\alpha}}),\qquad a\in\mathbb{N},\boldsymbol{\alpha}\in\mathcal{C}, 
\end{align*}
where 
\[ B_1(m,n)=\frac{1}{m!} \sum_{k=0}^{n-1} \binom{m}{k} B_ k.\]
\begin{Conjecture} \label{conj:uri_mult_threshold_shuffle} 
The uri multiplicities from Definition \ref{def:uri_mult} satisfy the threshold shuffle identities from Definition \ref{def:threshold_shuffle}.
\end{Conjecture}

Observe if we write down the threshold shuffle identities for the uri multiplicities, we would get a family of identities for Bernoulli numbers. Therefore we may rephrase Conjecture \ref{conj:uri_mult_threshold_shuffle}  by saying the Bernoulli numbers satisfy the threshold shuffle identities.

\begin{Remark} The examples we calculated with PARI/GP indicate that the uri multiplicities are in fact the only non trivial multiplicities as in \eqref{def:threshold_index} satisfying the threshold shuffle identities up to scaling. More precisely, for $t \in \mathbb{R} \setminus \{ 0\}$ set
\begin{align*}
\mathbf{B}_t(x,y)= \frac{y^2 x e^{x t}}{(1-y)(e^{x y t}-1)}= \sum_{m,n\geq0} B_t (m,n) x^m y^n,
\end{align*}
and define a family of multiplicities as in \eqref{def:threshold_index} by
\[
\mu^t_{a, \boldsymbol{\alpha}} = B_t\left(\ell(\boldsymbol{\alpha}), \ind{a}{\boldsymbol{\alpha}}\right),\qquad a\in\mathbb{N},\boldsymbol{\alpha}\in\mathcal{C}.\]
Of course, the case $t=1$ gives the uri multiplicities. We conjecture the following holds:
\begin{itemize}
\item[(i)]
For all $t\neq 0$ the multiplicities $\mu^t_{a, \boldsymbol{\alpha}}$ satisfy the threshold shuffle identities from Definition \ref{def:threshold_shuffle}.   
\item[(ii)] If a multiplicity $\mu_{a, \boldsymbol{\alpha}}^{\tilde B}$ given by \eqref{def:threshold_index} satisfies the threshold shuffle identities from Definition \ref{def:threshold_shuffle}, then  
$ \mu_{a,\boldsymbol{\alpha}}^{\tilde B} = \mu_{a,\boldsymbol{\alpha}}^{B}$ 
  for all  $a\in\mathbb{N}$, $\boldsymbol{\alpha}\in\mathcal{C}$
and
\[
\mathbf{B}(x,y) =  
\begin{cases} 
xy & \text{ if } \tilde B(2,1) =0, \\
\mathbf{B}_{2 \tilde B(2,1)} (x,y) & \text{ if } \tilde B(2,1) \neq 0 . 
\end{cases}
\]
 \end{itemize}
 
\end{Remark}	

	\phantomsection
	
	\addcontentsline{toc}{section}{References}
	\bibliographystyle{amsalpha}
	\bibliography{biblio}

\providecommand{\bysame}{\leavevmode\hbox to3em{\hrulefill}\thinspace}
\providecommand{\MR}{\relax\ifhmode\unskip\space\fi MR }
\providecommand{\MRhref}[2]{%
  \href{http://www.ams.org/mathscinet-getitem?mr=#1}{#2}
}
\providecommand{\href}[2]{#2}
\begin{thebibliography}{AKEFM22}

\bibitem[AKEFM22]{AKEFM22}
M.~J.~H. Al-Kaabi, K.~Ebrahimi-Fard, and D.~Manchon, \emph{Post-{L}ie {M}agnus expansion and {BCH}-recursion}, SIGMA Symmetry Integrability Geom. Methods Appl. \textbf{18} (2022), Paper No. 023, 16.

\bibitem[And21]{andrews:compositions}
G.~E. Andrews, \emph{Compositions and {Chebyshev} polynomials}, From operator theory to orthogonal polynomials, combinatorics, and number theory. A~volume in honor of Lance Littlejohn's 70th birthday, Cham: Birkh{\"a}user, 2021, pp.~1--14 (English).

\bibitem[BB23]{BaBu_cMES}
Henrik Bachmann and Annika Burmester, \emph{Combinatorial multiple {Eisenstein} series}, Res. Math. Sci. \textbf{10} (2023), no.~3, 32 (English), Id/No 35.

\bibitem[BCK24]{AGZT}
A.~Burmester, N.~Confurius, and U.~K\"uhn, \emph{{AGZT}-lectures on formal multiple zeta values}, 2024, \href{https://arxiv.org/abs/2406.13630}{arxiv:2406.13630}.

\bibitem[BGF22]{BGF}
J.~I. Burgos-Gil and J.~{Fr}es{\'a}n, \emph{Multiple zeta values: from numbers to motives}, (version September 2022), \url{http://javier.fresan.perso.math.cnrs.fr/mzv.pdf}.

\bibitem[BGST24]{BGST24}
C.~Bai, L.~Guo, Y.~Sheng, and R.~Tang, \emph{Post-groups, ({L}ie-){B}utcher groups and the {Y}ang-{B}axter equation}, Math. Ann. \textbf{388} (2024), no.~3, 3127--3167. \MR{4705761}

\bibitem[BK20]{BK20}
Henrik Bachmann and Ulf K{\"u}hn, \emph{A dimension conjecture for {{\(q\)}}-analogues of multiple zeta values}, Periods in quantum field theory and arithmetic. Based on the presentations at the research trimester on multiple zeta values, multiple polylogarithms, and quantum field theory, ICMAT 2014, Madrid, Spain, September 15--19, 2014, Cham: Springer, 2020, pp.~237--258 (English).

\bibitem[Bro12]{br}
F.~Brown, \emph{Mixed {T}ate motives over $\mathbb{Z}$}, Annals of Mathematics \textbf{175} (2012), 949--976.

\bibitem[Bro21]{brown-depth}
\bysame, \emph{Depth-graded motivic multiple zeta values}, Compos. Math. \textbf{157} (2021), 529--572.

\bibitem[BT18]{BaTa}
Henrik Bachmann and Koji Tasaka, \emph{The double shuffle relations for multiple {Eisenstein} series}, Nagoya Math. J. \textbf{230} (2018), 180--212 (English).

\bibitem[Bur21]{Bu21}
D.~Burde, \emph{Crystallographic actions on {L}ie groups and post-{L}ie algebra structures}, Communications in Mathematics \textbf{29} (2021), no.~1, 67--89.

\bibitem[Bur23]{BPhD}
A.~Burmester, \emph{An algebraic approach to multiple q-zeta values}, PhD Thesis, Universit{\"a}t Hamburg (2023).

\bibitem[Bur24]{Bu_balanced}
\bysame, \emph{Balanced multiple $q$-zeta values}, Advances in Mathematics \textbf{349} (2024).

\bibitem[Bur25]{Bu_extended}
\bysame, \emph{An extended bigraded {L}ie algebra arising from multiple zeta values}, In preparation (2025).

\bibitem[BvIM24]{fMES}
Henrik Bachmann, Jan-Willem van Ittersum, and Nils Matthes, \emph{{F}ormal multiple {E}isenstein series and their derivations}, 2024, \href{https://arxiv.org/abs/2312.04124}{arxiv:2312.04124}.

\bibitem[Cat24]{catoire2024}
Pierre Catoire, \emph{The {C}artier-{Q}uillen-{M}ilnor-{M}oore theorem in the {P}ost-{H}opf case}, 2024, \href{https://arxiv.org/abs/2401.09116}{arxiv:2401.09116}.

\bibitem[Con]{Confurius:phd_thesis}
N.~Confurius, \emph{in preparation}, PhD Thesis, Universit{\"a}t Hamburg.

\bibitem[Eca02]{Ecalle02}
J.~Ecalle, \emph{A tale of three structures: the arithmetics of multizetas, the analysis of singularities, the {L}ie algebra {ARI}}, Differential equations and the {S}tokes phenomenon, World Sci. Publ., River Edge, NJ, 2002, pp.~89--146.

\bibitem[Eca11]{Ecalle11}
\bysame, \emph{The flexion structure and dimorphy: flexion units, singulators, generators and the enumeration of multizeta irreducibles}, Asymptotics in Dynamics, Geometry and PDEs; Generalized Borel Summation Volume II, Publications of the Scuola Normale Superiore Volume 12, 2011, pp.~27--211.

\bibitem[EFLMK15]{EF15}
K.~Ebrahimi-Fard, A.~Lundervold, and H.~Munthe-Kaas, \emph{On the {L}ie enveloping algebra of a {P}ost-{L}ie {A}lgebra}, Journal of Lie Theory \textbf{25} (2015), no.~4, 1139--1165.

\bibitem[EFM18]{EFM18}
K.~Ebrahimi-Fard and I.~Mencattini, \emph{Post-{L}ie algebras, factorization theorems and isospectral flows}, Discrete mechanics, geometric integration and {L}ie-{B}utcher series, Springer Proc. Math. Stat., vol. 267, Springer, Cham, 2018, pp.~231--285.

\bibitem[Foi17]{Foissy_post_lie}
L.~Foissy, \emph{Extension of the product of a post-{L}ie algebra and application to the {SISO} feedback transformation group}, 2017, \href{https://arxiv.org/abs/1705.08132v1}{arxiv:1705.08132}.

\bibitem[Gon05]{Gon05}
A.~B. Goncharov, \emph{Galois symmetries of fundamental groupoids and noncommutative geometry}, Duke Math. J. \textbf{128} (2005), no.~2, 209--284.

\bibitem[HI17]{HI17}
M.~E. Hoffman and K.~Ihara, \emph{Quasi-shuffle products revisited}, J. Algebra \textbf{481} (2017), 293--326.

\bibitem[JZ24]{jaza_postlie}
J.-D. Jacques and L.~Zambotti, \emph{Post-{L}ie algebras of derivations and regularity structures}, 2024, \href{https://arxiv.org/abs/2306.02484v2}{arxiv:2306.02484}.

\bibitem[K{\"u}h19]{K-montreal}
U.~K{\"u}hn, \emph{Lie-algebras associated to multiple q-zeta values}, Expansions, Lie Algebras, and Invariants. Workshop at CRM, Montreal, 2019, \href{https://www.math.utoronto.ca/~drorbn/Talks/CRM-1907/Kuehn_montreal_2019.pdf}{ www.math.utoronto.ca/$\sim$drorbn/Talks/CRM-1907/Kuehn\_montreal\_2019.pdf}.

\bibitem[Li24]{Li24}
Y.~Li, \emph{On the sub-adjacent {H}opf algebra of the universal enveloping algebra of a post-{L}ie algebra}, 2024, \href{https://arxiv.org/abs/2408.01345v1}{arxiv:2408.01345}.

\bibitem[LNT07]{LNT}
J.-G. Luque, J.-C. Novelli, and J.-Y. Thibon, \emph{Period polynomials and {I}hara brackets}, J. Lie Theory \textbf{17} (2007), 229–--239.

\bibitem[LST22]{lishta22}
Yunnan Li, Yunhe Sheng, and Rong Tang, \emph{{P}ost-{H}opf algebras, relative {R}ota-{B}axter operators and solutions of the {Y}ang-{B}axter equation}, 2022, \href{https://arxiv.org/abs/2203.12174}{arxiv:2203.12174}.

\bibitem[MKL13]{MKL13}
H.~Z. Munthe-Kaas and A.~Lundervold, \emph{On post-{L}ie algebras, {L}ie-{B}utcher series and moving frames}, Found. Comput. Math. \textbf{13} (2013), no.~4, 583--613.

\bibitem[Rac00]{Ra}
G.~Racinet, \emph{Séries génératrices non-commutatives de polyzêtas et associateurs de {D}rinfeld}, Thesis, Laboratoire Amiénois de Mathématique Fondamentale et Appliquée (2000).

\bibitem[Reu93]{Reut}
C.~Reutenauer, \emph{Free {L}ie algebras}, Oxford University Press (1993).

\bibitem[Sch25]{Sch15}
L.~Schneps, \emph{{ARI}, {GARI}, {Z}ig and {Z}ag: {A}n introduction to {E}calle's theory of multiple zeta values}, 2015/25, \href{https://arxiv.org/abs/1507.01534}{arxiv:1507.01534v3}.

\bibitem[Swe69]{Sweedler}
M.~E. Sweedler, \emph{Hopf algebras}, Mathematics Lecture Note Series, W. A. Benjamin, Inc., New York, 1969.

\bibitem[Val07]{Val07}
B.~Vallette, \emph{Homology of generalized partition posets}, J. Pure Appl. Algebra \textbf{208} (2007), no.~2, 699--725.

\end{thebibliography}

\end{document}